\documentclass[
 aip,
 cha,
 amsmath,amssymb,
 preprint,
 floatfix
]{revtex4-2}

\usepackage[utf8]{inputenc}
\usepackage[english]{babel}
\usepackage{graphicx}
\usepackage{parskip}
\usepackage[dvipsnames]{xcolor}
\usepackage{enumitem}
\usepackage{mathtools}
\usepackage{amsthm}
\usepackage{setspace}

\usepackage{algpseudocode}
\usepackage[normalem]{ulem}
\usepackage{tikz}
\usepackage{multirow}

\usepackage{todonotes}

\def\bs{\boldsymbol}

\theoremstyle{definition}

\usepackage{calrsfs}
\DeclareMathAlphabet{\pazocal}{OMS}{zplm}{m}{n}

\newcommand{\pA}{\pazocal{A}}
\newcommand{\pB}{\pazocal{B}}
\newcommand{\pC}{\pazocal{C}}
\newcommand{\pD}{\pazocal{D}}
\newcommand{\pE}{\pazocal{E}}
\newcommand{\pF}{\pazocal{F}}
\newcommand{\pG}{\pazocal{G}}

\newcommand{\pM}{\pazocal{M}}
\newcommand{\pN}{\pazocal{N}}
\newcommand{\pO}{\pazocal{O}}
\newcommand{\pP}{\pazocal{P}}
\newcommand{\pQ}{\pazocal{Q}}

\newcommand{\pS}{\pazocal{S}}

\newcommand{\bbE}{\mathbb{E}}

\newcommand{\bbR}{\mathbb{R}}
\newcommand{\bbP}{\mathbb{P}}

\newcommand{\eps}{\varepsilon}
\newcommand{\be}{{\bs e}}

\newcommand{\bx}{{\bs x}}

\newcommand{\bz}{{\bs z}}
\newcommand{\bW}{{\bs W}}

\newcommand{\ta}{\text{a}}
\newcommand{\tc}{\text{c}}
\newcommand{\tf}{\text{f}}
\newcommand{\ti}{\text{i}}

\newcommand{\tr}{\text{r}}
\newcommand{\tu}{\text{u}}

\usepackage{dsfont}

\usepackage[hidelinks]{hyperref}
\usepackage{xr-hyper}
\begin{document}

\pagenumbering{roman} 

\title{Quantifying combined rate- and noise-induced transitions in a multiscale system}

\author{Ryan Deeley}
\email[Author to whom correspondence should be addressed: ]{ryan.deeley@uol.de}
\affiliation{Theoretische Physik/Komplexe Systeme, ICBM, Carl von Ossietzky Universit{\"a}t Oldenburg, Oldenburg 26129, Deutschland}
\author{Peter Ashwin}
\affiliation{Department of Mathematics \& Statistics, University of Exeter, Exeter EX4 4QF, United Kingdom}
\author{Ulrike Feudel}
\affiliation{Theoretische Physik/Komplexe Systeme, ICBM, Carl von Ossietzky Universit{\"a}t Oldenburg, Oldenburg 26129, Deutschland}

\vspace*{-1cm}\date{\today}

\begin{abstract}
Excitations in multiscale dynamical systems produce persistent but impermanent excursions away from a quiescent state. A novel mechanism for inducing excitations, through a combination of rate-dependent parametric and stochastic forcing, is discussed using the FitzHugh-Nagumo model as an example. In the absence of noise, a \textit{quasithreshold} separating the state space between supercritical and subcritical behaviour is defined with respect to \textit{slow manifolds} associated with the metastable excited state. This is formulated in an autonomous system framework, which can be obtained via a \textit{compactification of time}. Techniques from singular perturbation theory elucidate how \textit{canard trajectories} shape the geometric structure of the quasithreshold as the system's \textit{pullback attractor} moves through it. The probability with which weakly stochastically driven sample paths travel through the quasithreshold and undertake (i) solely noise-induced and (ii) combined rate- and noise-induced excitations is estimated numerically for each of the corresponding forcing scenarios, using the \textit{adaptive multilevel splitting} algorithm. The probability of weak stochastic forcing preventing a rate-induced excitation, which occurs deterministically, is likewise estimated. This quantifies the different rates at which the probability of each transition scenario scales exponentially with the noise strength, in accordance with sample path large deviation principles. It is revealed that a combination of (deterministically subcritical) parametric and stochastic forcing is exponentially more likely to induce excitations than stochastic forcing alone, which is attributed to the pullback attractor and quasithreshold converging towards one another.
\end{abstract}

\maketitle

\newpage


\newpage

\pagenumbering{arabic}

{\bf 
Sudden and long-lasting shifts in the dynamics of a nonlinear system in response to slow changes in input are called critical transitions or tipping phenomena, and great advances have been made in modelling their onset using nonautonomous systems. In particular, critical thresholds can be found in parameter space when varying an external input faster than the system’s response rate (R-tipping), or in state space following successive destabilizing stochastic perturbations (N-tipping). Few studies investigate transitions driven by a mixture of these mechanisms, yet in reality, they act in conjunction, and their effects cannot be fully decoupled. Moreover, if the system itself has components evolving on multiple timescales, these present an extra complication. In this study, we consider a multiscale FitzHugh-Nagumo system and show how noise introduces uncertainty to rate-dependent transition scenarios. The system is considered in a monostable (excitable) setup where transitions correspond to long-lived but impermanent excursions away from equilibrium. We describe a quasithreshold separating \textit{subcritical} paths from \textit{supercritical} excursions using concepts from singular perturbation theory and by introducing a concept of $L$-metastability: we say a set is $L$-metastable if it loses mass at a linear rate slower than $L>0$. We study combined rate- and noise-driven sample paths with supercritical excursions for rates below the deterministic critical threshold associated with R-tipping, and also investigate cases where noise prevents sample paths from undertaking supercritical excursions. Using a rare-event algorithm, we numerically approximate the scaling of the probability for each of these scenarios in the weak-noise limit. This quantification reveals that the combination of rate-dependent parametric and stochastic forcing drastically increases the probability of the system undertaking a supercritical excursion, in comparison to pure stochastic forcing. 
}

\section{Motivation}

Critical transitions (particularly ``tipping phenomena'') are of considerable interest in complexity science as they account for abrupt changes which can occur in the qualitative configuration of nonlinear systems subjected to time-dependent (or \textit{nonautonomous}) forcing that develops gradually and displays no obvious cause for driving the systematic change. Conceptually, these abrupt qualitative changes may refer either to excursions in the system's \textit{state space} or transitions to distinct quasi-stationary dynamical regimes. Examples of critical transitions have been studied in a wide variety of fields, such as climate science, neuroscience, ecology, economics, and physics; see, for example, \cite{scheffer2009critical,kuehn2011,ashwin2012}. There the shifts in the dynamics of the corresponding systems refer, for instance, to a breakdown of the Atlantic Meridional Overturning Circulation leading to a dramatically changed climate in Northwest Europe\cite{rahmstorf1995bifurcations,drijfhout2025shutdown,van2025european,van2025changing,van2025physics}, the melting of the Greenland or West Antarctic Ice Sheet leading to substantial sea level rise and freshwater flux~\cite{goelzer2020future,box2022greenland}, and the collapse of plant-pollinator networks~\cite{lever2014sudden,bascompte2023resilience} leading to a loss of biodiversity, to name only a few. 

We outline three complementary directions of study for understanding critical transitions in nonautonomous systems: (1) classifying dynamical mechanisms leading to transitions~\cite{ashwin2012,ashwin2017parameter,wieczorek2023,feudel2023}; (2) developing and examining early-warning signals of a transition~\cite{kuehn2011,dakos2012methods,bernuzzi2025,benyami2023,chapman2025}; and (3) quantifying the likelihood of noise-induced transitions using large deviation theory~\cite{freidlinwentzell}. \citet{ashwin2012} noted that critical transition phenomena can arise through independent mechanisms: slow passage through a bifurcation (B-tipping); noise-induced escape from a basin of attraction (N-tipping); or rate-induced breakdown of an adiabatic approximation for slow forcing (R-tipping). Each of these in isolation can induce either notable excursions in state space or transitions to alternative coexisting stable states\cite{izhikevich2000neural,khovanov2013noise,wieczorek2011}, though in practical situations multiple different mechanisms may be present to some extent. This raises the question of understanding and quantifying how the various mechanisms interact\cite{lindner2004moments,shi2016towards}.

Note that R-tipping can be irreversible or reversible, as discussed in detail in \citet{wieczorek2023}. Irreversible R-tipping - which requires some notion of bistability\cite{feudel2008complex} - describes the case where a system fails to track a moving stable state as a result of being ``overswept'' by the (quickly) moving realm of attraction of an alternative stable state, which it then moves towards and remains there indefinitely~\cite{ashwin2017parameter}. On the other hand, reversible R-tipping - which can occur in monostable systems\cite{vanselow24} - describes the case where a system fails to continuously \textit{closely} track a moving stable state by temporarily undertaking a notable excursion in state space away from the moving stable states, before eventually regaining close proximity and remaining there indefinitely. These excursions can be related to two branches of dynamical systems theory: excitable systems and metastable states, both of which are closely connected to the occurrence of multiple timescales in the system. Here excitability refers to a dynamical phenomenon which manifests itself as a relatively fast change in an excitatory variable, called an \textit{activator}, which itself initiates a slow change in an \textit{inhibitory} variable that eventually suppresses the excitation and leads the system to return to its original state. This excitation can be understood as a notable excursion in state space. There are various well-understood mechanisms leading to such excitations (for details cf. \cite{izhikevich2000neural,lindner2004effects,izhikevich2007dynamical}), which play a prominent role in chemical reactions~\cite{scott2000complex}, neuroscience~\cite{izhikevich2007dynamical}, climate science~\cite{ganopolski2001rapid,hobden2025regularization} and ecology~\cite{truscott1994ocean,truscott1994equilibria}. 

By contrast, metastability refers more generally to nonlinear systems possessing multiple coexisting dynamical regimes that are characterised by different dynamical properties and confined to certain disjoint regions in state space. Those regimes can be observed transiently and, as such, enable the overall dynamics to contain long time intervals in which the system dwells in one regime before eventually transitioning to another. Transition times are in general short compared to the dwelling times~\cite{froyland2013metastability}. The simplest example of systems with metastable (long-lived but transient) regimes are those exhibiting multistability, i.e., a multitude of coexisting attractors for a given set of parameters, subjected to noise. There, the system would then stay in a region of state space around a certain attractor for some time before it is kicked out by the noise~\cite{kraut2002multistability,kaneko2003chaotic} and moves toward another attractor. There are also many other mechanisms which lead to distinct metastable regimes in deterministic systems (for details cf. \cite{rossi25}). One important class of systems with distinct metastable regimes are multiscale systems exhibiting \textit{slow} dynamics in the neighbourhood of the critical manifold and \textit{fast} transitions between disconnected attracting branches of the critical manifold \cite{kuehn2015}. This is the class of systems considered in this work, in which the timescale separation between the different processes is explicitly given by a parameter $\eps$. 

For a wide class of stochastic systems, the probability of noise-induced escape from a given set within a given time satisfies \textit{large deviation principles} expressed in terms of a rate functional\cite{freidlinwentzell,dembozeitouni,berglundgentz2006}. The rate functional, termed the Freidlin-Wentzell action functional, serves as a measure for the amount of ``work'' required by the stochastic forcing to realise a particular path in the system. In accordance, the large deviation principles state that the probability of escape decays exponentially with the noise strength, as it approaches zero, at a rate determined by the infimum of the rate functional evaluated across all corresponding escape paths.      

The interaction of R-tipping and N-tipping has been the subject of several previous studies. For example, \citet{ritchie2016,ritchie2017} and \citet{slyman2023} examine classes of forced low-dimensional systems and quantify the effects of noise on a subcritical rate of parametric forcing. In particular, \citet{ritchie2017} analysed the Fokker-Planck equation\cite{risken1989fpe,chang1970practical} to develop approximate formulae for the probability of irreversibly deviating from the trajectory of the underlying deterministic nonautonomous system which remains close to (``tracks'') the moving stable equilibrium. Contrastingly, \citet{slyman2023} approached the problem from a large deviation theory perspective and considered the higher-dimensional Hamiltonian extension the stochastic nonautonomous system. There, they proved the existence of a heteroclinic orbit between the stable equilibrium and saddle point of the past and future limit systems of the underlying deterministic nonautonomous system, respectively. This represents the \textit{most probable} combined rate- and noise-induced transition path leading to irreversible tipping. Further, it was shown that the \textit{energetic cost/action} of this heteroclinic orbit was less than that of the most probable noise-induced transition path within the corresponding stochastic autonomous system. This is to say that, within this example system, the inclusion of the parametric forcing increases the likelihood of the system undergoing irreversible tipping. In another study, \citet{ritchie2016} investigated early warning signals for these combined rate- and noise-induced irreversible transitions and showed that these events were correlated with a delayed increase in the variance and autocorrelation of the state variable. The cause for the delay was explained in terms of the most probable transition path constructed there for finite noise (different from those considered in \citet{slyman2023}).  

For multiple timescale (multiscale) systems, several interacting dynamical processes evolve on a variety of timescales.\cite{kuehn2015} These are relevant for modelling a wide range of complex systems, but tipping phenomena can be complicated by the presence of more than one timescale, as several studies have shown \cite{wieczorek2011,bastiaansen2023climate,okeeffe2020,osullivan2023,vanselow24,borner2024saddle}. 
There is relatively limited understanding of the interaction of R-tipping and N-tipping in multiscale systems, and that is the motivation for this study. More precisely, we aim to assess the influence of multiple timescales, parametric forcing, and noise on the likelihood of critical transitions within a conceptual excitable model subjected to a combination of stochastic and deterministic drivers. The critical transitions themselves are characterised with respect to a \textit{metastability landscape} that we propose for the system. 

The concept of $L$-metastability and a related metastability landscape is introduced in Section~\ref{sec:metastability}, and outlined as a tool to classify transients in monostable multiscale systems. In Section~\ref{sec:fitzhughnagumo_autonomous}, we introduce an autonomous FitzHugh-Nagumo multiscale system with a fixed input current $I$, and discuss its multiple timescale structure and related metastability landscape. This is then extended in Section~\ref{sec:fitzhughnagumo_nonautonomous} to a nonautonomous case for forcing by a current $I(t)$ that is asymptotically constant. We use the compactification method of \citet{wieczorek2021compactification} to transform the nonautonomous system into an autonomous system and discuss its metastability landscape. In Section~\ref{sec:stochastic_FHN_autonomous} we investigate the vulnerability to noise-induced excitations in a stochastic extension of the unforced FitzHugh-Nagumo system. This is quantified as the rate of decay of probability of occurrence in the weak-noise limit, which is framed by a sample path large deviation principle\cite{freidlinwentzell}. The probability of noise-induced excitations is evaluated numerically using an adaptive multilevel splitting algorithm\cite{cerou2007}, with details discussed in Section \begin{NoHyper}\ref{sm:AMS}\end{NoHyper} of the supplementary material\cite{supplementary_material}. Then, in Section~\ref{sec:stochastic_FHN_nonautonomous}, we consider a comparable stochastic extension of the compactified forced FitzHugh-Nagumo system, and assess the likelihood with which the combined rate- and noise-driven sample paths transition to the excited state for rates of parametric forcing below the critical threshold associated with R-tipping. Conversely, we assess the likelihood that sample paths avoid the rate-induced transition to the excited state that would occur without noise. This is performed by numerically evaluating a committor function\cite{Vanden-Eijnden2006} defining the probability that a rate-dependent tracking or tipping scenario switches under the influence of noise, respectively. We finish with a discussion and conclusion in Section~\ref{sec:discuss}, which highlights how the vulnerability to critical transitions changes as the timescale separation strengthens and the forcing components work in concert.

\subsection{The metastability landscape of a deterministic dynamical system}
\label{sec:metastability}

\begin{figure}[ht]
    \centering
    \includegraphics[width=3.37in]{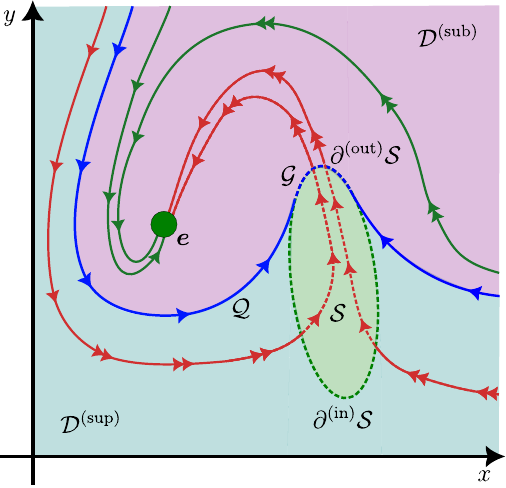}
    \caption{Schematic plot of a metastability landscape for a planar dynamical system. The green marker represents the globally asymptotically stable equilibrium point $\be$. The light green oval represents the $L$-metastable set $\pS$, with inflow boundary $\partial^\text{(in)}\pS$ and outflow boundary $\partial^\text{(out)}\pS$ corresponding to the dashed green and blue portion of $\partial\pS$, respectively. The supercritical set $\pD^{(\text{sup})}$ is indicated in light teal and the subcritical set $\pD^{(\text{sub})}$ is indicated in light purple. The quasithreshold $\pQ$ is shown in blue; the quasi-edge state $\pG$ is the dashed section of $\pQ$. The red trajectories indicate example supercritical excursions, where the slow motion through $\pazocal{S}$ corresponds to the dashed portion. The green trajectories illustrate solutions which remain subcritical.}
    \label{fig:metastability_landscape_schematic}
\end{figure}

Generally, we consider dynamical systems which may include stochastic as well as deterministic nonautonomous forcing. In this section, however, we introduce a framework for discussing certain aspects of metastability (a ``landscape'') in the context of deterministic autonomous dynamics. For a more general Perspective on metastability, see \cite{rossi25}. To explain the geometry of the metastability landscape, we conceptually introduce throughout the following certain \textit{domains} in state space, which are schematically illustrated in Figure \ref{fig:metastability_landscape_schematic}. Later, in Sections \ref{subsec:metastability_landscape_unforced} and \ref{subsec:metastability_landscape_forced}, we construct this landscape for a specific system, namely the FitzHugh-Nagumo model (in the absence and presence of parametric forcing, respectively).

Consider an autonomous dynamical system defined by the following ordinary differential equation initial value problem 
\begin{equation}\label{eqn:general_autonomous_ODE_IVP}
\frac{\text{d}\bx}{\text{d}t} = \bs b(\bx),\,\bx(0) = \bx^{(\ti)}\,,
\end{equation}
for $\bx(t)\in\bbR^k$ and smooth $\bs b:\bbR^k\rightarrow\bbR^k$. We write the solution of \eqref{eqn:general_autonomous_ODE_IVP} using the flow notation $\bs\phi_0(\bx^{(\ti)},t)\coloneqq \bx(t)$ and use $\ell()$ to denote  Lebesgue measure on $\bbR^k$. For any compact invariant set $\pA\subset \bbR^k$ let
$$
\pB(\pA)\coloneqq\{\bx^{(\ti)}\in \bbR^k:\omega(\bx^{(\ti)})\subseteq \pA\}
$$
denote the \textit{basin of attraction} of $\pA$, where 
\[
\omega(\bx^{(\ti)})\coloneqq\bigcap_{s\in\bbR}\text{cl}\big({\{\bs\phi_0(\bx^{(\ti)},t):t>s\}}\big)
\]
is the $\omega$-limit set\cite{glendinning94}. We assume there exists an equilibrium point $\be$ of \eqref{eqn:general_autonomous_ODE_IVP} which is asymptotically stable with basin $\pB(\pE)$, where $\pE=\{\be\}$.

For a bounded open subset $\pS$ of the phase space with smooth boundary $\partial \pS$ and outward unit normal $\hat{\bs n}()$, we define inflow and outflow boundaries as 
\begin{align*}    
\partial^\text{(in)}\pS & :=\{\bx\in \partial \pS~:~\bs\langle \bs b(\bx),\hat{\bs n}(\bx)\rangle<0\}\,, \\
\partial^\text{(out)}\pS & :=\{\bx\in \partial \pS~:~\bs\langle \bs b(\bx),\hat{\bs n}(\bx)\rangle>0\}\,,
\end{align*}
where $\langle\cdot\,,\cdot\rangle$ denotes the inner product. This generalises to piecewise smooth boundaries. If $\partial^\text{(out)}\pS=\emptyset$ then $\pS$ is forward-time invariant; if $\partial^\text{(in)}\pS=\emptyset$ then $\pS$ is backward-time invariant. In the case that $\partial^\text{(out)}\pS\neq\emptyset$, then $\pS$ may or may not contain an invariant set. More precisely, for each $\bx^{(\ti)}\in \pS$ we define the {\em residence time} of $\bx^{(\ti)}$ in $\pS$ as
$$
T^\text{(res)}(\bx^{(\ti)},\pS)\coloneqq\inf \{t>0~:~\bs\phi_0(\bx^{(\ti)},t)\not\in\pS\}
$$
and define the $T$-{\em persistent} subset
\[
\pS_T^\text{(pers)}\coloneqq\{\bx^{(\ti)}\in\pS~:~T^\text{(res)}(\bx^{(\ti)},\pS)> T\}\,,
\]
for every $T\in\bbR_+$. Note that $\pS_{T^{(1)}}^\text{(pers)}\subseteq \pS_{T^{(2)}}^\text{(pers)}$ if $T^{(1)}\geq T^{(2)}$ and so the measure $\ell(\pS_T^\text{(pers)})$ of the $T$-persistent subset will decrease monotonically with $T$. We say $\pS$ {\em loses mass at a rate slower than} $L\in\bbR_+$ if
\begin{align}
\ell(\pS_T^\text{(pers)})>\left(1-LT\right)\ell(\pS)~\mbox{ for every }~T\in(0,1/L)\,.
\label{eqn:losesmass}
\end{align}
If $\pS$ contains no invariant sets but loses mass at a rate slower than $L$, we say the set $\pS$ is $L$-{\em metastable}.  

We can use a metastable set to categorise other trajectories of the system. If a trajectory starting at $\bx^{(\ti)}$ passes through the $L$-metastable set $\pS$ in forward time, we say it undertakes a {\em supercritical} excursion, otherwise (when not infinitesimally close to a supercritical excursion) we consider it {\em subcritical}.
We can understand the coverage of supercritical excursions in the system by defining the set of initial conditions for \eqref{eqn:general_autonomous_ODE_IVP} that enter $\pS$ after finite positive time. This is the {\em supercritical set}
\begin{align}
\pD^{(\text{sup})}\coloneqq \{\bx^{(\ti)}\in\pB(\pE):\exists\, T\in\bbR_+\text{ with } \bs\phi_0(\bx^{(\ti)},T)\in \pS\}
\label{eqn:Dsup}
\end{align}
of initial conditions that induce supercritical excursions. Conversely, initial conditions in the {\em subcritical set}
\begin{align}
\pD^{(\text{sub})}&\coloneqq\{\bx^\text{(i)}\in\pB(\pE):\forall\,t\in\bbR_{\geq 0},\,\bs\phi_0(\bx^\text{(i)},t)\not\in\text{cl}(\pS)\} \label{eqn:Dsub}\\ 
&=\text{int}(\pB(\pE)\setminus\pD^\text{(sup)}) \notag 
\end{align}
give rise to solutions which remain subcritical in future time. In particular, the subcritical set comprises the initial conditions in the basin which: (i) do not visit $\pS$ in forward time under the flow of \eqref{eqn:general_autonomous_ODE_IVP}; and (ii) are away from the boundary of $\pD^\text{(sup)}$. Note that by construction $\pD^{(\text{sup})}$  is an open set which is backward-time invariant while $\pD^{(\text{sub})}$ is an open set which is forward-time invariant.

We define a {\em quasithreshold} $\pQ$ for the $L$-metastable set $\pS$ as the remaining initial conditions in the basin which are neither supercritical nor subcritical, 
\begin{align}\label{eqn:Q}
\pQ  & \coloneqq \pB(\pE)\,\setminus\,(\pD^\text{(sup)}\cup\pD^\text{(sub)}) \\ 
& = \{\bx^\text{(i)}\in\pB(\pE)\setminus\pD^\text{(sup)}:\exists\,T\geq 0\text{ with } \bs\psi(\bx^\text{(i)},T)\in\partial\pS\} \notag\,. 
\end{align}
The quasithreshold can be thought to comprise the ``critical solutions'' in the basin of attraction which, in forward time, visit $\text{cl}(\pS)$ but do not visit $\pS$. This behaviour is induced (exclusively) by initial conditions in the basin of attraction on the boundary of the supercritical set, i.e. 
\begin{align*}
\pQ&=\pB(\pE)\cap\partial\pD^\text{(sup)} \\ &=\pB(\pE)\cap\partial(\pB(\pE)\setminus\partial\pD^\text{(sup)})\,.    
\end{align*}
Moreover, if $\partial\,\text{int}(\pB(\pE)\setminus\pD^\text{(sup)}) = \partial (\pB(\pE)\setminus\pD^\text{(sup)})$, then $\pQ = \pB(\pE)\cap\partial\pD^\text{(sub)}$, implying that - within the basin of attraction - the quasithreshold, the boundary of the supercritical set, and the boundary of the subcritical set are all the same. 
The supercritical set, subcritical set, and quasithreshold give a disjoint partition of the basin of attraction, in the sense that
$$
\pB(\pE)=\pD^\text{(sup)}\cup\pD^\text{(sub)}\cup\pQ.
$$
We note some similarities and differences with regular thresholds that are on basin boundaries (cf. \citet{wieczorek2023}).
The quasithreshold $\pQ$ will generally be codimension one in phase space, like regular thresholds. However, unlike regular thresholds, the quasithreshold is generally not normally hyperbolic and repelling (globally) and not forward-time invariant (for instance, whenever $\text{cl}(\pS)\cap\pE=\emptyset$). In particular, the system need not have a regular edge state whose stable manifold can be used to characterise regular thresholds \cite{wieczorek2023}. 
Nonetheless, the subset of $\partial\pS$ which the trajectories along $\pQ$ intersect with in forward time is given by   
\begin{flalign*}    
\bigcup_{\bx^{(\ti)}\in\pQ}\big(\{\bs\phi_0(\bx^{(\ti)},t):t\geq 0\}\cap\partial\pS\big)&=\partial\pS\cap\bigcup_{\bx^{(\ti)}\in\pQ}\{\bs\phi_0(\bx^{(\ti)},t):t\geq 0\}\\ 
&=\partial\pS\cap\pQ \\ 
&\eqqcolon \pG
\end{flalign*}
and as such $\pG$ can be thought of as a \textit{quasi-edge state} for the quasithreshold in the basin of attraction.

We define a {\em metastability landscape} for an individual basin of attraction to be an $L$-metastable set for some $L\in\bbR_+$, the related supercritical and subcritical sets, and the quasithreshold. We use this to provide a more detailed description of transient dynamics within an individual basin of attraction. In particular, in cases where the system has several timescales, we can use it to describe subsets of the basin where many trajectories get ``stuck'' for a finite time with escape at rate $L$.  Note that if there are multiple attractors, each basin may or may not have regions of metastability.
For a discussion of the various types of dynamics associated with metastability more generally, see \citet{rossi25} and \citet{hancock24}. 

\subsection{The persistence of metastable regimes under weak stochastic forcing}\label{subsec:almost_invariant}

In this section we consider additional stochastic forcing on the conceptual dynamical system \eqref{eqn:general_autonomous_ODE_IVP}. In this context we wish to outline how, upon entering the $L$-metastable set $\pS$, the properties of becoming ``temporarily stuck there'' - as derived in the deterministic setting - similarly persist for stochastic sample paths under sufficiently weak noise. To this end, consider the dynamical system defined by the following stochastic differential equation initial value problem
\begin{equation}\label{eqn:general_autonomous_SDE_IVP}
    \text{d}\bx = \bs b(\bx)\,\text{d}t+\sigma\,\bs\Gamma(\bx)\,\text{d}\bW,\,\bs x(0)=\bx^\text{(i)}\,,
\end{equation}
for $\bx(t)\in\bbR^k$ and smooth $\bs\Gamma:\bbR^k\rightarrow\bbR^{k\times l}$. Here $\sigma\in\bbR_+$ describes the noise strength and $\bW$ is an $l$-dimensional Wiener process defined on a filtered probability space\cite{karatzas_shreve} with sample space $\bs\Omega$ and probability measure $\rho_{\bs\Omega}$. We write the solution of \eqref{eqn:general_autonomous_SDE_IVP} using the flow notation $\bs\phi_\sigma(\bx^\text{(i)},\bW,t)\coloneqq \bx(t)$. Note the dependence on the individual Wiener process $\bW$.

An intuitive result of sample path large deviation theory states that, in the weak-noise limit, sample paths of \eqref{eqn:general_autonomous_SDE_IVP} agree arbitrarily closely with the underlying deterministic trajectory of \eqref{eqn:general_autonomous_ODE_IVP} over a given finite time interval. Concretely, given any finite $T>0,\,\delta>0$, and $\bx^{(\ti)}\in\mathbb{R}^k$, it holds that 
\begin{equation}\label{eqn:asymptotic_agreement}
    \lim_{\sigma\downarrow 0}\,\bbP_{\bW\sim\rho_{\bs\Omega}}\bigg(\sup_{t\in[0,T]}\big\lVert\bs\phi_0(\bx^{(\ti)},t)-\bs\phi_\sigma(\bx^{(\ti)},\bW,t)\big\rVert<\delta\bigg)=1\,,
\end{equation}
see \citet{freidlinwentzell} (Chapter 2, Theorem 1.2). Now, for each $\bx^\text{(i)}\in\pS$, define the \textit{stochastic residence time} of $\bx^\text{(i)}$ driven by $\bW$ in $\pS$ as
\begin{equation*}
    T_\sigma^\text{(res)}(\bx^\text{(i)},\bW,\pS) \coloneqq \inf\{t>0:\bs\phi_\sigma(\bx^\text{(i)},\bW,t)\not\in\pS\}\,.
\end{equation*}
A particular consequence of the asymptotic property \eqref{eqn:asymptotic_agreement} is that
\begin{equation*}
    \lim_{\sigma\downarrow 0}\,\bbP_{\bW\sim\rho_{\bs\Omega}}\big(T_\sigma^\text{(res)}(\bx^\text{(i)},\bW,\pS)>T\big)=\begin{cases}
        1 & \text{if } \bx^{(\ti)}\in\{\bx\in\pS:T^\text{(res)}(\bx,\pS)>T\} \\ 0 & \text{if } \bx^{(\ti)}\in\{\bx\in\pS:T^\text{(res)}(\bx,\text{cl}(\pS))<T\} 
    \end{cases}\,.
\end{equation*}
Accordingly, the Lebesgue Dominated Convergence Theorem\cite{roydenfitzpatrick} implies that
\begin{flalign*}
    \liminf_{\sigma\downarrow 0}\,\int_{\pS}\bbP_{\bW\sim\rho_{\bs\Omega}}\big(T_\sigma^\text{(res)}(\bx^\text{(i)},\bW,\pS)>T\big)\,\text{d}\ell(\bx) & \geq \ell\big(\{\bx\in\pS:T^\text{(res)}(\bx,\pS)>T\}\big) \\  
    & = \ell(\pS_T^\text{(pers)}) \\
    & > (1-LT)\ell(\pS)\,,
\end{flalign*}
or equivalently,
\begin{flalign}\label{eqn:almost_invariance_analogy}
    \liminf_{\sigma\downarrow 0}\,\bbP\big(T_\sigma^\text{(res)}(\bx^\text{(i)},\bW,\pS)>T\big)& \geq \frac{\ell(\pS_T^\text{(pers)})}{\ell(\pS)} \notag\\ 
    & > (1-LT)\,,
\end{flalign}
where the probability is taken with respect to the measure $(\bW,\bx^\text{(i)})\sim\rho_{\bs\Omega}\times\text{Unif}(\pS)$.
Recall that these limits hold pointwise for each finite $T>0$. Still, this implies that with sufficiently weak noise strengths, $\pS$ loses mass in a probabilistic sense in system \eqref{eqn:general_autonomous_SDE_IVP} at a rate slower than one well approximated by $L$. It is, for example, a corollary of \eqref{eqn:almost_invariance_analogy} that
\[
\liminf_{\sigma\downarrow 0}\bigg(\inf_{t\in[0,1/L]}\big(\bbP(T_\sigma^\text{(res)}(\bx^\text{(i)},\bW,\pS)>t)-(1-Lt)\big)\bigg)\geq 0\,,\text{with }(\bW,\bx^\text{(i)})\sim\rho_{\bs\Omega}\times\text{Unif}(\pS)\,.
\]

The asymptotic result \eqref{eqn:almost_invariance_analogy} for $0<T\ll 1$ is reminiscent of what is meant for a set to be \textit{almost-invariant}\cite{froyland2005,froyland2009}. To define this concept in the context of stochastic differential equations, let $\bs\mu$ be an invariant measure of \eqref{eqn:general_autonomous_SDE_IVP}. The set $\pA\subseteq\bbR^k$ is said to be \textit{almost-invariant} on the time frame $T\in\bbR_+$ if 
\begin{equation*}
    \bbP\big(\bs\phi_\sigma(\bx^\text{(i)},\bW,T)\in\pA\big)\approx 1\,,\quad \text{with  }(\bW,\,\bs x^\text{(i)})\sim\rho_{\bs\Omega}\times\bs\mu\rvert_{\pA}\,,
\end{equation*}
which is analogous to the corresponding definition given for deterministic dynamical systems in \citet{froyland2009}. The symbol $\approx$ is taken to mean ``close to'', but not equal to. One could infer from \eqref{eqn:almost_invariance_analogy} that $\pS$ would be considered an \textit{almost-invariant} set for system \eqref{eqn:general_autonomous_SDE_IVP} on the timescale $0<T\ll 1/L$ when $\sigma>0$ is sufficiently small, provided that the change in the sampling density from Lebesgue measure to the invariant measure does not significantly affect the probability of remaining in $\pS$ over this time frame. This is a reasonable assumption whenever the invariant density is concentrated away from the boundary of $\pS$.

More generally, in the stochastic context, we say that the set $\pS$ \textit{loses mass in the stochastic invariant regime at a rate slower than} $L_\sigma\in\bbR_+$ if 
\begin{equation*}
    \bbP(T_\sigma^\text{(res)}(\bx^\text{(i)},\bW,\pS)>T)>1-L_\sigma T~\mbox{ for every }~T\in(0,1/L_\sigma)\,,
\end{equation*}
with $(\bW,\bx^\text{(i)})\sim\rho_{\bs\Omega}\times\bs\mu\lvert_\pS$.
If $\pS$ contains no invariant sets but loses mass in the invariant regime at a rate slower than $L_\sigma$, we say the set $\pS$ is $L_\sigma$-{\em metastable} for the stochastic system \eqref{eqn:general_autonomous_SDE_IVP} with noise strength $\sigma \in \bbR_+$.   

\section{An autonomous FitzHugh-Nagumo system}\label{sec:fitzhughnagumo_autonomous}

For our study of critical transitions in a multiscale setting, we consider the FitzHugh-Nagumo system. This was first derived in \textcite{FitzHugh1961} (and subsequently simulated electronically in \textcite{nagumo1962active}) as a simple model of a spiking neuron that includes excitable, bistable and oscillatory behaviour depending on parameters. See, for example, \textcite{izhikevich2007dynamical} for a review.
The FitzHugh-Nagumo system is an autonomous dynamical system defined by the following ordinary differential equation initial value problem
\begin{equation}\label{eqn:FHN_autonomous_slowtime}
    \begin{cases}
        \displaystyle\frac{\text{d}u}{\text{d}t} = \frac{1}{\eps}(-u^3 + u - v + I) \vspace{0.2cm}\\ 
        \displaystyle\frac{\text{d}v}{\text{d}t} = -\beta v + u
    \end{cases},\,\begin{pmatrix}
        u(0) \\ v(0)
    \end{pmatrix} = \begin{pmatrix}
        u^{(\ti)} \\ v^{(\ti)}
    \end{pmatrix}\,.
\end{equation}
The state variables $(u(t),v(t))\in\bbR^2$ are an \textit{activator-inhibitor} pairing, respectively describing evolving \textit{membrane potential} of neuronal cells and a variable controlling the recovery of the system following a stimulation. Here $I\in\bbR$ describes an external input current applied to the neuronal cell, $\beta\in\bbR$ is a parameter of the recovery procedure, and $0<\eps\ll 1$ describes the timescale separation between the neuronal spiking and recovery phases. For all of what follows, we fix $\beta = 0.75$.

We write the solution of the initial value problem (\ref{eqn:FHN_autonomous_slowtime}) using the flow notation
$$
\bs\phi_{0,\eps,I}^{(\text{FHN})}((u^{(\ti)},v^{(\ti)}),t):=(u(t),v(t))\,,\,\,\,t\in\mathbb{R}\,,
$$
and  the \textit{drift field} of \eqref{eqn:FHN_autonomous_slowtime} using the vectorised notation
\begin{equation}\label{eqn:autonomous_FHN_drift_field}
    \bs b_{\varepsilon,I}^{(\text{FHN})}(u,v) \coloneqq \begin{pmatrix}
        \varepsilon^{-1}(-u^3+u-v+I) \\ -\beta v + u
    \end{pmatrix}\,.
\end{equation}

\subsection{Fast-slow system analysis}\label{sec:fitzhughnagumo_autonomous_fastslow_analysis}

Upon considering the following equivalent formulation of \eqref{eqn:FHN_autonomous_slowtime},
\begin{equation*}
    \begin{cases}
        \eps\displaystyle\frac{\text{d} u}{\text{d}t} = -u^3 + u - v + I \vspace{0.2cm}\\ 
        \,\,\,\displaystyle\frac{\text{d} v}{\text{d}t} = -\beta v + u
    \end{cases},\,\begin{pmatrix}
        u(0) \\ v(0)
    \end{pmatrix} = \begin{pmatrix}
        u^{(\ti)} \\ v^{(\ti)}
    \end{pmatrix}
\end{equation*}
and then setting $\eps = 0$, we obtain the differential algebraic equation initial value problem
\begin{equation}\label{eqn:FHN_autonomous_slowsubsystem}
    \begin{cases} 
        0 = -u^3 + u - v + I\\
        \displaystyle\frac{\text{d} v}{\text{d}t} = -\beta v + u
    \end{cases},\,\begin{pmatrix}
        u(0) \\ v(0)
    \end{pmatrix} = \begin{pmatrix}
        u^{(\ti)} \\ v^{(\ti)}
    \end{pmatrix}\,.
\end{equation}
The resulting system $\eqref{eqn:FHN_autonomous_slowsubsystem}$ is referred to as the \textit{slow subsystem} of \eqref{eqn:FHN_autonomous_slowtime} since it provides a (first-order) approximation of the dynamics of \eqref{eqn:FHN_autonomous_slowtime} that occur on its \textit{slow} timescale. The dynamics of $\eqref{eqn:FHN_autonomous_slowsubsystem}$ are, by construction, constrained to the set 
\begin{equation}
\pS_{0,I}\coloneqq\{(u,v)\in\bbR^2:-u^3+u-v+I=0\},
\end{equation}
which is referred to as the \textit{critical manifold} of the system. Note that $\pS_{0,I}=\{(u,v)\in\mathbb{R}^2:v=h_I(u)\}$, where $h_I:\mathbb{R}\rightarrow\mathbb{R}$ defined via $h_I(u)=-u^3+u+I$ specifies the functional form of the critical manifold.

In the full system \eqref{eqn:FHN_autonomous_slowtime}, via an overshadowing evolution of the fast variable, the state is strongly attracted towards or repelled away from the critical manifold $\pS_{0,I}$ depending on its local stability properties. To arrive at the analogous (first-order) approximation of these dynamics which occur on the \textit{fast} timescale of the system, one has to consider the following rescaling of time $\tilde{u}(\tau) \coloneqq u(t(\tau)),\,\tilde{v}(\tau)\coloneqq v(t(\tau)),\,t(\tau) \coloneqq\eps\tau$ for the system $\eqref{eqn:FHN_autonomous_slowtime}$. This leads to the following rescaled initial value problem
\begin{equation}\label{eqn:FHN_autonomous_fasttime}
    \begin{cases}
        \displaystyle\frac{\text{d} \tilde{u}}{\text{d}\tau} = -\tilde{u}^3 + \tilde{u} - \tilde{v} + I \vspace{0.2cm}\\ 
        \displaystyle\frac{\text{d} \tilde{v}}{\text{d}\tau} = \eps(-\beta \tilde{v} + \tilde{u})
    \end{cases},\,\begin{pmatrix}
        \tilde{u}(0) \\ \tilde{v}(0)
    \end{pmatrix} = \begin{pmatrix}
        u^{(\ti)} \\ v^{(\ti)}
    \end{pmatrix}\,,
\end{equation}
for $\tau\in\bbR$ and $(\tilde{u}(\tau),\tilde{v}(\tau))\in\bbR^2$. Upon setting $\eps=0$ in \eqref{eqn:FHN_autonomous_fasttime}, one attains the following parameterised lower-dimensional differential equation initial value problem
\begin{equation}\label{eqn:FHN_autonomous_fastsubsystem}
    \begin{cases} 
        \,\,\, \displaystyle\frac{\text{d} \tilde{u}}{\text{d}\tau} = -\tilde{u}^3 + \tilde{u} - \tilde{v}+ I\\
        \tilde{v}(\tau)\equiv v^{(\ti)}
    \end{cases},\,\tilde{u}(0) = u^{(\ti)},
\end{equation}
known as the \textit{fast} or \textit{layer subsystem} of \eqref{eqn:FHN_autonomous_fasttime}.

The union of the set of equilibria of \eqref{eqn:FHN_autonomous_fastsubsystem}, across $v^\text{(i)}\in\bbR$, equals the critical manifold $\pS_{0,I}$. Further, given a subset $V$ of $\bbR$, a union of sets of attracting (resp. repelling) hyperbolic equilibria of \eqref{eqn:FHN_autonomous_fastsubsystem} across $v^\text{(i)}\in V$ is termed an \textit{attracting} (resp. \textit{repelling}) \textit{normally hyperbolic} subset of the critical manifold $\pS_{0,I}$. 

There are two examples of non-hyperbolic equilibria of \eqref{eqn:FHN_autonomous_fastsubsystem}, given by 
$$
\bs f_I^{(1)} := (u^{(\tf,1)},v_I^{(\tf,1)}):=(-\sqrt{3}/3,-2\sqrt{3}/9+I),~~\bs f_I^{(2)} := (u^{(\tf,2)},v_I^{(\tf,2)}):= (\sqrt{3}/3,2\sqrt{3}/9+I)\,.
$$
These correspond to the position of saddle-node bifurcations within $\eqref{eqn:FHN_autonomous_fastsubsystem}$, and are termed \textit{nondegenerate fold points} of the critical manifold $\pS_{0,I}$ (see Definition 3.3.2 in \citet{kuehn2015}). 

\begin{figure}[ht]
    \centering
    \includegraphics[width=3.37in]{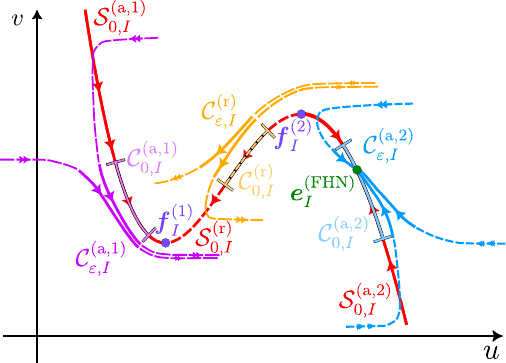}
    \caption{Schematic illustration of slow manifolds in the FitzHugh-Nagumo system \eqref{eqn:FHN_autonomous_slowtime}. The critical manifold $\pS_{0,I}$ contains normally attracting manifolds $\pS_{0,I}^{(\ta,1)},\,\pS_{0,I}^{(\ta,2)}$ [solid red curves], a normally repelling manifold $\pS_{0,I}^{(\tr)}$ [dashed red curve], and fold points $\bs f_I^{(1)},\,\bs f_I^{(2)}$ [filled purple markers]. The flow of the slow subsystem $\eqref{eqn:FHN_autonomous_slowsubsystem}$ on $\pS_{0,I}$ [red arrow heads] and the unique equilibrium point $\be_I^{(\text{FHN})}$ [filled green marker] is indicated. Compact submanifolds of $\pS_{0,I}^{(\ta,1)}$, $\pS_{0,I}^{(\tr)}$, and $\pS_{0,I}^{(\ta,2)}$ are exemplified by segments $\pC_{0,I}^{(\ta,1)}$ [light purple], $\pC_{0,I}^{(\tr)}$ [light orange], and $\pC_{0,I}^{(\ta,2)}$ [light blue], respectively. Trajectories of \eqref{eqn:FHN_autonomous_slowtime} with $0<\varepsilon\ll 1$ that closely shadow these submanifolds across their entire span are illustrated [purple, orange, and blue curves]. The solid portions $\pC_{\eps,I}^{(\ta,1)}$ [purple], $\pC_{\eps,I}^{(\tr)}$ [orange], and $\pC_{\eps,I}^{(\ta,2)}$ [blue] represent $\pO(\varepsilon)$ perturbations of $\pC_{0,I}^{(\ta,1)},\,\pC_{0,I}^{(\ta,1)}$ and $\pC_{0,I}^{(\ta,1)}$ respectively, and are examples of the ``slow manifolds'' corresponding to these submanifolds.}
    \label{fig:slow_manifolds_schematic}
\end{figure}

The union of the sets 
\begin{equation}
\pS_{0,I}^{(\text{a},1)}\coloneqq\{(u,v)\in\bbR^2:u<u^{(\tf,1)},v=h_I(u)\}
\end{equation}
and
\begin{equation}
\pS_{0,I}^{(\text{a},2)}\coloneqq\{(u,v)\in\bbR^2:u>u^{(\tf,2)},v=h_I(u)\}
\end{equation}
describe the largest attracting, normally hyperbolic subset of $\pS_{0,I}$. In contrast, the set
\begin{equation}
\pS_{0,I}^{(\text{r})}\coloneqq\{(u,v)\in\bbR^2:u^{(\tf,1)}<u<u^{(\tf,2)},v=h_I(u)\}
\end{equation}
describes the largest repelling, normally hyperbolic subset of $\pS_{0,I}$. The set $\pS_{0,I}^{(\text{r})}$ connects to the sets $\pS_{0,I}^{(\text{a},1)}$ and $\pS_{0,I}^{(\text{a},2)}$ via $\bs f_I^{(1)}$ and $\bs f_I^{(2)}$, respectively, and the critical manifold is fully comprised of the disjoint union
$$
\pS_{0,I} = \pS_{0,I}^{(\text{a},1)}\cup\{\bs f_I^{(1)}\}\cup\pS_{0,I}^{(\text{r})}\cup\{\bs f_I^{(2)}\}\cup\pS_{0,I}^{(\text{a},2)}.
$$
Consider the normally hyperbolic parts $\pS_{0,I}^{(\text{a},1)},\,\pS_{0,I}^{(\text{r})},\,\pS_{0,I}^{(\text{a},2)}$ of the critical manifold $\pS_{0,I}$. Fenichel theory (see for instance \cite{fenichel79,jones1995,kuehn2015}) implies that the flow of \eqref{eqn:FHN_autonomous_slowsubsystem} along compact subsets $\pC_{0,I}^{(\text{a},i)}$ of the attracting normally hyperbolic manifold $\pS_{0,I}^{(\text{a},i)}$ persists in \eqref{eqn:FHN_autonomous_slowtime} for $\varepsilon>0$ sufficiently small, along nearby, locally invariant\footnote{We say an embedded manifold $\pS$ with boundary $d \pS$ (taken with respect to the topology of the manifold) is {\em locally invariant}  \cite{kuehn2015} if for every $\bs x\in \pS\setminus d\pS$, there exists $\delta > 0$ whereby $\bs\phi(\bs x,t)\in\pS$ for all $t\in (-\delta,\delta)$, i.e. under the flow of the system you can only enter or leave the manifold $\pS$ at its boundary. This is equivalent to requiring that $\partial^\text{(in)}\pS\cup\partial^\text{(out)}\pS\cup\partial^\text{(gr)}\pS\subseteq d\pS$, for generalised inflow, outflow, and grazing boundaries
\begin{flalign*}
\partial^\text{(in)}\pS&\coloneqq\{\bs x\in\pS:\forall\,T>0,\,\exists\,(t_-,t_+)\in(-T,0)\times(0,T),\,\bs\phi(\bx,t_-)\not\in\pS,\,\bs\phi(\bx,t_+)\in\pS\}\,, \\ 
\partial^\text{(out)}\pS&\coloneqq\{\bs x\in\pS:\forall\,T>0,\,\exists\,(t_-,t_+)\in(-T,0)\times(0,T),\,\bs\phi(\bx,t_-)\in\pS,\,\bs\phi(\bs x,t_+)\not\in\pS\}\,, \\ 
\partial^\text{(gr)}\pS&\coloneqq\{\bs x\in\pS:\forall\,T>0,\,\exists\,(t_-,t_+)\in(-T,0)\times(0,T),\,\bs\phi(\bx,t_-)\not\in\pS,\,\bs\phi(\bs x,t_+)\not\in\pS\}\,.
\end{flalign*}\vspace*{-1cm}\textcolor{white}{.}}, attracting normally hyperbolic manifolds $\pC_{\eps,I}^{(\text{a},i)}$ ($i\in\{1,2\}$). Analogously, compact subsets $\pC_{0,I}^{(\tr)}$ of the repelling normally hyperbolic manifold $\pS_{0,I}^{(\text{r})}$ persist nearby as locally invariant, repelling normally hyperbolic manifolds $\pC_{\eps,I}^{(\text{r})}$. In this setup, the \textit{slow manifolds} $\pC_{\eps,I}^{(\cdot)}$ may be fully inflowing, or part inflowing and part overflowing. See Figure \ref{fig:slow_manifolds_schematic} for an illustration of the slow manifolds in the FitzHugh-Nagumo system \eqref{eqn:FHN_autonomous_slowtime}.

For any $I\in\bbR$, system \eqref{eqn:FHN_autonomous_slowtime} with $\beta = 3/4$ has a unique equilibrium given by 
\begin{align*}
\be_I^{(\text{FHN})}\coloneqq \big(u_I^{(\text{eq})},v_I^{(\text{eq})}\big)\coloneqq\bigg(\sqrt[3]{w_I^{(1)}}+\sqrt[3]{w_I^{(2)}},\bigg(\sqrt[3]{w_I^{(1)}}+\sqrt[3]{w_I^{(2)}}\bigg)/\beta\bigg),
\end{align*}
where $w_I^{(j)} = I/2 +(-1)^{j-1} \sqrt{I^2/4+1/27^2}$ for $j\in\{1,\,2\}$. For $I>\sqrt{3}(1/\beta-2/3)/3=2\sqrt{3}/9$ it holds that $u_I^{(\text{eq})}>u_I^{(\tf,2)}$; the equilibrium point therefore lies on the rightmost stable branch $\pS_{0,I}^\text{(a,2)}$ of the critical manifold in such cases. In addition, $\be_I^{(\text{FHN})}$ is a \textit{stable} equilibrium for the slow subsystem \eqref{eqn:FHN_autonomous_slowsubsystem}\footnote{In the context of a differential algebraic equation across a one-dimensional set, an equilibrium point is stable if the directional derivative of the drift field (in the direction of the tangent of the algebraic set) is negative at the equilibrium point.}, since $\text{d}v/\text{d}t = -\beta v + u < 0$ for $(u,v)\in\pS_{0,I}^{(\ta,2)}\big|_{v>v_I^{(\text{eq})}}$ and $\text{d}v/\text{d}t = -\beta v + u > 0$ for $(u,v)\in\pS_{0,I}^{(\ta,2)}\big|_{v<v_I^{(\text{eq})}}$ (see Figure \ref{fig:slow_manifolds_schematic}).

Altogether, the above theory implies that the equilibrium point $\be_I^{(\text{FHN})}$ is asymptotically stable for $\beta=3/4$, $I>2\sqrt{3}/9$ and $\eps>0$ sufficiently small. Accordingly, from this point onwards, we shall additionally restrict ourselves to $I>2\sqrt{3}/9$. Further, since there is no evidence of the existence of other compact invariant sets (in particular, limit cycles), we infer that $\be_I^{(\text{FHN})}$ is \textit{globally} asymptotically stable.

\subsection{A metastability landscape for unforced supercritical excursions}
\label{subsec:metastability_landscape_unforced}

\begin{figure}[ht]
    \centering
    \includegraphics[width=6.48in]{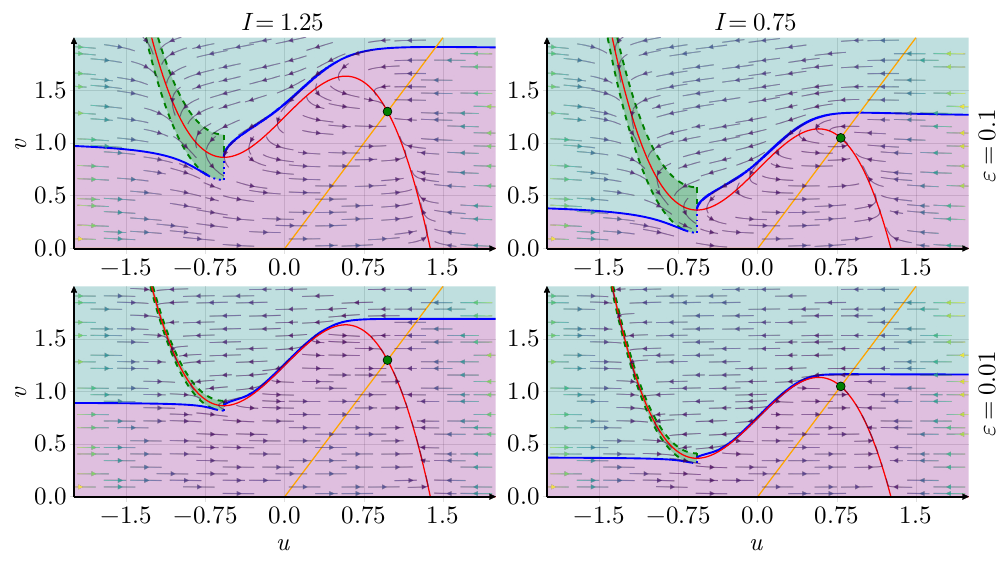}
    \caption{The metastability landscapes within the autonomous FitzHugh-Nagumo system for $(\eps,\,I)\in\{0.1,0.01\}\times\{1.25,0.75\}$. Each panel displays a phase portrait with the system drift [coloured flow lines - brighter colour denotes a greater magnitude of the derivative], equilibrium point $\be_I^{(\text{FHN})}$ [filled green marker], the $u$-nullcline (critical manifold) $\pS_{0,I}$ [red curve], the $v$-nullcline [orange curve], the $L$-metastable set $\tilde{\pS}_{\eps,I}$ [light green area], the supercritical set $\pD_{\eps,I}^{(\text{sup})}$ [light teal area], the subcritical set $\pD_{\eps,I}^{(\text{sub})}$ [light purple area], and the quasithreshold $\pQ_{\eps,I}$ [blue curve]. The inflow and outflow boundary of the $L$-metastable set are represented by the dashed green and dotted blue portions of $\partial\tilde{\pS}_{\eps,I}$, respectively. In particular, the geometry shows that $\pQ_{\eps,I}=\partial\pD_{\eps,I}^\text{(sup)}=\partial\pD_{\eps,I}^\text{(sub)}$. 
    } 
    \label{fig:autonomous_metastability_landscapes}
\end{figure}

We now interpret the multiple timescale dynamics of (\ref{eqn:FHN_autonomous_slowtime}) so as to relate the parts of the critical manifold defined above with the components of a metastability landscape discussed in Section~\ref{sec:metastability}, for different cases where the equilibrium $\be_I^{(\text{FHN})}\in \pS_{0,I}^{(\text{a},2)}$ is globally asymptotically stable. An illustration of the resulting metastability landscape, which is formulated algebraically throughout the following, is given in Figure \ref{fig:autonomous_metastability_landscapes}.   

To start with, consider a large but bounded absorbing (i.e., forward invariant) open region $\pM=\pM_u\times\pM_v\subset \bbR^2$ for $\eqref{eqn:FHN_autonomous_slowtime}$, and note that this is entirely contained within the basin $\pB:=\pB(\{\be_I^\text{(FHN)}\})=\bbR^2$. We base the metastability landscape around a neighbourhood of $\pS_{0,I}^{(\text{a},1)}$ that is large enough to contain all slow manifolds of \eqref{eqn:FHN_autonomous_slowtime} corresponding to compact subsets of $\pS_{0,I}^{(\text{a},1)}\,\cap\,\pM$, for all $\eps>0$ sufficiently small. This is given by the set 
\begin{equation}\label{eqn:L-metastable_set_unforced}
\tilde{\pS}_{\eps,\,I}\coloneqq\{(u,v)\in\pM:u<u_I^{(\tf,1)},\,\lvert v-h_I(u)\rvert<\eps^{2/3}\},
\end{equation}
which depends explicitly on the timescale separation $\varepsilon$ and input current $I$ (see Figure \ref{fig:autonomous_metastability_landscapes}). The particular structure of $\tilde{\pS}_{\eps,\,I}$ is motivated by the fact that, for planar fast-slow systems with a folded critical manifold, the delay in the slow variable beyond the fold point before the system transitions to the fast timescale scales like $\eps^{2/3}$ (see, i.e., \cite{krupa2001extending,kuehn2011}). Away from the fold point, the critical manifold is more strongly attracting (in the direction of the fast variable) and the slow manifolds deviate by a lesser extent.

For the system \eqref{eqn:FHN_autonomous_slowtime} we claim that $\tilde{\pS}_{\eps,I}$ will be an $L$-metastable set for $L$ of $\pO(1)$ as $\eps\downarrow 0$. This is because: (i) there is exponential attraction to a slow manifold $\pC_{\eps,I}^{(\text{a},1)}$ within $\tilde{\pS}_{\eps,I}$; (ii) the speed of motion on the slow manifold is at the slow timescale, i.e., of order 1 in $\eps$; and (iii) the length of $\tilde{\pS}_{\eps,I}$ is finite, positive, and independent of $\eps$. By contrast, an analogous neighbourhood of the repelling normally hyperbolic set $\pS_{0,I}^{(\text{r})}$ will be $L$-metastable only for $L$ of $\pO(1/\eps)$  as $\eps\downarrow 0$, because there is exponential repulsion away from a slow manifold $\pC_{\eps,I}^{(\text{r})}$ on the fast timescale, i.e., of order $1/\eps$ in $\eps$. 

As in Section~\ref{sec:metastability} we define the supercritical and subcritical subsets within the basin $\pB$ relative to this $L$-metastable set:
\begin{flalign*}
\pD_{\eps,I}^{(\text{sup})}&\coloneqq\{(u,v)\in\pB:\exists\, T\in\bbR_+\text{ with } \bs\phi_{0,\eps,I}^{(\text{FHN})}((u,v),T)\in \tilde{\pS}_{\eps,I}\}\,; \\[-1em]
\\
\pD_{\eps,I}^{(\text{sub})}&\coloneqq\{(u,v)\in\pB:\forall\, t\in\bbR_{\geq 0},\,\bs\phi_{0,\eps,I}^{(\text{FHN})}((u,v),t)\not\in \text{cl}(\tilde{\pS}_{\eps,I})\}\,.
\end{flalign*}
The corresponding quasithreshold is then defined as
\begin{flalign*}
\pQ_{\eps,I} &\coloneqq\{(u,v)\in\pB\setminus \pD_{\varepsilon,I}^\text{(sup)}:\exists\, T\in\bbR_{\geq 0}\text{ with } \bs\phi_{0,\eps,I}^{(\text{FHN})}((u,v),T)\in \text{cl}(\tilde{\pS}_{\eps,I})\}\,.
\end{flalign*}
For the FitzHugh-Nagumo system \eqref{eqn:FHN_autonomous_slowtime}, the supercritical set $\pD_{\eps,I}^{(\text{sup})}$ is constructed in terms of two particular system trajectories which separate the inflow and outflow boundaries of $\tilde{\pS}_{\eps,I}$, one of which is a normally hyperbolic and repelling \textit{maximal canard trajectory}\footnote{In this context, by ``maximal canard trajectory'', we mean a canard trajectory\cite{wieczorek2011} which shadows normally hyperbolic, repelling, slow manifolds for a maximal length of time.}. Figure \ref{fig:autonomous_metastability_landscapes} illustrates numerical solutions to these trajectories and the corresponding metastability landscapes of \eqref{eqn:FHN_autonomous_slowtime}, for a number of paired values of $\eps$ and $I$.     

\section{A nonautonomous FitzHugh-Nagumo system}\label{sec:fitzhughnagumo_nonautonomous}

In the case that $I$ varies with time, we can use the metastability landscape to help understand the response of the system to changes in $I$. In particular, under nonautonomous changes in $I$, the frozen-in quasithreshold $\pQ_{\eps,I}$ plays a role analogous to the regular threshold of \textcite{wieczorek2023}, but for a monostable system exhibiting metastability.

Figure~\ref{fig:autonomous_metastability_landscapes} demonstrates that $\be_{1.25}^{(\text{FHN})}\in\pD_{\eps,1.25}^{(\text{sub})}\setminus\pQ_{\eps,1.25}$ but $\be_{1.25}^{(\text{FHN})}\in\pD_{\eps,0.75}^{(\text{sup})}$, for each $\eps\in\{0.01,0.1\}$. This implies that if system \eqref{eqn:FHN_autonomous_slowtime} is at equilibrium with $I = 1.25$ and $0<\varepsilon\leq 0.1$, and the input current instantaneously decreases to $I = 0.75$, the state would proceed to undertake a supercritical excursion. More generally, if $I$ is varying smoothly in time and the system state fails to re-equilibrate at the rate at which the frozen-in quasithresholds $\pQ_{\eps,I}$ are shifting, the frozen-in quasithresholds can intersect with the position of the state and pass through it. The position of the state would then lie on the ``other side'' of the frozen-in quasithreshold; if at this point $I$ settles at its current value, the system would undertake a supercritical excursion. This describes the process behind rate-dependent transition events\cite{ashwin2012} which occur as a result of \textit{threshold instability}\cite{okeeffe2020,wieczorek2023}.              

To explore rate-dependent transition mechanisms enabled by quasithreshold instability in the FitzHugh-Nagumo system \eqref{eqn:FHN_autonomous_slowtime}, consider the nonautonomous dynamical system given by the following ordinary differential equation initial value problem
\begin{equation}\label{eqn:FHN_nonautonomous_slowtime}
    \begin{cases}
        \displaystyle\frac{\text{d} u}{\text{d}t} = \frac{1}{\eps}\big(-u^3 + u - v + J(rt)\big) \vspace{0.2cm}\\ 
        \displaystyle\frac{\text{d} v}{\text{d}t} = -\beta v + u
    \end{cases},\,\begin{pmatrix}
        u(t^{(\ti)}) \\ v(t^{(\ti)})
    \end{pmatrix} = \begin{pmatrix}
        u^{(\ti)} \\ v^{(\ti)}
    \end{pmatrix}
\end{equation}
for $(u(t),v(t))\in\bbR^2$. Here, $t^{(\ti)}\in\bbR$ describes an initial time at which the position of the state is known, and $r\in\bbR_+$ is a parameter determining the rate of evolution of the now time-dependent external input current $J:\bbR\rightarrow\bbR$. We write the solution of the nonautonomous system (\ref{eqn:FHN_nonautonomous_slowtime}) as a process
$$
\bs\Phi_{0,\eps,r}^{(\text{FHN})}((u^{(\ti)},v^{(\ti)}),t^{(\ti)},t) \coloneqq (u(t),v(t))
$$
and note that in general this depends on both the initial and final times, not just the elapsed time $t-t^{(\ti)}$.

For $J(t)$ we choose a smoothly varying monotonic function that is constant for $t$ outside a compact interval $[t_-,t_+]$, for $(t_-,t_+)\in\bbR^2$ with $t_-<t_+$,
\begin{flalign}\label{eqn:FHN_nonautonomous_ramping}
    J(t) &:= \begin{cases} I_- & \text{if } t\leq t_- \\  I_-+(I_+-I_-)\displaystyle\frac{\exp((t_--t)^{-1})}{\exp((t_--t)^{-1})+\exp((t-t_+)^{-1})} & \text{if } t_- < t < t_+ \\ I_+ & \text{if } t_+\leq t\end{cases}\,.
\end{flalign}
We assume $I_-,\,I_+>2\sqrt{3}/9$ with $I_-\neq I_+$. By construction of \eqref{eqn:FHN_nonautonomous_ramping}, if $\max(t,t^{(\ti)}) \leq t_-/r$, then
$$
\bs\Phi_{0,\eps,r}^{(\text{FHN})}((u^{(\ti)},v^{(\ti)}),t^{(\ti)},t)=\bs\phi_{0,\eps,I_-}^{(\text{FHN})}((u^{(\ti)},v^{(\ti)}),t-t^{(\ti)})\,.
$$
Similarly, if $\min(t,t^{(\ti)}) \geq t_+/r$, then 
$$
\bs\Phi_{0,\eps,r}^{(\text{FHN})}((u^{(\ti)},v^{(\ti)}),t^{(\ti)},t)=\bs\phi_{0,\eps,I_+}^{(\text{FHN})}((u^{(\ti)},v^{(\ti)}),t-t^{(\ti)})\,.
$$
For $t_-/r<t<t_+/r$ we must solve the nonautonomous system \eqref{eqn:FHN_nonautonomous_slowtime}. Hence, one can interpret system \eqref{eqn:FHN_nonautonomous_slowtime} as a concatenation of two distinct autonomous FitzHugh-Nagumo systems smoothly bridged together by a finite interval of explicitly nonautonomous dynamics lasting $(t_+-t_-)/r$ time units.       

\subsection{Compactification of time}

We recall the procedure of compactifying the time parameterisation \cite{wieczorek2021compactification} of the external input current forcing across $t \in (-\infty,\infty)$ into an equivalent parameterisation across the compact interval $s\in[-1,1]$, if $I(t)$ is asymptotically constant. This converts the system with explicit time-dependent drift into an autonomous form, with a compact absorbing portion of state space. It comes at the cost of increasing the dimensionality of the system by one degree, but allows one to cleanly formulate statements about rate-dependent transitions in terms of unstable manifolds of saddle-point invariant sets in autonomous settings \cite{wieczorek2023}.      

We define the function $g:
\bbR\rightarrow (-1,1)$ by
\begin{flalign}
    g(t):=\tanh(\alpha t/2)
\end{flalign}
and introduce a dynamical variable $s(t)\coloneqq g(t),\,t\in\bbR$. Here $\alpha > 0$ is known as the compactification parameter. It holds that 
\begin{equation*}
    \frac{\text{d}s}{\text{d}t} = \frac{\alpha}{2}\,\text{sech}^2(\alpha t/2) = \frac{\alpha}{2}\big(1-s^2\big)
\end{equation*}
Further, $g$ here is an invertible function, with inverse $g^{-1}:(-1,1)\rightarrow\bbR$ defined by
\begin{flalign*}
    g^{-1}(s)&=\frac{2}{\alpha}\,\text{artanh}(s) = \frac{1}{\alpha}\ln\left(\frac{1+s}{1-s}\right)\,,
\end{flalign*}
and accordingly one can apply the definition of $s(t)$ to deduce that $g^{-1}(s(t))=t$, or equivalently, 
$$
t = \frac{1}{\alpha}\ln\left(\frac{1+s(t)}{1-s(t)}\right)\,.
$$
A consequence of this is that the solutions of \eqref{eqn:FHN_nonautonomous_slowtime} can align exactly with the solutions for the $(u,v)$-variables in the following autonomous dynamical system yielded by the differential equation initial value problem
\begin{equation}\label{eqn:FHN_compactified_slowtime}
    \begin{cases}
        \displaystyle\frac{\text{d} u}{\text{d}t} = \eps^{-1}(-u^3 + u - v + \Lambda(s)) \vspace{0.2cm}\\ 
        \displaystyle\frac{\text{d} v}{\text{d}t} = -\beta v + u \vspace{0.2cm}\\ 
        \displaystyle\frac{\text{d}s}{\text{d}t} = r\alpha(1-s^2)/2
    \end{cases},\,\begin{pmatrix}
        u(0) \\ v(0) \\ s(0)
    \end{pmatrix}=\begin{pmatrix}
        u^{(\ti)} \\ v^{(\ti)} \\ s^{(\ti)}
    \end{pmatrix}\,,
\end{equation}
for $(u(t),v(t),s(t))\in \bbR^2\times [-1,1]\eqqcolon\pD$. Here the function $\Lambda:[-1,1]\rightarrow\mathbb{R}$ satisfies
\begin{flalign*}
    \Lambda(s) = \begin{cases} I_- & \text{if } s = -1 \\ J\big(g^{-1}(s)\big) & \text{if } -1 < s < 1 \\ I_+ & \text{if } s = 1\end{cases}\,,
\end{flalign*}
and we write the solution of the initial value problem \eqref{eqn:FHN_compactified_slowtime} using the flow notation 
\[
\bs\phi_{0,\eps,r}^\text{(c)}((u^\text{(i)},v^\text{(i)},s^\text{(i)}),t)\coloneqq (u(t),v(t),s(t)),\,t\in\mathbb{R}\,.
\]
To give an explicit example that will be applied later, the solution of \eqref{eqn:FHN_nonautonomous_slowtime} with $t^{(\ti)} = t_-/r$ and $(u^{(\ti)},v^{(\ti)}) = (u_{I_-}^{(\text{eq})},v_{I_-}^{(\text{eq})})$ aligns with the $(u,v)$-components of the solution of \eqref{eqn:FHN_compactified_slowtime} with $(u^{(\ti)},v^{(\ti)},s^{(\ti)})=(u_{I_-}^{(\text{eq})},v_{I_-}^{(\text{eq})},g(t_-))$, via the time transformation 
\[
\bs\Phi_{0,\eps,r}^{(\text{FHN})}((u_{I_-}^{(\text{eq})},v_{I_-}^{(\text{eq})}),t_-/r,t) = \begin{pmatrix}
    1 & 0 & 0 \\ 0 & 1 & 0
\end{pmatrix}\bs\phi_{0,\eps,r}^\text{(c)}((u_{I_-}^{(\text{eq})},v_{I_-}^{(\text{eq})},g(t_-)),t-t_-/r)\,.
\]

Since $s(t) \equiv \pm 1$ are fixed point solutions for the $s-$variable of $\eqref{eqn:FHN_compactified_slowtime}$, and $\Lambda(-1) = I_-$ while $\Lambda(1)=I_+$, initial conditions $(u^{(\ti)},v^{(\ti)},s^{(\ti)})$ for \eqref{eqn:FHN_compactified_slowtime} with $s^{(\ti)}=\pm1$ yield solutions whose $(u,v)$-variables align exactly with those of the autonomous system $\eqref{eqn:FHN_autonomous_slowtime}$, for $I = I_-$ and $I = I_+$, respectively. There is no rescaling of time to perform here, in that
\begin{flalign*}
\bs\phi_{0,\eps,I_-}^{(\text{FHN})}((u^{(\ti)},v^{(\ti)}),t) &= \begin{pmatrix}
    1 & 0 & 0 \\ 0 & 1 & 0
\end{pmatrix}\bs\phi_{0,\eps,r}^\text{(c)}((u^{(\ti)},v^{(\ti)},-1),t)\,, \\[0.75em] 
\bs\phi_{0,\eps,I_+}^{(\text{FHN})}((u^{(\ti)},v^{(\ti)}),t) &= \begin{pmatrix}
    1 & 0 & 0 \\ 0 & 1 & 0
\end{pmatrix}\bs\phi_{0,\eps,r}^\text{(c)}((u^{(\ti)},v^{(\ti)},1),t)\,.
\end{flalign*}
The additional inclusion of the invariant planes $s = -1$ and $s = 1$ provides the ``augmentation'' of the time parameterisation window $t \in (-\infty,\infty)$ that was alluded to earlier, in that the dynamics contained within the $s = -1$ and $s = 1$ invariant planes correspond to the dynamics of the past and future limiting autonomous systems that would be attained at $t ``=" -\infty$ and $t ``=" \infty$. In this sense, the augmented time window $t\in[-\infty,\infty]$ is brought into the compact interval $s\in[-1,1]$.

\subsection{Multiple timescale system analysis}\label{sec:nonautonomousmultipletimescaleanalysis}

\begin{figure}[ht]
    \centering
    \includegraphics[width=6.48in]{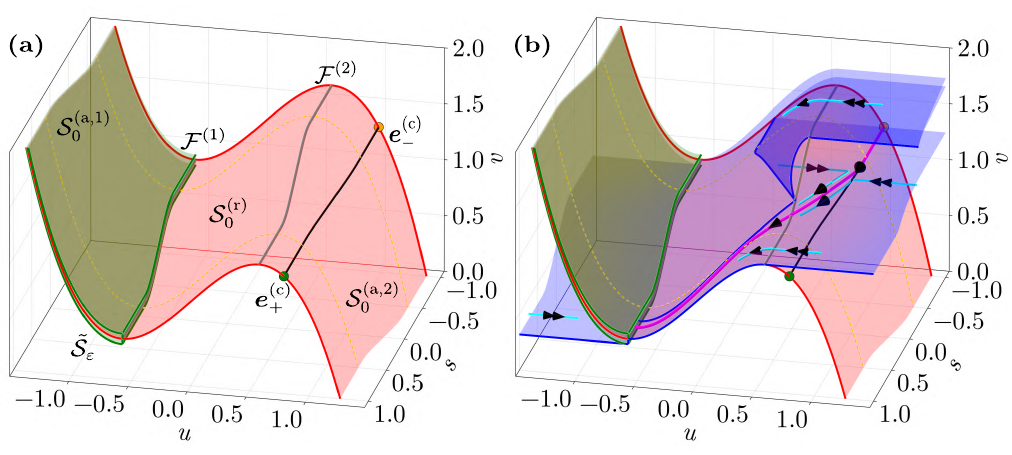}
    \caption{Schematic illustration of the metastability landscape constructed around the $L$-metastable set $\tilde{\pS}_\varepsilon$ for the basin $\tilde{\pB}$ of system \eqref{eqn:FHN_compactified_slowtime}, with $0<\varepsilon\ll 1$ and critical rate $r = r_\varepsilon^\text{(cr)}$. Both panels display: the three-dimensional $L$-metastable set $\tilde{\pS}_\eps$ [green block]; the folded two-dimensional critical manifold $\pS_0$ [red surface], comprising of the fold points $\pF^\text{(1)}$ [left grey curve, obscured by $\tilde{\pS}_\varepsilon$] and $\pF^\text{(2)}$ [right grey curve], the normally attracting sheets $\pS_0^\text{(a,1)}$ [left of $\pF^\text{(1)}$, obscured by $\tilde{\pS}_\varepsilon$] and $\pS_0^\text{(a,2)}$ [right of $\pF^{(2)}$], and the normally repelling sheet $\pS_0^\text{(r)}$ [between $\pF^{(1)}$ and $\pF^{(2)}$]; the stable node $\be_+^\text{(c)}$ [green point]; the saddle point $\be_-^\text{(c)}$ [orange point]; the moving frozen-in sinks $(u_{\Lambda(s)}^\text{(eq)},v_{\Lambda(s)}^\text{(eq)},s)$ [black curve]; and the vertical slices $s = g(t_-)$ and $s = g(t_+)$ [dashed yellow curves]. Panel (b) additionally shows part of the pullback attractor $A_{\varepsilon,r}(\be_-^\text{(c)})$ [magenta curve], and the full quasithreshold $\pQ_{\varepsilon,r}$ [blue multilayered surface]. The cyan curves indicate the flow of the system along certain parts of the quasithreshold. To view the various layers of the quasithreshold in isolation, and further understand how these surfaces separate supercritical and subcritical behaviour, see Figure \begin{NoHyper}\ref{fig:schematic_2d_quasithreshold}\end{NoHyper} in Section \begin{NoHyper}\ref{sm:decomposed_quasithreshold}\end{NoHyper} of the supplementary material\cite{supplementary_material}.}
\label{fig:schematic_2d_quasithreshold_whole}
\end{figure}

Eq. \eqref{eqn:FHN_compactified_slowtime} with $0<\varepsilon\ll 1$ describes a three-dimensional multiple timescale system, with one fast variable $u(t)$ and two slow variables $v(t)$ and $s(t)$. The first-order approximations of the dynamics of \eqref{eqn:FHN_compactified_slowtime} which play out on the system's \textit{slow} timescale $t$ are the solutions of the following differential algebraic equation initial value problem
\begin{equation}\label{eqn:FHN_compactified_slowsubsystem}
    \begin{cases}
        0 = -u^3 + u - v + \Lambda(s) \vspace{0.2cm}\\ 
        \displaystyle\frac{\text{d} v}{\text{d}t} = -\beta v + u \vspace{0.2cm}\\ 
        \displaystyle\frac{\text{d}s}{\text{d}t} = r\alpha(1-s^2)/2
    \end{cases},\,\begin{pmatrix}
        u(0) \\ v(0) \\ s(0)
    \end{pmatrix}=\begin{pmatrix}
        u^{(\ti)} \\ v^{(\ti)} \\ s^{(\ti)}
    \end{pmatrix}\,,
\end{equation}      
for $(u(t),v(t),s(t))\in\pD\coloneqq\bbR^2\times[-1,1]$. The variables of system \eqref{eqn:FHN_compactified_slowsubsystem} are restricted to the two-dimensional folded surface
\begin{equation}\label{eqn:FHN_compactified_CM}
\pS_0 \coloneqq \big\{(u,v,s)\in\pD:-u^3+u-v+\Lambda(s)=0\big\}\,,
\end{equation}
which describes the \textit{critical manifold} for \eqref{eqn:FHN_compactified_slowtime} (see Figure \ref{fig:schematic_2d_quasithreshold_whole}). Note as before that $\pS_0=\{(u,v,s)\in\pD:v=h(u,s)\}$, where $h:\bbR\times[-1,1]\rightarrow\bbR$ defined via $h(u,s) = -u^3+u+\Lambda(s)$ describes the functional form of the critical manifold. 

Contrastingly, the first-order approximation of the dynamics of \eqref{eqn:FHN_compactified_slowtime} which play out on the system's \textit{fast} timescale $\tau$ are the solutions of the differential equation initial value problem
\begin{equation}\label{eqn:FHN_compactified_fastsubsystem}
    \begin{cases}
        \displaystyle\frac{\text{d} u}{\text{d}\tau} = -u^3 + u - v + \Lambda(s) \vspace{0.2cm}\\ 
        v(\tau) \equiv v^{(\ti)} \vspace{0.2cm}\\ 
        s(\tau) \equiv s^{(\ti)}
    \end{cases},\,\begin{pmatrix}
        u(0) \\ v(0) \\ s(0)
    \end{pmatrix}=\begin{pmatrix}
        u^{(\ti)} \\ v^{(\ti)} \\ s^{(\ti)}
    \end{pmatrix}\,,
\end{equation}
for $(u(\tau),v(\tau),s(\tau))\in\pD$. There are infinitely many non-hyperbolic equilibrium points of system \eqref{eqn:FHN_compactified_fastsubsystem}, drawn out by two continuous curves
\begin{align*}
    \pF^{(1)} \coloneqq \{(u,v,s)\in\pD: u=u^{(\tf,1)},v=v_{\Lambda(s)}^{(\tf,1)}\}, \\ 
    \pF^{(2)} \coloneqq \{(u,v,s)\in\pD: u=u^{(\tf,2)},v=v_{\Lambda(s)}^{(\tf,2)}\}\,.
\end{align*}
Elements of $\pF^{(1)}\cup\pF^{(2)}$ are again examples of nondegenerate fold points\cite{kuehn2015} of the critical manifold $\pS_0$, which can also be understood as saddle-node bifurcations.

The union of the sets 
$$
\pS_0^{(\ta,1)} \coloneqq \big\{(u,v,s)\in\pD:u<u^{(\tf,1)},\,v=h(u,s)\big\}
$$
and
$$
\pS_0^{(\ta,2)} \coloneqq \big\{(u,v,s)\in\pD:u>u^{(\tf,2)},\,v=h(u,s)\big\}
$$
provide the largest attracting, normally hyperbolic subset of the critical manifold $\pS_0$. On the other hand, the set 
$$
\pS_0^{(\tr)} \coloneqq \big\{(u,v,s)\in\pD:u^{(\tf,1)} < u < u^{(\tf,2)},\,v=h(u,s)\big\}
$$
provides the largest repelling, normally hyperbolic subset of the critical manifold $\pS_0$. The set $\pS_0^{(\tr)}$ connects to the sets $\pS_0^{(\ta,1)}$ and $\pS_0^{(\ta,2)}$ via $\pF^{(1)}$ and $\pF^{(2)}$ respectively, and the critical manifold is fully comprised of the disjoint union of these parts, i.e.
\begin{equation*}
    \pS_0 = \pS_0^{(\ta,1)} \cup \pF^{(1)} \cup \pS_0^{(\tr)} \cup \pF^{(2)} \cup \pS_0^{(\ta,2)}\,.
\end{equation*}
Panel (a) of Figure \ref{fig:schematic_2d_quasithreshold_whole} depicts the critical manifold $\pS_0$ within the three-dimensional state space $\pD$, and showcases its normally attracting, normally repelling, and non-hyperbolic subsets. 

There are two equilibrium points of \eqref{eqn:FHN_compactified_slowtime}, namely $\bs e_-^{(\text{c})} \coloneqq (u_{I_-}^{(\text{eq})},v_{I_-}^{(\text{eq})},-1)$ within the invariant plane $s = -1$ and $\bs e_+^{(\text{c})} \coloneqq (u_{I_+}^{(\text{eq})},v_{I_+}^{(\text{eq})},1)$ within the invariant plane $s = 1$, both lying in the attracting sheet $\pS_0^{(\ta,2)}$ (see Figure \ref{fig:schematic_2d_quasithreshold_whole}). Given how the Jacobian of \eqref{eqn:FHN_compactified_slowtime} relates to that of \eqref{eqn:FHN_autonomous_slowtime}, one can deduce - for $\eps>0$ sufficiently small - that the equilibrium point $\bs e_-^{(\text{c})}$ has two real and negative eigenvalues (corresponding to eigenvectors whose $s$-components equal zero) and one real and positive eigenvalue (with value $r\alpha$, corresponding to an eigenvector $(0,0,1)$). On the other hand, for $\varepsilon > 0$ sufficiently small, the equilibrium point $\bs e_+^{(\text{c})}$ has three real and negative eigenvalues (two corresponding to eigenvectors whose $s$-components equal zero, and the third with value $-r\alpha$ corresponding to an eigenvector $(0,0,1)$). Therefore, $\bs e_-^{(\text{c})}$ and $\bs e_+^{(\text{c})}$ are saddle-point and stable-node type equilibria, respectively, when $0<\eps\ll 1$. We shall denote the unstable manifold of $\bs e_-^{(\text{c})}$ by $W_{\eps,r}^{(\text{u})}(\bs e_-^{(\text{c})})$, and the basin of attraction of $\bs e_+^{(\text{c})}$ by $\tilde{\pB}:=\pB(\{\be_+^\text{(c)}\}) = \bbR^2\times (-1,1]$. The resticted unstable manifold $A_{\varepsilon,r}(\be_-^\text{(c)})\coloneqq W_{\eps,r}^{(\text{u})}(\bs e_-^{(\text{c})})\setminus\{\be_-^\text{(c)}\}$ is referred to throughout the following sections as the \textit{pullback attractor}\cite{kloeden2011nonautonomous} of system \eqref{eqn:FHN_compactified_slowtime}.

\subsection{A metastability landscape for parametrically forced supercritical excursions}\label{subsec:metastability_landscape_forced}

Within the basin $\tilde{\pB}:=\pB(\{\be_+^\text{(c)}\})$ of the parametrically forced FitzHugh-Nagumo model \eqref{eqn:FHN_compactified_slowtime} there are two disconnected \textit{stable sheets} of the critical manifold, $\pS_0^{(\ta,1)}$ and $\pS_0^{(\ta,2)}$, which together produce a coexistence of qualitatively distinct transient dynamics. In this section, we formulate a metastability landscape within $\tilde{\pB}$ to explicitly define the subset of this basin which yields solutions that represent \textit{supercritical excursions} away from $\bs e_+^{(\text{c})}$. A schematic illustration of the resulting metastability landscape is given in panel (b) of Figure \ref{fig:schematic_2d_quasithreshold_whole} for a \textit{critical} rate $r$ with which the pullback attractor intersects the quasithreshold but does not intersect the supercritical set. This critical rate, which is later reasoned to be unique in Section \ref{subsec:nonautonomous_bifurcation_analysis}, is denoted by $r_\varepsilon^\text{(cr)}$. 

As in section \ref{subsec:metastability_landscape_unforced}, we base the metastability landscape around a neighbourhood of $\pS_0^{(\ta,1)}$. To construct a suitable neighbourhood, we apply an adiabatic approximation for the forcing variable $s(t)$ and continue the $L$-metastable set $\tilde{\pS}_{\eps,I}$, defined in Eq. \eqref{eqn:L-metastable_set_unforced}, across the $s$-axis by considering the union
\begin{equation}\label{eqn:nonautonomous_metastable_set}
\tilde{\pS}_\eps\coloneqq\bigcup_{s\in(-1,1]}\tilde{\pS}_{\eps,\Lambda(s)}\times\{s\}=\{(u,v,s)\in\pM\times(-1,1]:u<u^{(\tf,1)},\,\lvert v- h(u,s)\rvert<\eps^{2/3}\}\,.
\end{equation}
The set $\tilde{\pS}_\eps$ is illustrated in panel (a) of Figure \ref{fig:schematic_2d_quasithreshold_whole}. For a given $0<\eps\ll 1$, when the rate $r>0$ is sufficiently small, the neighbourhood $\tilde{\pS}_\eps$ will contain all \textit{slow manifolds} of \eqref{eqn:FHN_compactified_slowtime} corresponding to compact submanifolds of $\pS_0^{(\ta,1)}\cap\big(\pM\times(-1,1]\big)$. This is because individual solutions of \eqref{eqn:FHN_compactified_slowtime} converge to those of \eqref{eqn:FHN_autonomous_slowtime} as $r\downarrow 0$ by an application of Fenichel's second theorem (see part 2 of Theorem 3.1.4 in Ref. \cite{kuehn2015}). This itself suggests that $\tilde{\pS}_\eps$ is an $L$-metastable set for system \eqref{eqn:FHN_compactified_slowtime} with $L \overset{\eps\downarrow 0}{=} \pO(1)$ when $r = r (\varepsilon)$ is chosen sufficiently small for each $0<\varepsilon\ll 1$, by reasoning analogous to points (i)-(iii) in section \ref{subsec:metastability_landscape_unforced}.  

The inclusion of the above paragraph acts solely as a motivation. It is now important to clarify that, for what follows, we will actually be interested in considering noninfinitesimal (``order 1'') values of the rate $r$. For the particular paired values of the $\varepsilon$ and $r$ parameters considered in section \ref{subsec:nonautonomous_bifurcation_analysis} and beyond, numerical solutions of system \eqref{eqn:FHN_compactified_slowtime} corroborate that the slow manifolds associated with $\pS_0^{(\ta,1)}\cap\big(\pM\times(-1,1]\big)$ are indeed contained in $\tilde{\pS}_\eps$ (see Figure \begin{NoHyper}\ref{fig:schematic_2d_quasithreshold}\end{NoHyper} of the supplementary material\cite{supplementary_material}). As a result, in these instances, the set $\tilde{\pS}_\varepsilon$ persists to capture transient regimes which are long-lived relative to the fast timescale of \eqref{eqn:FHN_compactified_slowtime}.

With this in mind, the supercritical and subcritical sets within the basin $\tilde{\pB}$ of system \eqref{eqn:FHN_compactified_slowtime} are then defined with respect to this $L$-metastable set $\tilde{\pS}_\eps$ in the standard way, 
\begin{flalign*}
\pD_{\eps,r}^{(\text{sup})}&\coloneqq\{(u,v,s)\in\tilde{\pB}:\exists\, T\in\bbR_+\text{ with } \bs\phi_{0,\eps,r}^\text{(c)}((u,v,s),T)\in\tilde{\pS}_\eps\}\,, \\[-1em]
\\ 
\pD_{\eps,r}^{(\text{sub})}&\coloneqq\{(u,v,s)\in\tilde{\pB}:\forall\, t\in\bbR_{\geq 0},\,\bs\phi_{0,\eps,r}^\text{(c)}((u,v,s),t)\not\in\text{cl}(\tilde{\pS}_\eps)\}\,,
\end{flalign*}
and the corresponding quasithreshold is defined as
\begin{equation*}
    \pQ_{\eps,r}\coloneqq\{(u,v,s)\in\tilde{\pB}\setminus\pD_{\varepsilon,r}^\text{(sup)}:\exists\, T\in\bbR_+\text{ with } \bs\phi_{0,\eps,r}^\text{(c)}((u,v,s),T)\in\text{cl}(\tilde{\pS}_\eps)\}\,. 
\end{equation*}
Recall that the quasithreshold satisfies
\begin{flalign*}
\pQ_{\varepsilon,r} = \tilde{\pB} \cap \partial\pD_{\eps,r}^{(\text{sup})} & = \tilde{\pB} \cap \partial\big(\tilde{\pB}\setminus\pD_{\eps,r}^{(\text{sup})}\big) \\ 
& \supseteq \tilde{\pB}\cap\partial\pD_{\eps,r}^\text{(sub)}\,.
\end{flalign*}
The superset symbol becomes an equality symbol if (and only if) $\partial\,\text{int}\big(\tilde{\pB}\setminus\pD_{\eps,r}^{(\text{sup})}\big) = \partial\big(\tilde{\pB}\setminus\pD_{\eps,r}^{(\text{sup})}\big)$, i.e., there are no ``thin'' parts of the non-supercritical set within the basin of attraction which are not in the closure of its interior. This \textit{is} the case for the proposed metastability landscape of the basin $\tilde{\pB}$ of system \eqref{eqn:FHN_compactified_slowtime}; as such, the quasithreshold simply describes the shared boundary of the supercritical and subcritical sets in the basin $\tilde{\pB}$ of system \eqref{eqn:FHN_compactified_slowtime}. See panel (b) of Figure \ref{fig:schematic_2d_quasithreshold_whole} for a schematic illustration of $\pD_{\varepsilon,r}^\text{(sup)},\,\pD_{\varepsilon,r}^\text{(sub)}$ and $\pQ_{\varepsilon,r}$ in the metastability landscape of the basin $\tilde{\pB}$ of system \eqref{eqn:FHN_compactified_slowtime} with $r = r_\varepsilon^\text{(cr)}$.

The spatial configuration of the quasithreshold $\pQ_{\varepsilon,r}$ of the basin $\tilde{\pB}$ is determined in part by trajectories of \eqref{eqn:FHN_compactified_slowtime} which closely follow $\pS_0^{(\tr)}$ before tangentially grazing the $L$-metastable $\tilde{\pS}_\varepsilon$ set along the curve of fold points $\pF^{(1)}$. Near $\pF^{(1)}$, the union of these trajectories creates a surface which separates the initial conditions of \eqref{eqn:FHN_compactified_slowtime} into those which are either attracted towards or repelled away from $\tilde{\pS}_\eps$ (see panels (a)-(c) of Figure \begin{NoHyper}\ref{fig:schematic_2d_quasithreshold}\end{NoHyper} in the supplementary material\cite{supplementary_material}). There is a second set of trajectories of \eqref{eqn:FHN_compactified_slowtime} which tangentially graze $\tilde{\pS}_\eps$ at other positions separating its inflow and outflow boundaries, near to the \textit{kink points} of $\tilde{\pS}_\varepsilon$ along 
$$
\{(u,v,s)\in\pD: u=u^{(\tf,1)},v=v_{\Lambda(s)}^{(\tf,1)}-\varepsilon^{2/3}\}\,.
$$
The union of these trajectories similarly creates a surface which separates the initial conditions of \eqref{eqn:FHN_compactified_slowtime} into those which are either attracted towards or repelled away from $\tilde{\pS}_\eps$ (see panel (d) of Figure \begin{NoHyper}\ref{fig:schematic_2d_quasithreshold}\end{NoHyper} in the supplementary material\cite{supplementary_material}). The supercritical set is then determined in full by the backward time propagation of these surfaces into the parts of the basin $\tilde{\pB}$ away from $\pF^{(1)}$. 

Unlike within the unforced system \eqref{eqn:FHN_autonomous_slowtime}, the configuration of the quasithreshold is complicated by an interaction between the surface which grazes $\pF^{(1)}$ and \textit{canard trajectories}, which begin to emerge within \eqref{eqn:FHN_compactified_slowtime} for $r\in\bbR_+$ sufficiently large. These canard trajectories - which follow an attracting slow manifold of $\pS_0^{(\ta,2)}$, pass close by a fold point along $\pF^{(2)}$, and subsequently follow a repelling slow manifold of $\pS_0^{(\tr)}$ \cite{szmolyan2001} - create an additional threshold, local to $\pS_0^{(\ta,2)}$, separating supercritical and subcritical behaviour. Globally, the quasithreshold takes the form of an intricate multilayered surface which combines these three distinct threshold components. A more detailed set-theoretic description of the individual components comprising the quasithreshold surface - in the metastability landscape of the basin $\tilde{\pB}$ of system \eqref{eqn:FHN_compactified_slowtime} with $r = r_\varepsilon^\text{(cr)}$ - is given in Section \begin{NoHyper}\ref{sm:decomposed_quasithreshold}\end{NoHyper} of the supplementary material\cite{supplementary_material}. There, they are decomposed and schematically illustrated in isolation throughout Figure \begin{NoHyper}\ref{fig:schematic_2d_quasithreshold}\end{NoHyper}. 

The following two sections are dedicated to explaining how canard trajectories both: (i) come about within system \eqref{eqn:FHN_compactified_slowtime}; and (ii) influence the global structure of the quasithreshold of the metastability landscape of $\tilde{\pB}$ by contributing an additional component. Specifically, within Section \ref{subsec:geometric_singular_perturbation_theory} we describe the \textit{slow} dynamics of \eqref{eqn:FHN_compactified_slowtime} local to the critical manifold $\pS_0$. Then, within Sections \ref{subsec:geometric_singular_perturbation_theory} and \ref{subsec:nonautonomous_bifurcation_analysis}, we comment on how the slow dynamics connects with the \textit{fast} dynamics of \eqref{eqn:FHN_compactified_slowtime} away from the critical manifold to form the multilayered quasithreshold in the full-dimensional state space.  

\subsection{Singular canard trajectories and R-tipping quasithresholds}\label{subsec:geometric_singular_perturbation_theory}

To gain insight into the rates $r$ for which canard trajectories are present in \eqref{eqn:FHN_compactified_slowtime}, but also their geometric configuration in state space, it is useful to first study the flow of \eqref{eqn:FHN_compactified_slowsubsystem} (the slow subsystem of \eqref{eqn:FHN_compactified_slowtime}). Subsequently, one can draw certain conclusions about the flow of \eqref{eqn:FHN_compactified_slowtime} for $0<\varepsilon\ll1$ from that of \eqref{eqn:FHN_compactified_slowsubsystem} by applying persistence results from singular perturbation theory (stemming from Fenichel theory). To study the flow of the slow subsystem, we will first derive a \textit{reduced form} for \eqref{eqn:FHN_compactified_slowsubsystem} in terms of a planar system which equivalently describes its dynamics. We then identify the dynamical object within this planar system which corresponds to the pullback attractor $A_{\varepsilon,r}(\be_-^\text{(c)})$, and describe its changing properties throughout the various stages in its transition from subcritical to supercritical behaviour as the $r$-value increases. We explain how this change in behaviour is understood via the technique of desingularisation. In the following section \ref{subsec:nonautonomous_bifurcation_analysis}, the persistence of these different behavioural regimes within the full system \eqref{eqn:FHN_compactified_slowtime} with $0<\varepsilon\ll 1$ is addressed.

To derive a reduced form for \eqref{eqn:FHN_compactified_slowsubsystem}, note that its algebraic constraint $-u^3 + u - v + \Lambda(s)=0$ can be differentiated with respect to $t$ to give
\begin{equation*}
    -3u^2\frac{\text{d}u}{\text{d}t}+\frac{\text{d}u}{\text{d}t}-\frac{\text{d}v}{\text{d}t}+\frac{\text{d}\Lambda}{\text{d}s}\frac{\text{d}s}{\text{d}t}=0\,.
\end{equation*}
For circumstances in which $1-3u(t)^2\neq 0$, this can then be rearranged to explicitly express $\text{d}u/\text{d}t$ as 
\begin{equation}\label{eqn:slowsubsystem_u_derivative}
    \frac{\text{d}u}{\text{d}t} = \frac{-\beta v+u-r\alpha\Lambda'(s)(1-s^2)/2}{1-3u^2}\,.
\end{equation}
Here $\Lambda'(s)$ denotes the derivative of $\Lambda(s)$ with respect to $s$, which is given by
\begin{flalign*}
    \Lambda'(s) = \begin{cases} 0 & \text{if } s = -1 \\ \displaystyle\frac{\text{d}J(t)}{\text{d}t}\bigg\rvert_{t=g^{-1}(s)}\frac{\text{d}g^{-1}}{\text{d}s} & \text{if } -1 < s < 1 \\ 0 & \text{if } s = 1\end{cases}\,,
\end{flalign*}
where the derivatives of the functions $J$ and $g^{-1}$ are explicitly calculable. Since the algebraic constraint $-u^3 + u - v + \Lambda(s)=0\iff v=-u^3+u+\Lambda(s)$ of \eqref{eqn:FHN_compactified_slowsubsystem} yields an explicit unique expression for $v$ in terms of $u$ and $s$, this can be substituted into \eqref{eqn:slowsubsystem_u_derivative} to describe $\text{d}u/\text{d}t$ solely in terms of $u$ and $s$ as
\begin{equation*}
    \frac{\text{d}u}{\text{d}t} = \frac{-\beta(-u^3+u+\Lambda(s))+u-r\alpha\Lambda'(s)(1-s^2)/2}{1-3u^2}\,.
\end{equation*}
This allows the differential algebraic equation initial value problem \eqref{eqn:FHN_compactified_slowsubsystem} to be re-expressed in the following reduced form
\begin{equation}\label{eqn:FHN_compactified_slowsubsystem_reduced}
    \begin{cases}
        \displaystyle\frac{\text{d}u}{\text{d}t} = \frac{-\beta(-u^3+u+\Lambda(s))+u-r\alpha\Lambda'(s)(1-s^2)/2}{1-3u^2} \vspace{0.2cm}\\ 
        \displaystyle\frac{\text{d}s}{\text{d}t} = r\alpha(1-s^2)/2
    \end{cases}\,,\,\begin{pmatrix}
        u(0) \\ s(0)
    \end{pmatrix}=\begin{pmatrix}
        u^{(\ti)} \\ s^{(\ti)}
    \end{pmatrix}\,,
\end{equation}
for $(u(t),s(t))\in\bbR\times[-1,1],\,1-3u(t)^2\neq 0$. When paired with $v(t) = -u(t)^3+u(t)+\Lambda(s(t))$, the reduced system \eqref{eqn:FHN_compactified_slowsubsystem_reduced} provides a complete description of the dynamics of \eqref{eqn:FHN_compactified_slowsubsystem} (the singularities at $1-3u^2=0$ are common to each description). 

For purposes of future reference, let 
\[\pF^{(1,\text{red})}\coloneqq\{u^{(\tf,1)}\}\times[-1,1]\,,\quad\pF^{(2,\text{red})}\coloneqq\{u^{(\tf,2)}\}\times[-1,1]\,,\]
which respectively correspond to the subsets $\pF^{(1)}$ and $\pF^{(2)}$ of non-hyperbolic fold points of the critical manifold $\pS_0$ of \eqref{eqn:FHN_compactified_slowtime}. Further, let 
\[\pS_0^{(\text{a},1,\text{red})}\coloneqq(-\infty,u^{(\tf,1)})\times[-1,1]\,,\quad\pS_0^{(\text{a},2,\text{red})}\coloneqq(u^{(\tf,2)},\infty)\times[-1,1]\,,\]
which respectively correspond to the normally hyperbolic attracting subsets $\pS_0^{(\text{a},1)}$ and $\pS_0^{(\text{a},2)}$ of $\pS_0$. Finally, let 
\[\pS_0^{(\text{r},\text{red})}\coloneqq(u^{(\tf,1)},u^{(\tf,2)})\times[-1,1]\,,\]
which corresponds to the normally hyperbolic repelling subset $\pS_0^{(\text{r})}$ of $\pS_0$.

In correspondence with the equilibria of \eqref{eqn:FHN_compactified_slowtime}, there is a saddle point 
\begin{equation*}
    \bs e_{-}^{(\tc,\text{red})}\coloneqq (u_{I_-}^{(\text{eq})},-1)     
\end{equation*} 
contained within an invariant set $\{(u,s)\in\pS_0^{(\ta,2,\text{red})}:s=-1\}$, and a stable node 
\begin{equation*}
    \bs e_{+}^{(\tc,\text{red})}\coloneqq (u_{I_+}^{(\text{eq})},1)
\end{equation*}
contained within an invariant set $\{(u,s)\in\pS_0^{(\ta,2,\text{red})}:s=1\}$ for the reduced system $\eqref{eqn:FHN_compactified_slowsubsystem_reduced}$. Let $W_{0,r}^{(\text{u})}(\bs e_-^{(\tc,\text{red})})$ denote the unstable manifold of the saddle point $\bs e_-^{(\tc,\text{red})}$ in $\eqref{eqn:FHN_compactified_slowsubsystem_reduced}$, which is a reduced form of the singular limit of $W_{\eps,r}^{(\tu)}(\bs e_-^{(\text{c})})$ as $\varepsilon \downarrow 0$. Accordingly, analyzing how the corresponding pullback attractor $A_{0,r}(\be_-^\text{(c,red)})\coloneqq W_{0,r}^{(\text{u})}(\bs e_-^{(\tc,\text{red})})\setminus\{\be_-^\text{(c,red)}\}$ changes with $r$ unveils the mechanism behind the transition from subcritical to supercritical behaviour which can occur in the full system $\eqref{eqn:FHN_compactified_slowtime}$ given suitable values of $I_-$ and $I_+$. To this end, from this point onward we shall fix 
\begin{align}\label{eqn:IminusIplus_values}
    t_-=-0.5\,, t_+ = 0.5\,, && I_- = 1.25\,, I_+ = 0.75\,, && \alpha = 2\ln3\,,
\end{align}
without full loss of generality. With parameter values as in \eqref{eqn:IminusIplus_values}, the pullback attractor $A_{0,r}(\be_-^\text{(c,red)})$ of system \eqref{eqn:FHN_compactified_slowsubsystem_reduced} assumes three qualitatively distinct forms as the $r$-value grows, which require the following classification of fold points in \eqref{eqn:FHN_compactified_slowsubsystem_reduced} in order to describe. Defining 
\begin{equation*}
    n_r(u,s)\coloneqq-\beta(-u^3+u+\Lambda(s))+u-r\alpha\Lambda'(s)(1-s^2)/2\,,
\end{equation*}
a fold point $(u,s)\in\pF^{(1,\text{red})}\cup\pF^{(2,\text{red})}$ is referred to either as a \textit{jump point}\cite{szmolyan2001} if $n_r(u,s)\neq 0$, or a \textit{folded singularity}\cite{szmolyan2001} if $n_r(u,s)=0$. 

\begin{figure}[ht]
    \centering
    \includegraphics[width=5.1in]{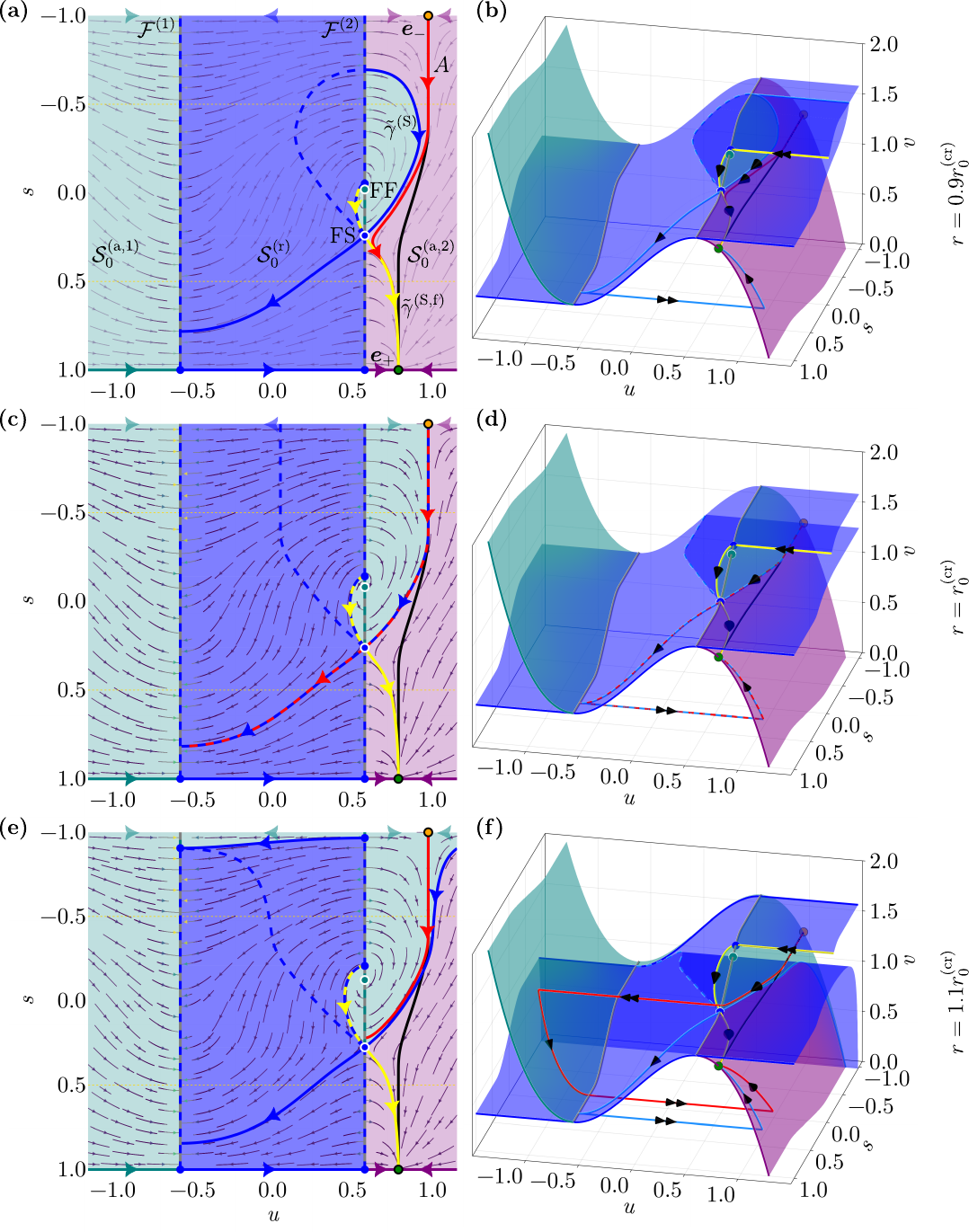}
    \caption{Quasithresholds $\pQ_{0,r}$ in the singular system \eqref{eqn:FHN_compactified_slowtime} with $\varepsilon=0$, through sequential stages in the tracking-tipping transition of the pullback attractor $A_{0,r}(\be_-^\text{(c,red)})$ around $r = r_0^{(\text{cr})}$. Left column: phase portraits of the reduced slow subsystem \eqref{eqn:FHN_compactified_slowsubsystem_reduced} for $r=0.9r_0^{(\text{cr})},\,r=r_0^{(\text{cr})},\,r=1.1r_0^{(\text{cr})}$. The moving sinks $(u_{\Lambda(s)}^\text{(eq)},s)$ [black], the pullback attractor $A_{0,r}(\be_-^\text{(c,red)})$ [red], the faux folded saddle singular canard trajectory $\tilde{\gamma}_{0,r}^\text{(S,f)}$ from $\pS_0^{(\text{r,red})}$ to $\pS_0^{(\text{a},2,\text{red})}$ [yellow], and the folded saddle singular canard trajectory $\tilde{\gamma}_{0,r}^{(\text{S})}$ from $\pS_0^{(\text{a},2,\text{red})}$ to $\pS_0^{(\text{r,red})}$ [blue] are shown. The teal, purple, and blue regions show - in reduced coordinates - the parts of the critical manifold $\pS_0$ that are contained in the supercritical set, the subcritical set, and the quasithreshold, respectively. Right column: the components of the quasithresholds $\pQ_{0,r}$ associated with the fast and slow flow [blue surfaces], alongside the critical manifold $\pS_0$ [folded surface coloured as in the left column]. Candidate trajectories containing the folded saddle singular canard trajectory [blue], the pullback attractor [red] and the faux folded saddle singular canard trajectory [yellow] are shown. Note in particular that $\pQ_{0,r}\neq\tilde{\pB}\cap\partial\pD_{0,r}^\text{(sub)}$ in each case, due to the ``thin part'' of $\tilde{\pB}\setminus\pD_{0,r}^\text{(sup)}$ emanating from $\pF_0^{(2)}$.}
\label{fig:singular_metastability_landscapes}
\end{figure}

Now, for $r>0$ sufficiently small, the pullback attractor $A_{0,r}(\be_-^\text{(c,red)})$ forms a heteroclinic connection between $\bs e_-^{(\tc,\text{red})}$ and $\bs e_+^{(\tc,\text{red})}$. For these $r$ values the pullback attractor $A_{0,r}(\be_-^\text{(c,red)})$ is fully contained within $\pS_0^{(\text{a},2,\text{red})}$. See panel (a) in Figure \ref{fig:singular_metastability_landscapes} for an illustration of this. However, at a unique value $r = r_0^{(\text{cr})}\in\bbR_+$, the pullback attractor $A_{0,r}(\be_-^\text{(c,red)})$ ceases to form a heteroclinic connection between $\bs e_-^{(\tc,\text{red})}$ and $\bs e_+^{(\tc,\text{red})}$. Instead, it travels from $\pS_0^{(\text{a},2,\text{red})}$ across to $\pS_0^{(\text{r},\text{red})}$, passing through a folded singularity within $\pF^{(2,\text{red})}$ at finite velocity. The pullback attractor $A_{0,r}(\be_-^\text{(c,red)})$ never leaves $\pS_0^{(\text{r},\text{red})}$ after reaching it but eventually ceases to exist within $\bbR\times[-1,1]$ (at finite time) by limiting towards a jump point within $\pF^{(1,\text{red})}$ at velocity tending to infinity. See panel (b) in Figure \ref{fig:singular_metastability_landscapes} for an illustration of this. Finally, for $r>r_0^{(\text{cr})}$, the pullback attractor $A_{0,r}(\be_-^\text{(c,red)})$ neither forms a heteroclinic orbit between $\bs e_-^{(\tc,\text{red})}$ and $\bs e_+^{(\tc,\text{red})}$ nor travels across to $\pS_0^{(\text{r},\text{red})}$. It remains fully contained within $\pS_0^{(\text{a},2,\text{red})}$ but ceases to exist within $\bbR\times[-1,1]$ (at finite time) by limiting towards a jump point within $\pF^{(2,\text{red})}$ at velocity tending to infinity. See panel (c) in Figure $\ref{fig:singular_metastability_landscapes}$ for an illustration of this. For the given parameterisation \eqref{eqn:IminusIplus_values}, one has $r_0^\text{(cr)}\approx 0.5277897$.  

The sequential progression in $r$ between these forms for $A_{0,r}(\be_-^\text{(c,red)})$ is labelled the \textit{simple slow} case of ``tracking-tipping transition'' within \citet{osullivan2023}, and can be explained via a global bifurcation analysis of a desingularised version of system \eqref{eqn:FHN_compactified_slowsubsystem_reduced}. System \eqref{eqn:FHN_compactified_slowsubsystem_reduced} possesses singularities at its fold points, because $1-3u^2=0$ for all $(u,s)\in\pF^{(1,\text{red})}\cup\pF^{(2,\text{red})}$, and hence the time-derivative of the $u$-variable is ill-defined there. One can \textit{desingularise} the system by introducing the following intricate rescaling of time $z\mapsto t(z)$, whose functional form is the solution to the following differential equation
\begin{equation}\label{eqn:intricate_time_rescaling}
    \frac{\text{d}t}{\text{d}z}=3u(t)^2-1,\,\,t(0)=0\,,
\end{equation}
for $t(z)\in\mathbb{R},\,z\in(z_-,z_+)$ where $-\infty\leq z_-<0<z_+\leq\infty$. Here, $u(t)$ is the $u$-component of a solution to \eqref{eqn:FHN_compactified_slowsubsystem_reduced} with an initial condition away from the fold points, i.e., $u^\text{(i)}\not\in \{u^{(\tf,1)},\,u^{(\tf,2)}\}$; hence, $t(z)$ is specific to each trajectory and its derivative depends on the proximity of the trajectory to the fold points. One can then introduce 
\begin{equation}\label{eqn:intricate_change_of_variables}
 \hat{u}(z)=u(t(z))\,,\quad\hat{s}(z)=s(t(z))
\end{equation}
as the associated state variables which are the time-rescaled versions of $u$ and $s$; it follows that $\hat{u}$ and $\hat{s}$ satisfy the following system of differential equations
\begin{equation}\label{eqn:FHN_compactified_slowsubsystem_reduced_desingularised}
    \begin{cases}
        \displaystyle\frac{\text{d}\hat{u}}{\text{d}z} = \beta\big(-\hat{u}^3+\hat{u}+\Lambda(\hat{s})\big)-\hat{u}+r\alpha\Lambda'(\hat{s})(1-\hat{s}^2)/2 \vspace{0.2cm}\\ 
        \displaystyle\frac{\text{d}\hat{s}}{\text{d}z} = r\alpha(1-\hat{s}^2)(3\hat{u}^2-1)/2
    \end{cases},
    \,\begin{pmatrix}
       \hat{u}(0) \\ \hat{s}(0)
    \end{pmatrix}=\begin{pmatrix}
       \hat{u}^{(\ti)} \\ \hat{s}^{(\ti)}
    \end{pmatrix} 
\end{equation}
for $(\hat{u}(z),\hat{s}(z))\in\bbR\times[-1,1]$. It also follows by definitions \eqref{eqn:intricate_time_rescaling}-\eqref{eqn:intricate_change_of_variables} that, across an interval for $z$ in which a trajectory of \eqref{eqn:FHN_compactified_slowsubsystem_reduced_desingularised} satisfies 
$$
(\hat{u}(z),\hat{s}(z))\in\pS_0^{(\ta,1,\text{red})}\cup\pS_0^{(\ta,2,\text{red})}\iff 3\hat{u}(z)^2-1>0\,,
$$
the functional rescaling of time $t(z)$ with respect to which this trajectory is defined satisfies $\text{d}t/\text{d}z>0$. Hence, the ``direction of flow'' of the trajectories of \eqref{eqn:FHN_compactified_slowsubsystem_reduced_desingularised} within $\pS_0^{(\ta,1,\text{red})}\cup\pS_0^{(\ta,2,\text{red})}$ matches the ``direction of flow'' of the corresponding trajectories in system $\eqref{eqn:FHN_compactified_slowsubsystem_reduced}$. On the other hand, across an interval for $z$ in which a trajectory of \eqref{eqn:FHN_compactified_slowsubsystem_reduced_desingularised} satisfies 
$$
(\hat{u}(z),\hat{s}(z))\in\pS_0^{(\tr,\text{red})}\iff 3\hat{u}(z)^2-1<0\,,
$$
the functional rescaling of time $t(z)$ with respect to which this trajectory is defined satisfies $\text{d}t/\text{d}z<0$. Hence, the ``direction of flow'' of the trajectories of \eqref{eqn:FHN_compactified_slowsubsystem_reduced_desingularised} within $\pS_0^{(\tr,\text{red})}$ reverses the ``direction of flow'' of the corresponding trajectories in system $\eqref{eqn:FHN_compactified_slowsubsystem_reduced}$. 

We refer to \citet{osullivan2023} and Section \begin{NoHyper}\ref{sm:bifofreduced}\end{NoHyper} of the supplementary material\cite{supplementary_material} for a discussion of the bifurcations which take place within the reduced desingularised system \eqref{eqn:FHN_compactified_slowsubsystem_reduced_desingularised} as the rate $r$ is increased. In particular, there \textit{folded saddle singular canard trajectories} $\tilde{\gamma}_{0,r}^{(\text{S})}$ and \textit{faux folded saddle singular canard trajectories} $\tilde{\gamma}_{0,r}^{(\text{S,f})}$ that exist within system \eqref{eqn:FHN_compactified_slowsubsystem_reduced} are introduced, which are plotted in Figure \ref{fig:singular_metastability_landscapes}. Faux canard trajectories\cite{benoit1983systemes} are those which - in the reverse sense to regular canard trajectories - initially follow the repelling part of the critical manifold, pass close by a fold point, and subsequently follow the attracting part of the critical manifold. 

One can pair the dynamics of the slow subsystem \eqref{eqn:FHN_compactified_slowsubsystem} with the dynamics of the fast subsystem \eqref{eqn:FHN_compactified_fastsubsystem} as a means of constructing metastability landscapes for the basin $\tilde{\pB}$ within the singular limit of system \eqref{eqn:FHN_compactified_slowtime} with $\varepsilon = 0$. These are illustrated in Figure \ref{fig:singular_metastability_landscapes} - across panels (b), (d), and (f) for subcritical, critical, and supercritical rates, respectively - and explained directly below.

The metastability landscapes for the basin $\tilde{\pB}$ within the singular system are taken with respect to the limit of the $L$-metastable sets $\tilde{\pS}_\eps$ as $\varepsilon\downarrow 0$, namely  
$$
\tilde{\pS}_0\coloneqq\bigcap_{\varepsilon\in\bbR_+}\tilde{\pS}_\eps =\pS_0^{(\ta,1)}\cap\big(\pM\times (-1,1]\big)\,,
$$
The corresponding supercritical set is defined by
$$
\pD_{0,r}^{(\text{sup})}\coloneqq\{\bs z\in\tilde{\pB}:\forall\,\bs\gamma\in\mathcal{C}_0([0,1])\text{ with }(\bs\gamma(0),\bs\gamma(1))=(\bs z,\be_+^{(\tc)}),\,\text{Im}(\bs\gamma)\cap\tilde{\pS}_0\neq\emptyset\}\,,
$$
where $\mathcal{C}_0([0,1])$ is the set of functions whose images are \textit{candidate orbits}\footnote{A candidate orbit of the singular limit of the system \eqref{eqn:FHN_compactified_slowtime} with $\varepsilon=0$, parameterised on the interval $[0,1]$, is defined\cite{kuehn2011} as the image of a homeomorphic function $\bs\gamma_0:[0,1]\rightarrow\bbR^2\times[-1,1]$ with the following properties: there exists some partition $0 = t_0 < t_1 < t_2 < ... < t_{p-1} < t_p = 1$ of $[0,1]$ whereby the segment $\bs\gamma_0((t_{j-1},t_{j}))$ aligns with a trajectory of either the slow or fast subsystems of \eqref{eqn:FHN_compactified_slowtime} (see \eqref{eqn:FHN_compactified_slowsubsystem} and \eqref{eqn:FHN_compactified_fastsubsystem}), for each $j\in\{1,2,...,p\}$; in addition, the orientation (``direction of flow'') along each segment is the same as that of the corresponding trajectory of the slow or fast subsystem.} of the singular limit of the system \eqref{eqn:FHN_compactified_slowtime} with $\varepsilon=0$, parameterised on $[0,1]$. Therefore, $\pD_{0,r}^{(\text{sup})}$ is the set of initial conditions which yield candidate orbits that can only reach $\be_+^{(\text{c})}$ after entering $\tilde{\pS}_0$. On the other hand, the corresponding subcritical set is defined by
$$
\pD_{0,r}^{(\text{sub})}\coloneqq\text{int}\big(\tilde{\pB}\setminus\pD_{0,r}^{(\text{sup})}\big)\,,
$$
which is comprised of the initial conditions within $\tilde{\pB}$ that: (i) yield candidate orbits that can reach $\be_+^{(\text{c})}$ without entering $\tilde{\pS}_0$; and (ii) are away from the boundary of the supercritical set $\pD_{0,r}^\text{(sup)}$. The corresponding quasithreshold of the singular system is then defined as 
$$
\pQ_{0,r}\coloneqq \tilde{\pB}\cap\partial\pD_{0,r}^{(\text{sup})}\,.
$$
For the parametrically forced FitzHugh-Nagumo system considered here, an explicit formulation of $\pD_{0,r}^{(\text{sup})}$ is given in Section \begin{NoHyper}\ref{sm:singular_supercritical_set}\end{NoHyper} of the supplementary material\cite{supplementary_material}.

\subsection{Nonautonomous bifurcation analysis (nonsingular canard trajectories and R-tipping quasithresholds)}\label{subsec:nonautonomous_bifurcation_analysis}

\begin{figure}[ht]
    \centering
    \includegraphics[width=6in]{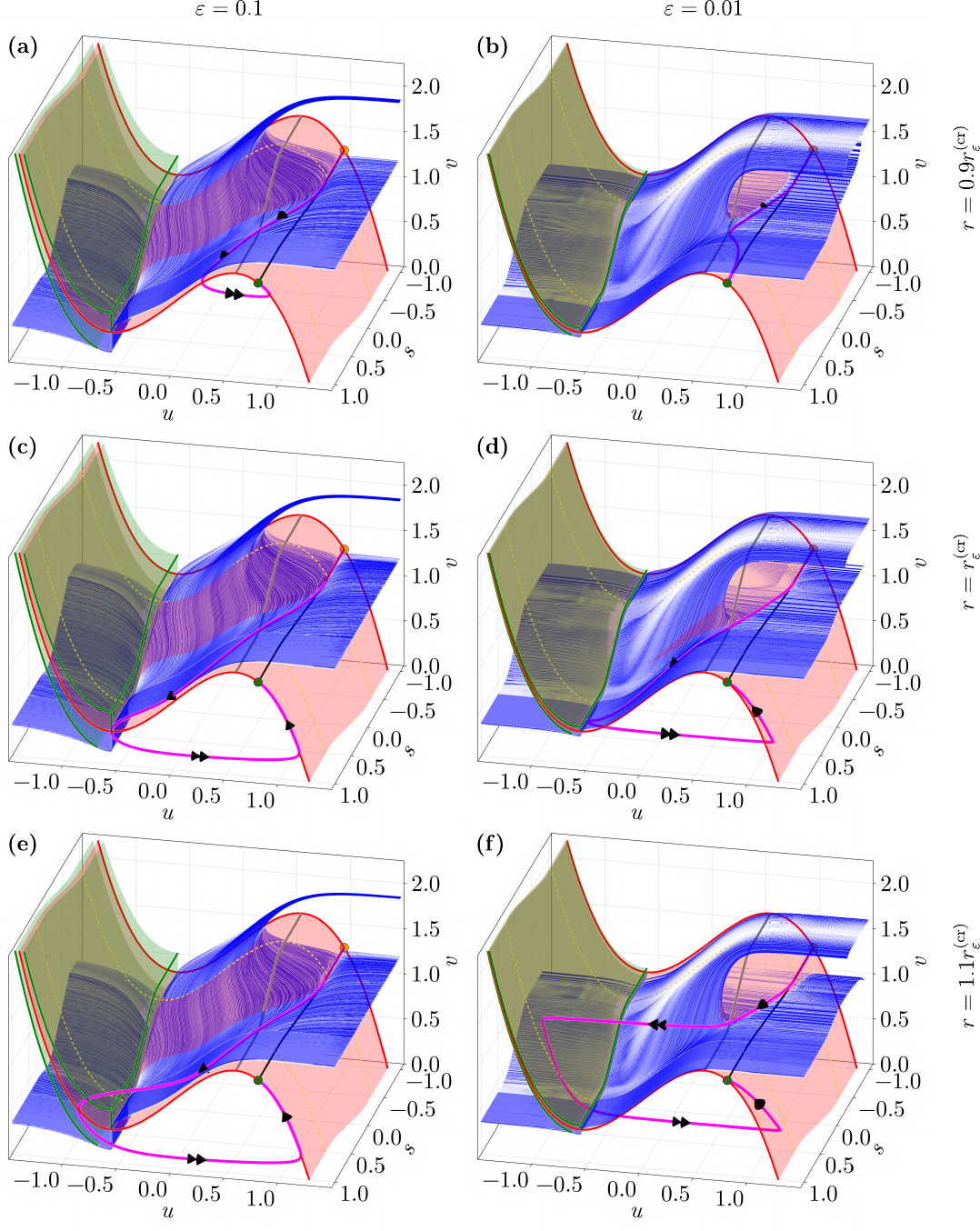}
    \caption{Quasithresholds $\pQ_{\eps,r}$ in system \eqref{eqn:FHN_compactified_slowtime}, through sequential stages in the tracking-tipping transition of the pullback attractor $A_{\varepsilon,r}(\be_-^\text{(c)})$ at $r = r_\varepsilon^{(\text{crit})}$. The left and right columns correspond to $\varepsilon = 0.1$ and $\varepsilon = 0.01$, respectively, while the top, middle and bottom rows correspond to $r = 0.9r_\varepsilon^{(\text{crit})}$, $r = r_\varepsilon^{(\text{crit})}$, and $r = 1.1r_\varepsilon^{(\text{crit})}$, respectively. Every individual panel contains the $L-$metastable set $\tilde{\pS}_{\eps}$ [green surface with emboldened outline], the critical manifold [red transparent surface], the pullback attractor $A_{\varepsilon,r}(\be_-^\text{(c)})$ [magenta curve], and the moving sinks of the frozen-in system $\eqref{eqn:FHN_autonomous_slowtime}$ with $I = \Lambda(s)$ [black curve]. Also shown is a numerical approximation of the quasithreshold $\pQ_{\eps,r}$ [blue multilayered surface].}
    \label{fig:nonautonomous_metastability_landscapes}
\end{figure}

Recall that within this section we aim to outline the interaction between the dynamics of the singularly perturbed, parametrically forced FitzHugh-Nagumo system \eqref{eqn:FHN_compactified_slowtime} which occurs on its fast and slow timescales, in order to explain the geometric structure of the quasithreshold $\pQ_{\varepsilon,r}$ with $0<\varepsilon\ll 1$. We approach this by breaking down the individual components of the quasithreshold and describing them sequentially throughout the following. Ultimately, numerical approximations of the quasithreshold $\pQ_{\varepsilon,r}$ - in the metastability landscape of the basin $\tilde{\pB}$ of system \eqref{eqn:FHN_compactified_slowtime} - are collected in Figure \ref{fig:nonautonomous_metastability_landscapes} across different rates with which the pullback attractor $A_{\varepsilon,r}(\be_-^\text{(c)})$ corresponds to a subcritical, critical, and supercritical trajectory.   

Fenichel theory, as discussed in section \ref{sec:fitzhughnagumo_autonomous_fastslow_analysis}, directly implies that the flow of \eqref{eqn:FHN_compactified_slowsubsystem} along compact subsets of $\pS_0^{(\ta,1)}\cup\pS_0^{(\tr)}\cup\pS_0^{(\ta,2)}$ provides a first-order approximation for trajectories of \eqref{eqn:FHN_compactified_slowtime} with $\varepsilon > 0$ sufficiently small. However, it is also true that the folded saddle singular canard trajectory $\tilde{\gamma}_{0,r}^{(\text{S})}$ of \eqref{eqn:FHN_compactified_slowsubsystem} - which travels through the non-hyperbolic fold points $\pF^{(2)}$ - perturbs to a family of canard trajectories within \eqref{eqn:FHN_compactified_slowtime} for $\varepsilon>0$ sufficiently small \cite{szmolyan2001}. This is to say that there are slow manifolds of $\eqref{eqn:FHN_compactified_slowtime}$ corresponding to compact subsets of $\tilde{\gamma}_{0,r}^{(\text{S})}\cap\pS_0^{(\ta,2)}$ and $\tilde{\gamma}_{0,r}^{(\text{S})}\cap\pS_0^{(\tr)}$ which connect to one another via the flow of \eqref{eqn:FHN_compactified_slowtime} near a \textit{folded saddle singularity} along $\pF^\text{(2)}$. In fact, there is a \textit{maximal canard trajectory} $\tilde{\gamma}_{\varepsilon,r}^\text{(S)}$ of \eqref{eqn:FHN_compactified_slowtime} which closely tracks $\pS_0^{(\ta,2)}$ by approximately traversing the entire path outlined by $\tilde{\gamma}_{0,r}^{(\text{S})}\cap\pS_0^{(\ta,2)}$, and upon crossing $\pF^{(2)}$, proceeds to closely track $\pS_0^{(\tr)}$ by approximately traversing the entire path outlined by $\tilde{\gamma}_{0,r}^{(\text{S})}\cap\pS_0^{(\tr)}$. This maximal canard trajectory $\tilde{\gamma}_{\varepsilon,r}^\text{(S)}$ ultimately moves ``downward'' through $\pF^{(1)}$ with $\text{d}u/\text{d}t=0$ and $\text{d}v/\text{d}t<0$, at a position which shall henceforth be denoted by 
\[\bs z^{(1)} \coloneqq (u^{(\tf,1)},v_{\Lambda(s^{(1)})}^{(\tf,1)},s^{(1)})\,,\quad \tilde{\gamma}_{\varepsilon,r}^\text{(S)}\cap\pF^{(1)}=\{\bz^{(1)}\}\,,\]
and afterwards returns directly to $\pS_0^{(\ta,2)}$ via an overwhelming increase in the $u$-variable. For examples of $\tilde{\gamma}_{\varepsilon,r}^\text{(S)}$ at the critical rate $r = r_\varepsilon^\text{(cr)}$, see $A_{\eps,r_\varepsilon^\text{(cr)}}(\be_-^{(\text{c})})$ within the middle row of Figure \ref{fig:nonautonomous_metastability_landscapes}. The position $\bz^{(1)}$ at which $\tilde{\gamma}_{\varepsilon,r_\varepsilon^\text{(cr)}}^\text{(S)}$ intersects the line of fold points $\pF^{(1)}$ is also indicated in the annotations of Figure \begin{NoHyper}\ref{fig:schematic_2d_quasithreshold}\end{NoHyper} of the supplementary material\cite{supplementary_material}.

In a close vicinity of the maximal canard trajectory $\tilde{\gamma}_{\varepsilon,r}^\text{(S)}$, the system exhibits sensitive divergent behaviour which shapes the structure of the quasithreshold $\pQ_{\eps,r}$. For initial conditions along $\pF^{(1)}$ in a sufficiently small neighbourhood around $\bs z^{(1)}$, the \textit{backward-time} trajectories of \eqref{eqn:FHN_compactified_slowtime} initially closely track $\pS_0^{(\tr)}$ by approximately traversing the entire path outlined by $\tilde{\gamma}_{0,r}^{(\text{S})}\cap\pS^{(\tr)}$, thereby also following the maximal canard trajectory. However, upon crossing $\pF^{(2)}$, and subsequently closely tracking $\pS_0^{(\ta,2)}$ for a transient period by approximately traversing an initial part of the path outlined by $\tilde{\gamma}_{0,r}^{(\text{S})}\cap\pS^{(\ta,2)}$, these trajectories are suddenly repelled away from $\pS_0^{(\ta,2)}$ via an overwhelming increase or decrease of the $u$-coordinate. Hence, these backward-time trajectories eventually peel away from $\tilde{\gamma}_{\varepsilon,r}^\text{(S)}$ somewhere along $\pS_0^\text{(a,2)}$ via an \textit{ensemble-splitting}-type effect, as depicted in panel (b) of Figure \ref{fig:schematic_2d_quasithreshold_whole} and all panels of Figure \ref{fig:nonautonomous_metastability_landscapes}. 

The closer the initial conditions on $\pF^{(1)}$ are to $\bs z^{(1)}$, the longer the backward-time trajectories of \eqref{eqn:FHN_compactified_slowtime} follow the maximal canard trajectory $\tilde{\gamma}_{\varepsilon,r}^{(S)}$ along $\pS_0^{(\ta,2)}$ for before being repelled away. For initial conditions $(u^{(\tf,1)},v_{\Lambda(s^\text{(i)})}^{(\tf,1)},s^\text{(i)})\in\pF^{(1)}$ with $0 < s^\text{(i)} - s^{(1)} \ll 1$, the backward-time trajectories are eventually repelled away from $\pS_0^{(\ta,2)}$ via an overwhelming \textit{increase} in the $u$-coordinate towards $u = +\infty$. These trajectories appear across Figure \ref{fig:nonautonomous_metastability_landscapes}, and partly make up the blue surface depicted in panel (a) of Figure \begin{NoHyper}\ref{fig:schematic_2d_quasithreshold}\end{NoHyper} of the supplementary material\cite{supplementary_material}. Conversely, for initial conditions $(u^{(\tf,1)},v_{\Lambda(s^\text{(i)})}^{(\tf,1)},s^\text{(i)})\in\pF^{(1)}$ with $0 < s^{(1)}- s^\text{(i)} \ll 1$, the backward-time trajectories are eventually repelled away from $\pS_0^{(\ta,2)}$ via an overwhelming \textit{decrease} in the $u$-coordinate towards $\pS_0^{(\tr)}$ (or $u = -\infty$, if $r > r_\varepsilon^\text{(cr)}$ and the backward-time trajectories follow $\tilde{\gamma}_{\varepsilon,r}^\text{(S)}$ long enough to reach $u \gtrsim 2\sqrt{3}/3$ - cf. Section \begin{NoHyper}\ref{sm:singular_supercritical_set}\end{NoHyper} of the supplementary material\cite{supplementary_material}). Upon reaching $\pS_0^\text{(r)}$, such backward-time trajectories move ``upward'' with $\text{d}u/\text{d}t=0$ and $\text{d}v/\text{d}t>0$, before closely tracking $\pS_0^{(\tr)}$ via the nearby backward-time slow manifolds which brings the system toward $\pF^{(2)}$ (cf. the left column of Figure \ref{fig:singular_metastability_landscapes}). Upon crossing $\pF^{(2)}$, these backward-time trajectories exhibit an overwhelming \textit{increase} in the $u$-coordinate towards $u = +\infty$. These trajectories are partly visible across Figure \ref{fig:nonautonomous_metastability_landscapes}, but for clarity, also appear in panel (b) of Figure \begin{NoHyper}\ref{fig:schematic_2d_quasithreshold}\end{NoHyper} of the supplementary material\cite{supplementary_material}. 

Altogether, this ensemble-splitting of trajectories of \eqref{eqn:FHN_compactified_slowtime} creates a folded, i.e., multilayered with respect to the $v$-coordinate, surface which separates the state space between regions which - under the flow of \eqref{eqn:FHN_compactified_slowtime} - connect in forward time to slow manifolds associated with parts of $\pS_0^{(\ta,2)}$ on either side of $\tilde{\gamma}_{0,r}^{(\text{S})}$. The slow manifolds associated with the part of $\pS_0^{(\ta,2)}$ coloured teal in Figure \ref{fig:singular_metastability_landscapes} induce supercritical behaviour, because they ultimately reach $\pF^{(2)}$ and connect to trajectories which ``shoot off'' towards $\pS_0^\text{(a,1)}$ on the fast timescale before settling to the equilibrium $\be_+^\text{(c)}$. On the other hand, the slow manifolds associated with the part of $\pS_0^\text{(a,2)}$ coloured purple in Figure \ref{fig:singular_metastability_landscapes} - which shall henceforth be denoted as $\pB_{0,r}(\{\be_+^{(\text{c})}\})$ - induce subcritical behaviour, because they account for trajectories which settle directly to the equilibrium $\be_+^\text{(c)}$ without ``shooting off'' towards $\pS_0^\text{(a,1)}$. Standard examples of trajectories corresponding to these supercritical and subcritical scenarios are illustrated across panels (a) - (b) of Figure \begin{NoHyper}\ref{fig:schematic_2d_quasithreshold}\end{NoHyper} of the supplementary material\cite{supplementary_material}.


An additional section of the quasithreshold $\pQ_{\eps,r}$ is made up of an ensemble of backward-time trajectories of \eqref{eqn:FHN_compactified_slowtime} initialised across the fold points within $\pF^{(1)}$ beyond a small neighbourhood of $\bs z^{(1)}$. Specifically, certain\footnote{In fact, all, whenever $r\leq r_\varepsilon^\text{(cr)}$ and $0<\varepsilon\ll 1$.} points along $\pF^{(1)}$ yield trajectories which: in backward time, connect to nearby slow manifolds which move along $\pS_0^{(\tr)}$; in forward time, are repelled away from $\pF^{(1)}$ towards $\pS_0^{(\ta,2)}$ via an overwhelming increase in the $u$-coordinate and subsequently connect to the slow manifolds associated with $\pB_{0,r}(\{\be_+^{(\text{c})}\})$. 
An example of such a trajectory is in panel (c) of Figure \begin{NoHyper}\ref{fig:schematic_2d_quasithreshold}\end{NoHyper} of the supplementary material\cite{supplementary_material}. 
The ensemble of the trajectories closely tracking $\pS_0^{(\tr)}$ in backward time creates a surface (see the portions of the panels of Figure \ref{fig:nonautonomous_metastability_landscapes} local to $\pF^{(1)}$) that separates the state space between regions which - under the flow of \eqref{eqn:FHN_compactified_slowtime} - exhibit supercritical and subcritical behaviour. ``Above'' the surface, the flow moves over $\pF^{(1)}$ and visits $\tilde{\pS}_\varepsilon$; ``below'' the surface, the system is repelled away from $\pF^{(1)}$ and connects to the slow manifolds associated  with $\pB_{0,r}(\{\be_+^{(\text{c})}\})$. Examples of trajectories corresponding to these supercritical and subcritical scenarios are illustrated across panels (a) - (c) of Figure \begin{NoHyper}\ref{fig:schematic_2d_quasithreshold}\end{NoHyper} of the supplementary material\cite{supplementary_material}. 


Part of this surface extends ``horizontally'' beyond $\pF^{(2)}$, via the backward-time trajectories which reach this set and directly undertake an overwhelming increase in the $u$-coordinate towards $u=+\infty$ (besides those that transiently follow the maximal canard trajectory $\tilde{\gamma}_{\varepsilon,r}^\text{(S)}$). A very narrow channel of subcritical behaviour is enclosed between supercritical behaviour by this horizontal extension and the folded surface resulting from the ensemble-splitting along the maximal canard trajectory $\tilde{\gamma}_{\varepsilon,r}^\text{(S)}$. This refers to the parts of the blue surfaces in panels (b) and (c) of Figure \begin{NoHyper}\ref{fig:schematic_2d_quasithreshold}\end{NoHyper} of the supplementary material\cite{supplementary_material} which lie above/below one another in the right panel of Figure \ref{fig:schematic_2d_quasithreshold_whole}. In the $\varepsilon\downarrow 0$ limit, these corresponding parts coalesce (cf. the right column panels of Figure \ref{fig:singular_metastability_landscapes}). An example of a subcritical trajectory from this narrow region is illustrated in panel (b) of Figure \begin{NoHyper}\ref{fig:schematic_2d_quasithreshold}\end{NoHyper} of the supplementary material\cite{supplementary_material}. This subcritical channel between the layers of the quasithreshold is too narrow to be apparent in Figure \ref{fig:nonautonomous_metastability_landscapes} under the displayed axis limits, but becomes clearly visible after restricting to the relevant section (not shown). 

The final component of the quasithreshold $\pQ_{\eps,r}$ is given by an ensemble of trajectories of \eqref{eqn:FHN_compactified_slowtime} which travel tangentially along $\tilde{\pS}_\varepsilon$ at the positions that separate its inflow and outflow boundaries, and directly connect to the slow manifolds associated with $\pB_{0,r}(\{\be_+^{(\text{c})}\})$ afterwards. This tangential grazing occurs close to the \textit{kink points} of the $L$-metastable set along $\{(u,v,s)\in\pD: u=u^{(\tf,1)},v=v_{\Lambda(s)}^{(\tf,1)}-\varepsilon^{2/3}\}$. As before, this creates a surface (see the portions of the panels of Figure \ref{fig:nonautonomous_metastability_landscapes} below $\tilde{\pS}_\varepsilon$) which separates the state space between regions that, under the flow of \eqref{eqn:FHN_compactified_slowtime}, exhibit supercritical and subcritical behaviour. Examples of corresponding supercritical and subcritical trajectories are illustrated in panel (d) of Figure \begin{NoHyper}\ref{fig:schematic_2d_quasithreshold}\end{NoHyper} of the supplementary material\cite{supplementary_material}.  

Local to the attracting sheet $\pS_0^\text{(a,2)}$, the surface layers of the quasithreshold are built around the position of the maximal canard trajectory $\tilde{\gamma}_{\varepsilon,r}^\text{(S)}$, which is expected to coincide exactly once with the pullback attractor $A_{\varepsilon,r}(\be_-^\text{(c)})$ in between stages of $\tilde{\gamma}_{\varepsilon,r}^\text{(S)}$ lying locally ``above'' and then ``below'' $A_{\varepsilon,r}(\be_-^\text{(c)})$. This is based on the behaviour of their limiting counterparts, the folded saddle singular canard trajectory and the singular pullback attractor (recall Section \ref{subsec:geometric_singular_perturbation_theory}). Consequently, we conjecture that there exists a unique critical rate $r=r_\varepsilon^\text{(cr)}\in\bbR_+$ which separates subcritical and supercritical behaviour of the pullback attractor in the following sense:
\begin{equation*}
    \begin{cases}
        A_{\varepsilon,r}(\be_-^\text{(c)})\subset \pD_{\varepsilon,r}^\text{(sub)} & \text{if } r\in(0,r_\varepsilon^\text{(cr)}) \\
        A_{\varepsilon,r}(\be_-^\text{(c)})\cap\pQ_{\varepsilon,r}\neq\emptyset \text{\,\,\,and\,\,\,} A_{\varepsilon,r}(\be_-^\text{(c)})\subset \tilde{\pB}\setminus\pD_{\varepsilon,r}^\text{(sup)} & \text{if } r=r_\varepsilon^\text{(cr)} \\ 
        A_{\varepsilon,r}(\be_-^\text{(c)})\cap \pD_{\varepsilon,r}^\text{(sup)} \neq \emptyset & \text{if } r\in(r_\varepsilon^\text{(cr)},\infty)
    \end{cases}\,.
\end{equation*}
This aligns with the conclusions drawn and results found for the simple slow case of ``tracking-tipping transition'' studied in \citet{osullivan2023} - see panel a) of Figure 8 in \cite{osullivan2023} - and is also corroborated by numerical experiments conducted here for system \eqref{eqn:FHN_compactified_slowtime}. Indeed, we find that $r_{\eps=0.01}^\text{(cr)} \approx 0.560606$ and $r_{\eps=0.1}^\text{(cr)} \approx 2.287962$. Recall that the critical rate for the singular case was $r_0^\text{(cr)}\approx 0.5277897$. Away from $\varepsilon = 0$, there is a notable discrepancy between $r_0^{(\text{crit})}$ and $r_\varepsilon^{(\text{crit})}$, though it will hold that $\displaystyle\lim_{\eps\downarrow 0}r_\varepsilon^{(\text{crit})}=r_0^{(\text{crit})}$.

As alluded to at the beginning of this section, the quasithresholds $\pQ_{\eps,r}$ can be approximated numerically for system \eqref{eqn:FHN_compactified_slowtime}, by methods which are briefly outlined in Section \begin{NoHyper}\ref{sm:quasithreshold_numerical_approximation}\end{NoHyper} of the supplementary material\cite{supplementary_material}. These approximations are gathered in full in Figure \ref{fig:nonautonomous_metastability_landscapes} for two timescale separation values, $\varepsilon = 0.01$ and $\varepsilon = 0.1$; for each case, rates $r$ equal to, just below, and just above $r_\varepsilon^{(\text{crit})}$ are considered. The structure and positioning of $\pQ_{\varepsilon=0.01,r}$ closely aligns with that of $\pQ_{\varepsilon=0,r}$ (cf. the right panels of Figure \ref{fig:singular_metastability_landscapes}), for subcritical, critical and supercritical $r$ values. On the other hand, there is notable discrepancies between the positioning of $\pQ_{\eps=0.1,r}$ and that of $\pQ_{\varepsilon=0,r}$. This is related to a combination of the facts that: the singular limit of system \eqref{eqn:FHN_compactified_slowtime} with $\varepsilon = 0$ is better approximated by \eqref{eqn:FHN_compactified_slowtime} with $\varepsilon = 0.01$ than it is with $\varepsilon = 0.1$; and the singular critical rate $r_0^\text{(cr)}$ is considerably closer to $r_\varepsilon^\text{(cr)}$  for $\varepsilon = 0.01$ than it is for $\varepsilon = 0.1$.

\section{A stochastic autonomous FitzHugh-Nagumo system}\label{sec:stochastic_FHN_autonomous}

Another forcing mechanism for inducing supercritical excursions (starting from a position of equilibrium) is the incorporation of stochasticity into system \eqref{eqn:FHN_autonomous_slowtime}. For instance, rather than having the input current evolve deterministically throughout time as in Section \ref{sec:fitzhughnagumo_nonautonomous}, it can be set to the value of a noise process which evolves unpredictably throughout time. This introduces fluctuations in the evolution of the system relative to the dynamics of \eqref{eqn:FHN_autonomous_slowtime}, which have a different effect at each instance of time, resulting in unpredictable deviations away from the deterministic flow lines. In particular, this can cause the system to be driven out of a set which is forward-time invariant under the deterministic flow. This becomes relevant for supercritical excursions if the stochastic forcing were to drive the system out of the subcritical set $\pD_{\eps,I}^{(\text{sub})}$ into the supercritical set $\pD_{\eps,I}^{(\text{sup})}$, since then the system may subsequently travel to the $L$-metastable set $\tilde{\pS}_{\eps,I}$ via a rough shadowing of the deterministic dynamics. Such scenarios are most generally referred to as noise-induced escapes or noise-induced transitions\cite{horsthemke_noise-induced_2006}. We shall now document - throughout this section - the likelihood with which noise-induced transitions to $\tilde{\pS}_\varepsilon$ occur (across finite-time intervals) within the autonomous FitzHugh-Nagumo system subjected to a stochastically evolving input current.

A simple way to feed stochasticity into the input current $I$ is to augment its value with a white-in-time Gaussian noise process. This motivates considering the dynamical system yielded by the following stochastic ordinary differential equation initial value problem of It\^{o}-type,
\begin{equation}\label{eqn:FHN_autonomous_slowtime_stochastic}
    \begin{cases}
        \text{d} u = \eps^{-1}\big(-u^3 + u - v + I\big)\,\text{d}t+\sigma\eps^{-1/2}\,\text{d}W \vspace{0.2cm}\\ 
        \text{d} v = \big(-\beta v + u\big)\,\text{d}t
    \end{cases},\,\begin{pmatrix}
        u(0) \\ v(0)
    \end{pmatrix} = \begin{pmatrix}
        u^{(\ti)} \\ v^{(\ti)}
    \end{pmatrix}\,, 
\end{equation} 
for $(u(t),v(t))\in\bbR^2,\,t\in\bbR_{\geq 0}$. Relative to the underlying deterministic dynamics of \eqref{eqn:FHN_autonomous_slowtime}, the evolution equation of the $u$-variable has been augmented with an \textit{additive} white-in-time Gaussian noise process. This is provided by the increments of the Wiener process $W:[0,\infty)\rightarrow\mathbb{R}$ that is defined on a filtered probability space\cite{karatzas_shreve} with sample space $\Omega$ and probability measure $\rho_\Omega$. Here, $0<\sigma\ll 1$ is a positive scalar parameter which scales the amplitude of the noise process and shall henceforth be referred to as the noise strength. To maintain balance between the deterministic and stochastic components of the system in the theoretical limit $\varepsilon\downarrow 0$, the noise increments are scaled by $\varepsilon^{-1/2}$. This is the appropriate scaling to ensure that the stochasticity on $u$ becomes neither asymptotically negligible nor asymptotically dominant - relative to its deterministic component - as the timescale separation parameter $\varepsilon$ limits to zero. This becomes apparent when applying a rescaling of time to consider the dynamics of \eqref{eqn:FHN_autonomous_slowtime_stochastic} on the fast time $\tau = t/\varepsilon$, which reads as  
\begin{equation}\label{eqn:FHN_autonomous_drift_fasttime_stochastic}
    \begin{cases}
        \text{d} \tilde{u} = \big(-\tilde{u}^3 + \tilde{u} - \tilde{v} + I\big)\,\text{d}\tau+\sigma\,\text{d}\tilde{W} \vspace{0.2cm}\\ 
        \text{d} \tilde{v} = \eps\big(-\beta \tilde{v} + \tilde{u}\big)\,\text{d}\tau
    \end{cases},\,\begin{pmatrix}
        \tilde{u}(0) \\ \tilde{v}(0)
    \end{pmatrix} = \begin{pmatrix}
        u^{(\ti)} \\ v^{(\ti)}
    \end{pmatrix}\,,
\end{equation} 
where $\tilde{W}(\tau)\coloneqq\varepsilon^{-1/2}W(\varepsilon\tau)$ describes a (rescaled) Wiener process, for all $\varepsilon\in\bbR_+$.

We write the solution of the stochastic initial value problem \eqref{eqn:FHN_autonomous_slowtime_stochastic} using the flow notation
$$
\bs\phi_{\sigma,\eps,I}^{(\text{FHN})}((u^{(\ti)},v^{(\ti)}),W,t)=(u(t),v(t))\,,\,\,\,t\geq 0\,,
$$
and the \textit{diffusion matrix} of \eqref{eqn:FHN_autonomous_slowtime_stochastic} using the matrix notation
\begin{equation*}
    \bs\Gamma_\varepsilon^{(\text{FHN})}(u,v) \equiv \bs\Gamma_\varepsilon^{(\text{FHN})} \coloneqq \begin{pmatrix}
        \eps^{-1/2} & 0 \\ 0 & 0
    \end{pmatrix}\,.
\end{equation*}
Note the dependence of $\bs\phi_{\sigma,\eps,I}^{(\text{FHN})}((u^{(\ti)},v^{(\ti)}),W,t)$ on the particular Wiener process $W$. 

\subsection{The occurrence and consequences of noise-induced supercritical excursions}\label{subsec:stochastic_autonomous_FHN_occurrence_consequences}

The first goal of this section to introduce how prevalent the occurrence of noise-induced supercritical excursions is within \eqref{eqn:FHN_autonomous_slowtime_stochastic}. Consider a sample path $\bs\phi_{\sigma,\eps,I}^{(\text{FHN})}(\be_I^\text{(FHN)},W,t)$ of \eqref{eqn:FHN_autonomous_slowtime_stochastic} initialised at the equilibrium point $\be_I^\text{(FHN)}$ and driven by the Wiener process $W\in\Omega$. Define the \textit{first-entrance time} of this sample path into the $L$-metastable set $\tilde{\pS}_{\varepsilon,I}$ as
\begin{flalign*}
    \zeta_{\sigma,\eps,I}^\text{(FHN)}(\be_I^\text{(FHN)},W,\tilde{\pS}_{\varepsilon,I}) 
    \coloneqq &\inf\{t>0:\bs\phi_{\sigma,\eps,I}^{(\text{FHN})}(\be_I^\text{(FHN)},W,t)\in\tilde{\pS}_{\varepsilon,I}\} 
\end{flalign*}
It will hold that $\zeta_{\sigma,\eps,I}^\text{(FHN)}\big(\be_I^{(\text{FHN})},W,\tilde{\pS}_{\eps,I}\big)<\infty$ almost surely, i.e. sample paths of \eqref{eqn:FHN_autonomous_slowtime_stochastic} with $(u^{(\ti)},v^{(\ti)})=(u_I^\text{(eq)},v_I^\text{(eq)})$ enter $\tilde{\pS}_{\eps,I}$ after (unbounded) finite-time with probability one (cf. \citet{dembozeitouni},  Theorem $5.7.11(a)$). Given some finite $T\in\bbR_+$, we say that the sample path $\bs\phi_{\sigma,\eps,I}^\text{(FHN)}(\be_I^{(\text{FHN})},W,t)$ undertakes a \textit{noise-induced supercritical excursion} within $T$ time units if 
$$
\zeta_{\sigma,\eps,I}^\text{(FHN)}(\be_I^{(\text{FHN})},W,\tilde{\pS}_{\eps,I})< T. 
$$
The probability that a sample path undertakes a noise-induced supercritical excursion within $T$ time units is, however, less than one; in fact, this probability will tend to zero in the weak-noise limit $\sigma\rightarrow 0$ (cf. \citet{dembozeitouni}, Exercise 5.6.25). 

\begin{figure}[ht]
    \centering
    \includegraphics[width=6.48in]{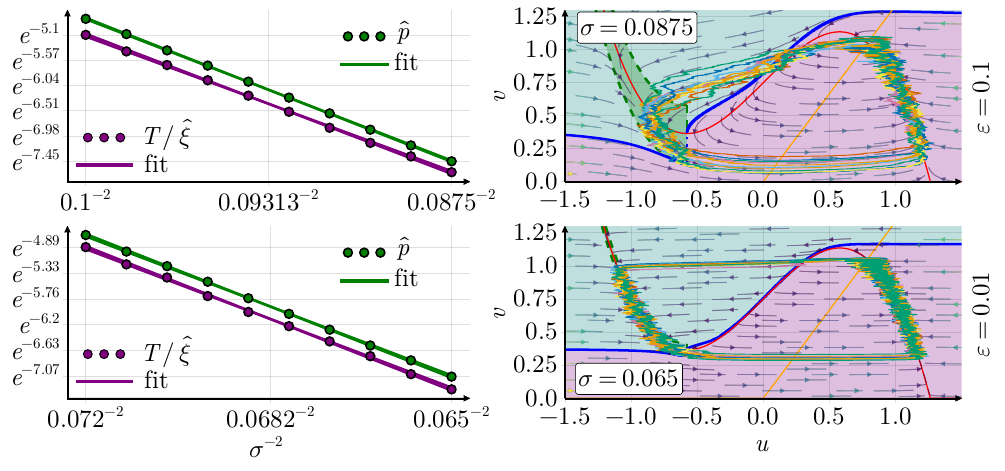}
    \caption{The vulnerability to noise-induced supercritical excursions in the stochastic autonomous FitzHugh-Nagumo system \eqref{eqn:FHN_autonomous_slowtime_stochastic}. Left: numerical estimates $\hat{p}=\hat{p}_{\sigma,\varepsilon,I,T}^\text{(sup)}(\be_I^\text{(FHN)})$ of $p_{\sigma,\varepsilon,I,T}^\text{(sup)}(\be_I^\text{(FHN)})$ [green points] and scaled reciprocal values of numerical estimates $\hat{\xi}=\hat{\xi}_{\sigma,\eps,I}^{(\text{sup})}$ of $\xi_{\sigma,\eps,I}^{(\text{sup})}$ [purple points], plotted using a logarithmic scale on the vertical axis. Estimates were generated via the adaptive multilevel splitting algorithm. The green and purple intersecting curves describe weighted least-squares exponential fits to $\hat{p}$ and $\hat{\xi}$ (before rescaling), respectively. The resulting \textit{slopes} produce the estimates $(\Delta \hat{V}_{\eps=0.1,I,T}^\text{(sup)},\Delta \hat{V}_{\eps=0.01,I,T}^\text{(sup)})\coloneqq(0.08654,0.05421)$ and $(\Delta \hat{V}_{\eps=0.1,I}^\text{(sup)},\Delta \hat{V}_{\eps=0.01,I}^\text{(sup)})\coloneqq(0.08332,0.05422)$ for $\Delta V_{\eps,I,T}^\text{(sup)}$ and $\Delta V_{\eps,I}^\text{(sup)}$ respectively (for weights and associated confidence intervals, see Tables \ref{tab:autonomous_AMS_statistical_output} - \ref{tab:confidence_intervals_autonomous_system} of the supplementary material\cite{supplementary_material}). Right: an ensemble of ten sample paths $\bs\phi_{\sigma,\eps,I}^\text{(FHN)}(\be_I^{(\text{FHN})},W,t)$ of system $\eqref{eqn:FHN_autonomous_slowtime_stochastic}$ satisfying $\zeta_{\sigma,\eps,I}^\text{(FHN)}(\be_I^{(\text{FHN})},W,\tilde{\pS}_{\eps,I})<T$, for $(\varepsilon,\sigma) = (0.1,0.0875)$ and $(\varepsilon,\sigma) = (0.01,0.065)$ [coloured noisy paths]. The sample paths are plotted across the time-interval $t\in [t^{(\text{reactive})},t^{(\text{return})}]$ constructed with respect to $\tilde{\pS}_{\eps,I}$ and the neighbourhood $\pN_{\eps,\delta}(\be_I^{(\text{FHN})})$ with $\delta = 0.1$ (see Eqs. \eqref{eqn:L-metastable_set_unforced}, \eqref{eqn:reactive_time}-\eqref{eqn:FHN_autonomous_neighbourhood}). The backgrounds of these panels display the metastability landscape of the underlying deterministic system. The parameter values considered throughout are $I=0.75$ and $T = 50.0$.}    \label{fig:excitable_sample_paths_autonomous}
\end{figure}

The rate at which this probability tends to zero in the weak noise limit is discussed and evaluated in the following section \ref{subsec:probability_autonomous}, but prior to that, we briefly comment here on how noise-induced supercritical excursions lead to \textit{sustained} continuous periods where the system is in a different state away from equilibrium. To this end, begin by noting that whenever the first-entrance time $\zeta_{\sigma,\eps,I}^\text{(FHN)}=\zeta_{\sigma,\eps,I}^\text{(FHN)}\big(\be_I^{(\text{FHN})},W,\tilde{\pS}_{\eps,I}\big)$ is finite, the corresponding sample path must lie on the boundary of $\tilde{\pS}_{\varepsilon,I}$ at this time $t = \zeta_{\sigma,\varepsilon,I}^\text{(FHN)}$, i.e.  
\[ 
\bs\phi_{\sigma,\eps,I}^\text{(FHN)}(\be_I^{(\text{FHN})},W,\zeta_{\sigma,\eps,I}^\text{(FHN)})\in\partial\tilde{\pS}_{\eps,I}\,.
\]
Directly after $t = \zeta_{\sigma,\varepsilon,I}^\text{(FHN)}$, the sample path enters $\tilde{\pS}_{\eps,I}$ for the first time (almost surely). While the deterministic inflow and outflow boundaries of the $L$-metastable set $\tilde{\pS}_{\eps,I}$,
or the residence times of trajectories (sample paths) within $\tilde{\pS}_{\varepsilon,I}$, do not exactly persist in the context of system \eqref{eqn:FHN_autonomous_slowtime_stochastic} due to the stochastic deviations away from the solutions of \eqref{eqn:FHN_autonomous_slowtime}, they will instead closely persist in a probabilistic sense under weak noise strengths\cite{vanden2006towards,vanden2010transition}. For instance, the discussion in Section \ref{subsec:almost_invariant} outlines that - under sufficiently weak noise - the stochastic system \eqref{eqn:FHN_autonomous_slowtime_stochastic} will lose mass from $\tilde{\pS}_{\varepsilon,I}$ at a similar rate to the loss of mass from $\tilde{\pS}_{\varepsilon,I}$ in the deterministic system \eqref{eqn:FHN_autonomous_slowtime}. Further, the last-exit position of the sample path from the set $\tilde{\pS}_{\eps,I}$ before returning to the attracting branch $\pS_0^\text{(a,1)}$ is very likely to be contained in the outflow boundary of $\tilde{\pS}_{\varepsilon,I}$ as a result of the asymptotic agreement result in \eqref{eqn:asymptotic_agreement}. On passing through the outflow boundary, the system is then likely to travel to a close neighbourhood of $\be_I^{(\text{FHN})}$ without any large deviations away from the deterministic return mechanism. See Figure \ref{fig:excitable_sample_paths_autonomous} for an ensemble of numerically simulated noise-induced supercritical excursions, which roughly corroborate these latter statements. 

In order to formally deduce the time frames across which $\tilde{\pS}_{\varepsilon,I}$ is an almost-invariant set, one would have to inspect the invariant densities of system \eqref{eqn:FHN_autonomous_slowtime_stochastic}, though numerical observations (the sample paths in Figure \ref{fig:excitable_sample_paths_autonomous}) suggest that the set $\tilde{\pS}_{\varepsilon,I}$ loses mass in the stochastic invariant regime at rates $L_\sigma\approx 27$ for $(\varepsilon,\sigma)=(0.1,0.0875)$ and $L_\sigma\approx 603$ for $(\varepsilon,\sigma)=(0.01,0.065)$. This being said, Figure \ref{fig:excitable_sample_paths_autonomous} actually showcases that the noise-induced supercritical excursions found under these noise strengths are consistently being kicked in and out of $\tilde{\pS}_\varepsilon$ between its initial and final occupation times. Thus, while the noise strengths are weak enough to make the probability of noise-induced supercritical excursions within finite time already very low, at the same time they are strong enough for the corresponding sample paths to consistently diffuse beyond the $\varepsilon^{2/3}$-radius channel of the critical manifold which contains the deterministic slow manifolds. With weaker noise strengths, it \textit{will} become highly probable for the sample paths showcasing noise-induced supercritical excursions to consistently remain within this channel, and as a result the $L_\sigma$-value (that is inversely proportional to characteristic residence times) is expected to decrease. It may be possible to prove that the rate of loss of mass from $\tilde{\pS}_{\varepsilon,I}$ in the stochastic invariant regime satisfies $L_\sigma\overset{\varepsilon\downarrow 0}{=}\pO(1)$, with certain noise strengths $\sigma$ that depend on the timescale separation parameter $\varepsilon$. 

For sharper estimates of related stochastic residence times, upper and lower bounds for the probability that the first-exit time from a parameterised neighbourhood of an attracting slow manifold is less than some finite positive value are given in \citet{berglundgentz2006} and \citet{kuehn2015}. Put simply, these demonstrate that sample paths are likely to remain within a neighbourhood of an attracting slow manifold on exponential time scales with respect to $\sigma^{-2}$ (unless the outflow boundary is reached).

\subsection{The probability of noise-induced supercritical excursions within finite-time}\label{subsec:probability_autonomous} 

To quantify the decay in the probability with which sample paths $\bs\phi_{\sigma,\eps,I}^\text{(FHN)}(\be_I^{(\text{FHN})},W,t)$ of \eqref{eqn:FHN_autonomous_slowtime_stochastic} undertake noise-induced supercritical excursions within a given finite $T$ time units, in the weak-noise limit $\sigma\rightarrow 0$, one can apply the large deviation principle (LDP) developed within sample path large deviation theory\cite{freidlinwentzell,dembozeitouni}. Note that the event $\zeta_{\sigma,\eps,I}^\text{(FHN)}(\be_I^{(\text{FHN})},W,\tilde{\pS}_{\eps,I})<T$ corresponds to the sample path $\bs\phi_{\sigma,\eps,I}^\text{(FHN)}(\be_I^{(\text{FHN})},W,t)$ being contained within the set of (continuous) paths which start at $\be_I^{(\text{FHN})}$ and enter $\tilde{\pS}_{\eps,I}$ within $T$ time units, namely   
\begin{equation*}
    \pP_{\eps,I,T}^{(\text{sup})} \coloneqq \big\{\bs\psi\in C^0([0,T],\bbR^2):\bs\psi(0)=\be_I^{(\text{FHN})},\,\inf\{t\in[0,T]:\bs\psi(t)\in\tilde{\pS}_{\eps,I}\}<T\big\}
\end{equation*}
The aforementioned large deviation principle (Theorem 5.6.7 of \citet{dembozeitouni}, see also \citet{freidlinwentzell} and \citet{berglundgentz2006}) then quantifies the scaling of the probability of these noise-induced supercritical excursions in the weak noise limit $\sigma\rightarrow 0$ as
\begin{equation}\label{eqn:asymptotic_probability_scaling}
    p_{\sigma,\eps,I,T}^\text{(sup)}(\be_I^\text{(FHN)})\coloneqq\bbP_{W\sim\rho_\Omega}\big(\zeta_{\sigma,\eps,I}^\text{(FHN)}(\be_I^{(\text{FHN})},W,\tilde{\pS}_{\eps,I})<T\big) \overset{\sigma\downarrow 0}{\asymp}\text{exp}\bigg(-\sigma^{-2}\displaystyle\inf_{\bs\psi\in\pP_{\eps,I,T}^{(\text{sup})}}S_T[\bs\psi]\bigg)\,.
\end{equation}
The symbol $\asymp$ denotes logarithmic asymptotic equivalence, meaning that the quotient of the logarithms of the right- and left-hand sides of \eqref{eqn:asymptotic_probability_scaling} converges to 1 as $\sigma\downarrow 0$. Here the \textit{rate function} $S_T$ represents the \textit{Freidlin-Wentzell action functional}, which in the (singular noise) case of system \eqref{eqn:FHN_autonomous_slowtime_stochastic} is given for $\bs\psi(t)\in C^0([0,T],\bbR^2)$ in terms of the following variational principle
\begin{gather}\label{eqn:FW_action_functional}
S_T[\bs\psi]=\inf_{\bs\vartheta\in\Theta(\bs\varphi)}\frac{1}{2}\int_0^T\big\lVert \bs\vartheta(s)\big\rVert^2\,\text{d}s\,.
\end{gather}
Here, $\Theta(\bs\psi)$ is the set of $\bs\vartheta\in L_2([0,T],\bbR^2)$ such that 
\begin{gather*}
\bs\psi(t)=\bx^{(\ti)} + \int_0^t\bs b_{\varepsilon,I}^{(\text{FHN})}(\bs\psi(s))+\bs\Gamma_\varepsilon^{(\text{FHN})}\,\bs\vartheta(s)\,\text{d}s,\,\forall \,t\in[0,T]\,.
\end{gather*}
In physics terminology, the Freidlin-Wentzell action functional $S_T$ measures the ``energy cost'', ``work'' or ``action'' required by the stochastic forcing (represented by $\bs\vartheta$
) to realise the path $\bs\psi(t)$ in system $\eqref{eqn:FHN_autonomous_slowtime_stochastic}$. The large deviation principle \eqref{eqn:asymptotic_probability_scaling} captures the fact that, in the weak noise limit, noise-induced supercritical excursions occur almost surely via the mechanism(s) with least action (see Eq. 6 in \citet{grafke2019}), i.e., there exists a ``most probable'' transition path(s) which asymptotically dominates the support of the corresponding transition events. The probability itself decays \textit{exponentially} with $\sigma^{-2}$ at a rate given by the least action value.    

It is useful to quantify the asymptotic \textit{rate} 
\begin{equation}\label{eqn:autonomous_asymptotic_rate_finite_time}
\Delta V_{\eps,I,T}^{(\text{sup})}\coloneqq \inf_{\bs\psi\in\pP_{\eps,I,T}^{(\text{sup})}}S_T[\bs\psi]> 0\,,
\end{equation}
with which $p_{\sigma,\eps,I,T}^\text{(sup)}(\be_I^\text{(FHN)})$ decays exponentially to 0 with $\sigma^{-2}$, as a means of assessing how the relative risk of noise-induced supercritical excursions changes across different system parameterisations. To elaborate on this point, since \eqref{eqn:asymptotic_probability_scaling} also yields a large deviation principle for the \textit{ratio} of probabilities, quantifying the asymptotic rates in \eqref{eqn:autonomous_asymptotic_rate_finite_time} enables one to deduce whether the likelihood of system $\eqref{eqn:FHN_autonomous_slowtime_stochastic}$ undertaking noise-induced supercritical excursions either increases or decreases as, say, the timescale separation strengthens. For each $T^{(1)},\,T^{(2)}\in\bbR_+$, with $T^{(1)}<T^{(2)}$, and $\bs\psi_{T^{(1)}}\in\pP_{\eps,I,T^{(1)}}^\text{(sup)}$, one has that the candidate path defined by 
$$
\bs\psi_{T^{(2)}}(t) \coloneqq \begin{cases}
    \be_I^\text{(FHN)} & \text{if }\,0\leq t\leq T^{(2)}-T^{(1)} \\ 
    \bs\psi_{T^{(1)}}(t-(T^{(2)}-T^{(1)})) & \text{if }\,T^{(2)}-T^{(1)}<t\leq T^{(2)}
\end{cases}
$$
satisfies $\bs\psi_{T^{(2)}}\in\pP_{\eps,I,T^{(2)}}^\text{(sup)}$ with $S_{T^{(1)}}[\bs\psi_{T^{(1)}}]=S_{T^{(2)}}[\bs\psi_{T^{(2)}}]$. This can be used to show that
$$ 
\Delta V_{\eps,I,T^{(2)}}^\text{(sup)}\leq\Delta V_{\eps,I,T^{(1)}}^\text{(sup)}\,,
$$
i.e. the asymptotic rate in \eqref{eqn:autonomous_asymptotic_rate_finite_time} is a monotonically decreasing function of $T$. Further, the asymptotic rate limits to a finite positive value as $T\rightarrow\infty$, namely the quantity
\begin{equation}\label{eqn:autonomous_asymptotic_rate}
\Delta V_{\eps,I}^{(\text{sup})}\coloneqq \inf_{T\in(0,\infty)}\Delta V_{\eps,I,T}^{(\text{sup})}=\inf_{T\in(0,\infty)}\inf_{\bs\psi\in\pP_{\eps,I,T}^{(\text{sup})}}S_T[\bs\psi]\,,
\end{equation}
which is an example of a \textit{quasipotential barrier height}\cite{ashwin2023}. The quasipotential barrier height is instructive in itself as it provides\cite{dembozeitouni} the rate at which the \textit{expected value} of $\zeta_{\sigma,\eps,I}^\text{(FHN)}(\be_I^\text{(FHN)},W,\tilde{\pS}_{\eps,I})$ diverges to infinity in the weak-noise limit, through the asymptotic logarithmic equivalence   
$$
\xi_{\sigma,\eps,I}^\text{(sup)}(\be_I^\text{(FHN)})\coloneqq \bbE_{W\sim\rho_\Omega}\big[\zeta_{\sigma,\eps,I}^\text{(FHN)}(\be_I^\text{(FHN)},W,\tilde{\pS}_{\eps,I})\big]\overset{\sigma\downarrow 0}{\asymp}\exp\big(\sigma^{-2}\Delta V_{\eps,I}^\text{(sup)}\big)\,.
$$
The values of $\Delta V_{\eps,I,T}^\text{(sup)}$ and $\Delta V_{\eps,I}^\text{(sup)}$ can be approximated numerically through a series of estimates for $p_{\sigma,\eps,I,T}^\text{(sup)}(\be_I^\text{(FHN)})$ and $\xi_{\sigma,\eps,I}^\text{(sup)}(\be_I^\text{(FHN)})$ across a range of weak-noise values, followed by a fitting of the estimates to extract the rates of exponential scaling. The probabilities $p_{\sigma,\eps,I,T}^\text{(sup)}(\be_I^\text{(FHN)})$ are estimated via the trajectory-adaptive multilevel splitting algorithm (introduced in Section \begin{NoHyper}\ref{sm:AMS}\end{NoHyper} of the supplementary material\cite{supplementary_material}), while the mean first-exit times $\xi_{\sigma,\eps,I}^\text{(sup)}(\be_I^\text{(FHN)})$ are estimated via direct Monte-Carlo sampling methods. The estimates are gathered in Figure \ref{fig:excitable_sample_paths_autonomous} for timescale separation values $\varepsilon = 0.1$ and $\varepsilon = 0.01$ and input current $I = 0.75$. The outcome of the fits (the procedure behind which is introduced in Section \begin{NoHyper}\ref{sm:wlsqfits}\end{NoHyper} of the supplementary material\cite{supplementary_material}) shows, in particular, that the probability of undertaking noise-induced supercritical excursions within $T$ time units in system \eqref{eqn:FHN_autonomous_slowtime_stochastic} decays exponentially at a rate $\Delta V_{\eps,I,T}^\text{(sup)}\overset{\text{o.o.m}}{\sim}10^{-2}$, for every $T\in[50,\infty)$. In addition, since 
\begin{equation}\label{eqn:comparative_asymptotic_rates_autonomous}
\Delta \hat{V}_{\eps=0.01,I=0.75}^\text{(sup)}<\Delta \hat{V}_{\eps=0.1,I=0.75}^\text{(sup)}\,,
\end{equation}
we learn that \eqref{eqn:FHN_autonomous_slowtime_stochastic} becomes more vulnerable to noise-induced supercritical excursions within $T$ time-units with a stronger timescale separation $\eps = 0.01$ than with a weaker timescale separation $\eps = 0.1$, for again every $T\in [50,\infty)$. 

Figure \ref{fig:excitable_sample_paths_autonomous} also illustrates ensembles of sample paths $\bs\phi_{\sigma,\eps,I}^\text{(FHN)}(\be_I^{(\text{FHN})},W,t)$ of \eqref{eqn:FHN_autonomous_slowtime_stochastic} which enter the set $\tilde{\pS}_{\eps,I}$ within the finite-time interval $t \in(0,50)$, for timescale separation values $\varepsilon = 0.1$ and $\varepsilon=0.01$ and weak noise strengths corresponding to the limit of what is practically attainable. These were generated via Monte Carlo rejection sampling using the semi-implicit Euler-Maruyama method\cite{higham2002}, which is briefly outlined in Section \begin{NoHyper}\ref{sm:timestepping}\end{NoHyper} of the supplementary material\cite{supplementary_material}. The sample paths are shown between the last time the system is in a neighbourhood $\pN_{\eps,\delta}(\be_I^{(\text{FHN})})$ around the equilibrium $\be_I^{(\text{FHN})}$ \textit{before} entering $\tilde{\pS}_{\eps,I}$ at $t = \zeta_{\sigma,\eps,I}^\text{(FHN)}(\be_I^{(\text{FHN})},W,\tilde{\pS}_{\eps,I})$, and the first time the system re-enters the neighbourhood $\pN_{\eps,\delta}(\be_I^{(\text{FHN})})$ \textit{after} entering $\tilde{\pS}_{\eps,I}$ at $t = \zeta_{\sigma,\eps,I}^\text{(FHN)}(\be_I^{(\text{FHN})},W,\tilde{\pS}_{\eps,I})$. These times correspond to
\begin{flalign}
    t^\text{(reactive)} &\coloneqq \sup\big\{0<t<\zeta_{\sigma,\eps,I}^\text{(FHN)}(\be_I^{(\text{FHN})},W,\tilde{\pS}_{\eps,I}):\bs\phi_{\sigma,\eps,I}^\text{(FHN)}(\be_I^{(\text{FHN})},W,t)\in\pN_{\eps,\delta}(\be_I^{(\text{FHN})})\big\}\,, \label{eqn:reactive_time} \\ 
    t^\text{(return)} &\coloneqq \inf\big\{t>\zeta_{\sigma,\eps,I}^\text{(FHN)}(\be_I^{(\text{FHN})},W,\tilde{\pS}_{\eps,I}):\bs\phi_{\sigma,\eps,I}^\text{(FHN)}(\be_I^{(\text{FHN})},W,t)\in\pN_{\eps,\delta}(\be_I^{(\text{FHN})})\big\}\,, \label{eqn:return_time}
\end{flalign}
respectively. For the neighbourhood around the equilibrium, we choose the ellipse
\begin{equation}\label{eqn:FHN_autonomous_neighbourhood}
    \pN_{\eps,\delta}(\be_I^{(\text{FHN})}) = \{(u,v)\in\mathbb{R}^2:(u-u_I^\text{(eq)})^2+(v-v_I^\text{(eq)})^2/\eps^2<\delta^2\}\,,
\end{equation}
for $\delta>0$ sufficiently small so as to satisfy $\pN_{\eps,\delta}(\be_I^{(\text{FHN})})\subset \pD_{\varepsilon,I}^\text{(sub)}$. Here the dependence on $\eps$ reflects the fact that, in the lead-up to escaping the subcritical set, the standard deviation away from the equilibrium $\be_I^{(\text{FHN})}$ that is undertaken by the sample paths in the $v$-coordinate becomes smaller as $\eps$ decreases. This itself suggests that $\displaystyle\lim_{\eps\downarrow 0}\Delta V_{\eps,I}^\text{(sup)}$ may align with the quasipotential barrier height of the reduced one-dimensional system for $u(t)$ when taking $\varepsilon=0$ and $\tilde{v}\equiv v_I^\text{(eq)}$ in \eqref{eqn:FHN_autonomous_drift_fasttime_stochastic}. This refers to the quantity 
\[
-2\displaystyle\int_{u_I^\text{(eq)}}^{w_I} -\tilde{u}^3+\tilde{u}-v_I^\text{(eq)}+I\,\text{d}\tilde{u}\,,
\]
where $w_I$ satisfies $h_I(w_I)=v_I^\text{(eq)}$ and $u^\text{(f,1)}<w_I<u^\text{(f,2)}$, which to six significant digits equals 0.0509187 for the case $I=0.75$.  

\section{A stochastic nonautonomous FitzHugh-Nagumo system}\label{sec:stochastic_FHN_nonautonomous}

Adding stochastic forcing to the rate-dependent system \eqref{eqn:FHN_compactified_slowtime} causes trajectories to deviate from their deterministic counterparts, and these deviations can accumulate to either \textit{induce} or \textit{prevent} rate-dependent supercritical excursions. Consider a system governed by \eqref{eqn:FHN_compactified_slowtime} with additive noise. A trajectory that starts in the subcritical set $\pD_{\eps,r}^\text{(sub)}$ may be pushed across the quasithreshold $\pQ_{\eps,r}$ by the noise and thereby enter the supercritical set $\pD_{\eps,r}^\text{(sup)}$. It may then approximately follow the underlying deterministic dynamics to $\tilde{\pS}_{\eps}$, so that the system temporarily occupies the alternative state. Conversely, a trajectory that starts in $\pD_{\eps,r}^\text{(sup)}$ may be pushed by the noise into $\pD_{\eps,r}^\text{(sub)}$ before it reaches $\tilde{\pS}_\eps$; it can then settle near $\be_+^\text{(c)}$ without first becoming temporarily trapped in the alternative state. When the occurrence of these scenarios depends jointly on the rate of the parametric forcing and on the noise, we refer to them as combined rate- and noise-induced transitions~\cite{ritchie2017,slyman2023}. In this section, we quantify how likely such transitions are in the parametrically forced FitzHugh--Nagumo system driven by a stochastically evolving input current.

We incorporate \textit{randomness} into the input current of \eqref{eqn:FHN_compactified_slowtime} in precisely the same way as in Section \ref{sec:stochastic_FHN_autonomous}, by perturbing its (now, time-dependent) value with white-in-time Gaussian fluctuations. This motivates considering the dynamical system yielded by the following stochastic ordinary differential equation initial value problem of It\^{o} type,    
\begin{equation}\label{eqn:FHN_compactified_slowtime_stochastic}
    \begin{cases}
        \text{d} u = \eps^{-1}\big(-u^3 + u - v + \Lambda(s)\big)\,\text{d}t+\sigma\eps^{-1/2}\,\text{d}W \vspace{0.2cm}\\ 
        \text{d} v = (-\beta v + u)\,\text{d}t \vspace{0.2cm}\\ 
        \text{d}s = r\alpha(1-s^2)\,\text{d}t
    \end{cases},\,\begin{pmatrix}
        u(0) \\ v(0) \\ s(0)
    \end{pmatrix} = \begin{pmatrix}
        u^{(\ti)} \\ v^{(\ti)} \\ s^{(\ti)}
    \end{pmatrix} 
\end{equation} 
for $(u(t),v(t),s(t))\in\pD$. Within \eqref{eqn:FHN_compactified_slowtime_stochastic}, the noise intensity (proportional to $\sigma\in\bbR_+$) is scaled with the timescale separation parameter $\varepsilon$ in the same way as in \eqref{eqn:FHN_autonomous_slowtime_stochastic}.

We write the solution of the stochastic initial value problem \eqref{eqn:FHN_compactified_slowtime_stochastic} using the flow notation
$$
\bs\phi_{\sigma,\varepsilon,r}^\text{(c)}((u^{(\ti)},v^{(\ti)},s^{(\ti)}),W,t)=(u(t),v(t),s(t))\,,\,\,\,t\geq 0\,,
$$
the \textit{drift field} of \eqref{eqn:FHN_compactified_slowtime_stochastic} using the vectorised notation
\begin{equation}\label{eqn:compactified_drift_field}
\bs b_{\varepsilon,r}^\text{(c)}(u,v,s)\coloneqq\begin{pmatrix}
    \eps^{-1}\big(-u^3 + u - v + \Lambda(s)\big) \\ 
    (-\beta v + u) \\ 
    r\alpha(1-s^2)
\end{pmatrix}\,,
\end{equation}
and the \textit{diffusion matrix} of \eqref{eqn:FHN_compactified_slowtime_stochastic} using the matrix notation
\begin{equation*}
    \bs\Gamma_\varepsilon^{(\text{c})}(u,v) \equiv \bs\Gamma_\varepsilon^{(\text{c})} \coloneqq \begin{pmatrix}
        \eps^{-1/2} & 0 & 0 \\ 0 & 0 & 0 \\ 0 & 0 & 0 
    \end{pmatrix}\,.
\end{equation*}
Note the dependence of $\bs\phi_{\sigma,\varepsilon,r}^\text{(c)}((u^{(\ti)},v^{(\ti)},s^{(\ti)}),W,t)$ on the particular Wiener process $W$. 

\subsection{A committor function for combined rate- and noise-induced transitions}

If a dynamical system subjected to both parametric and stochastic forcing exhibits a transition event which is \textit{sensitive} to the rate of the parametric forcing and \textit{reliant} on the presence of the stochastic forcing, it is said to be a combined rate- and noise-induced transition event. For system \eqref{eqn:FHN_compactified_slowtime_stochastic}, such transitions could only manifest themselves across intervals of time which end \textit{after} the last time the parametric forcing $\Lambda$ attains its initial stationary value $I_-$, and begin \textit{before} the time $\Lambda$ attains its final stationary value $I_+$. In other words, the time between the onset and completion of any transition must include a period in which the parametric forcing is explicitly time-dependent for its occurrence to be influenced by the rate $r$. 
The parametric forcing $\Lambda(s)$ attains $I_-$ and $I_+$ with $s \leq s_-$ and $s \geq s_+$, respectively, whereby
\begin{align*}
s_-  \coloneqq g(t_-)=\tanh(\alpha t_-/2)\,,&&s_+\coloneqq g(t_+)=\tanh(\alpha t_+/2)\,.
\end{align*} 
With this in mind, to pursue combined rate- and noise-induced supercritical excursions, we will examine sample paths of \eqref{eqn:FHN_compactified_slowtime_stochastic} with a subcritical rate $r<r_\eps^\text{(crit)}$, which are initialised on the pullback attractor at the instance when the parametric forcing kicks in, that is, at 
\[
\bz_-^\text{(c)}\coloneqq (u_{I_-}^\text{(eq)},v_{I_-}^\text{(eq)},s_-)\,.
\]
As $\bz_-^\text{(c)}$ lies in the subcritical set $\pD_{\eps,r}^\text{(sub)}$ for $r<r_\varepsilon^\text{(cr)}$, stochastic forcing is required for any supercritical excursion to occur. Futher, if the stochastic forcing were to lead the sample path $\bs\phi_{\sigma,\varepsilon,r}^{(c)}(\bz_-^\text{(c)},W,t)$ to travel to $\tilde{\pS}_\varepsilon$, the preceeding time would necessarily include a period in which the parametric forcing $\Lambda$ is explicitly time-dependent (given the $s$-component of $\bz_-^\text{(c)}$). This being said, it is possible for the sample path to only visit $\tilde{\pS}_\varepsilon$ after a stage where $\Lambda$ has reached its final stationary value $I_+$ and the $u$- and $v$-components have settled to a close neighbourhood $\pN_{\varepsilon,\delta}(\be_{I_+}^\text{(FHN)})$ of $\be_{I_+}^\text{(FHN)}$, i.e., after a time $t\in\bbR_+$ in which 
\begin{flalign}
    \bs\phi_{\sigma,\varepsilon,r}^\text{(c)}(\bz_-^\text{(c)},W,t)&\in\pN_{\varepsilon,\delta}(\be_{I_+}^\text{(FHN)})\times [s_+,1] \label{eqn:stochastic_end_point_tracking} \\ 
    &\eqqcolon \pN_{\varepsilon,\delta}(\be_+^\text{(c)})\,.\label{eqn:stochastic_nonautonomous_nbh}
\end{flalign}
In these cases, the corresponding supercritical excursions are (insensitive to the rate $r$ and instead) solely induced by the stochastic forcing. This is because observing ``end-point tracking'' as in \eqref{eqn:stochastic_end_point_tracking} signals reaching the quasi-equilibrium state of the future autonomous regime which the sample paths inevitably (with probability one) arrive at, for any rate $r$; from this stage onwards, the dynamics of the system are fully independent of the rate $r$. As such, the sample path $\bs\phi_{\sigma,\varepsilon,r}^\text{(c)}(\bz_-^\text{(c)},W,t)$ exhibits a combined rate- and noise-induced supercritical excursion whenever it enters the $L$-metastable set \textit{before} settling to the quasi-equilibrium state in the future autonomous regime. That is, whenever 
\begin{equation}\label{eqn:supercritical_excursion_event}
\zeta_{\sigma,\varepsilon,r}^\text{(c)}(\bz_-^\text{(c)},W,\tilde{\pS}_\varepsilon) < \zeta_{\sigma,\varepsilon,r}^\text{(c)}(\bz_-^\text{(c)},W,\pN_{\varepsilon,\delta}(\be_+^\text{(c)}))\,\,\,\,\,\text{with}\,\,\,\,\,r<r_\varepsilon^\text{(cr)}\,,
\end{equation}
where 
\begin{flalign*}
    \zeta_{\sigma,\varepsilon,r}^\text{(c)}(\bz_-^\text{(c)},W,\tilde{\pS}_\varepsilon)&\coloneqq\text{inf}\{t>0:\bs\phi_{\sigma,\varepsilon,r}^\text{(c)}(\bz_-^\text{(c)},W,t)\in\tilde{\pS}_\varepsilon\}\,, \\ 
    \zeta_{\sigma,\varepsilon,r}^\text{(c)}(\bz_-^\text{(c)},W,\pN_{\varepsilon,\delta}(\be_+^\text{(c)}))&\coloneqq\text{inf}\{t>0:\bs\phi_{\sigma,\varepsilon,r}^\text{(c)}(\bz_-^\text{(c)},W,t)\in\pN_{\varepsilon,\delta}(\be_+^\text{(c)})\}\,,
\end{flalign*}
denote the first-entrance times of the sample path $\bs\phi_{\sigma,\varepsilon,r}^\text{(c)}(\bz_-^\text{(c)},W,t)$ into the sets $\tilde{\pS}_\varepsilon$ and $\pN_{\varepsilon,\delta}(\be_+^\text{(c)})$, respectively. Outcome \eqref{eqn:supercritical_excursion_event} occurs with probability
\begin{equation}\label{eqn:supercritical_excursion_probability}
    q_{\sigma,\varepsilon,r}^\text{(sup)}(\bz_-^\text{(c)})\coloneqq\bbP_{W\sim\rho_\Omega}\big(\zeta_{\sigma,\varepsilon,r}^\text{(c)}(\bz_-^\text{(c)},W,\tilde{\pS}_\varepsilon) < \zeta_{\sigma,\varepsilon,r}^\text{(c)}(\bz_-^\text{(c)},W,\pN_{\varepsilon,\delta}(\be_+^\text{(c)}))\big)\,.
\end{equation}
Conversely, there are cases where the combination of stochastic and rate-dependent parametric forcing prevents the occurrence of a supercritical excursion. For this, we examine sample paths of \eqref{eqn:FHN_compactified_slowtime_stochastic} which are again initialised at $\bz_-^\text{(c)}$, but this time with a supercritical rate $r>r_\varepsilon^\text{(cr)}$. Since $\bz_-^\text{(c)}$ lies in the supercritical set $\pD_{\varepsilon,r}^\text{(sup)}$ for $r > r_\varepsilon^\text{(cr)}$, stochastic forcing is required for the sample path to first settle to the quasi-equilibrium state $\pN_{\varepsilon,\delta}(\be_+^\text{(c)})$ in the future autonomous regime before visiting the $L$-metastable set $\tilde{\pS}_\varepsilon$. Further, all scenarios that culminate in settling to $\pN_{\varepsilon,\delta}(\be_+^\text{(c)})$ from $\bz_-^\text{(c)}$ will be directly preceded by the influence of $(t_+-t_-)/r$ units of explicit time-dependent forcing. As such, the sample path $\bs\phi_{\sigma,\varepsilon,r}^\text{(c)}(\bz_-^\text{(c)},W,t)$ exhibits a combined rate- and noise-induced subcritical path whenever it settles to the future quasi-equilibrium state without having previously visited the $L$-metastable set, i.e., whenever
\begin{equation}\label{eqn:subcritical_path_event}
\zeta_{\sigma,\varepsilon,r}^\text{(c)}(\bz_-^\text{(c)},W,\tilde{\pS}_\varepsilon) > \zeta_{\sigma,\varepsilon,r}^\text{(c)}(\bz_-^\text{(c)},W,\pN_{\varepsilon,\delta}(\be_+^\text{(c)}))\,\,\,\,\,\text{with}\,\,\,\,\,r> r_\varepsilon^\text{(cr)}\,.
\end{equation}
This particular outcome occurs with probability 
\begin{equation}\label{eqn:subcritical_path_probability}
q_{\sigma,\varepsilon,r}^\text{(sub)}(\bz_-^\text{(c)})\coloneqq\bbP_{W\sim\rho_\Omega}\big(\zeta_{\sigma,\varepsilon,r}^\text{(c)}(\bz_-^\text{(c)},W,\tilde{\pS}_\varepsilon) > \zeta_{\sigma,\varepsilon,r}^\text{(c)}(\bz_-^\text{(c)},W,\pN_{\varepsilon,\delta}(\be_+^\text{(c)}))\big)
\end{equation}
The quantities $q_{\sigma,\varepsilon,r}^\text{(sup)}(\bz_-^\text{(c)})$ and $q_{\sigma,\varepsilon,r}^\text{(sub)}(\bz_-^\text{(c)})$ are themselves examples of committor functions\cite{Vanden-Eijnden2006}.

Recall 
that for $r<r_\varepsilon^\text{(cr)}$ the trajectory $\bs\phi_{0,\varepsilon,r}^\text{(c)}(\bz_-^\text{(c)},t)$ reaches $\pN_{\varepsilon,\delta}(\be_+^\text{(c)})$ and subsequently settles to $\be_+^\text{(c)}$ without having previously visited the closure of the $L$-metastable set, $\text{cl}(\tilde{\pS}_\varepsilon)$. Alternatively, for $r>r_\varepsilon^\text{(cr)}$ the trajectory $\bs\phi_{0,\varepsilon,r}^\text{(c)}(\bz_-^\text{(c)},t)$ visits $\tilde{\pS}_\varepsilon$ before reaching $\text{cl}(\pN_{\varepsilon,\delta}(\be_+^\text{(c)}))$ and settling to $\be_+^\text{(c)}$. Since these scenarios occur deterministically without noise, it follows - from the asymptotic alignment result in Eq. \eqref{eqn:asymptotic_agreement} - that the probability of combined rate- and noise-induced (i) supercritical excursions and (ii) subcritical paths will both approach zero in the weak-noise limit, that is, 
\begin{flalign}\label{eqn:limiting_committors}
    \lim_{\sigma\downarrow 0}q_{\sigma,\varepsilon,r}^\text{(sup)}(\bz_-^\text{(c)})=0\,\,\,\,\,\text{for}\,\,\,\,\,r < r_\varepsilon^\text{(cr)}\,, \notag \\ 
    \lim_{\sigma\downarrow 0}q_{\sigma,\varepsilon,r}^\text{(sub)}(\bz_-^\text{(c)})=0\,\,\,\,\,\text{for}\,\,\,\,\,r > r_\varepsilon^\text{(cr)}\,.
\end{flalign}
Since $q_{\sigma,\varepsilon,r}^\text{(sup)}(\bz_-^\text{(c)})=1-q_{\sigma,\varepsilon,r}^\text{(sub)}(\bz_-^\text{(c)})$, result \eqref{eqn:limiting_committors} implies that $\displaystyle\lim_{\sigma\downarrow 0}q_{\sigma,\varepsilon,r}^\text{(sup)}(\bz_-^\text{(c)})=1$ for $r>r_\varepsilon^\text{(cr)}$. 

The following section \ref{subsec:combined_rate_noise_probability} discusses and evaluates the rates at which the probabilities $q_{\sigma,\varepsilon,r}^\text{(sup)}(\bz_-^\text{(c)})$ and $q_{\sigma,\varepsilon,r}^\text{(sub)}(\bz_-^\text{(c)})$ decay to zero in the weak-noise limit, for given subcritical and supercritical rates, respectively. Before this, however, we clarify a consequence of the dynamically forced setting in which combined rate- and noise-induced transitions are considered, which contrasts with the properties of the solely noise-induced transitions outlined in Section \ref{sec:stochastic_FHN_autonomous}. Whether a combined rate- and noise-induced transition occurs or not is determined by the time  
\[
t=T_{\sigma,\varepsilon,r}^\text{(c)}(\bz_-^\text{(c)},W)\coloneqq \text{min}\big(\zeta_{\sigma,\varepsilon,r}^\text{(c)}(\bz_-^\text{(c)},W,\tilde{\pS}_\varepsilon),\zeta_{\sigma,\varepsilon,r}^\text{(c)}(\bz_-^\text{(c)},W,\pN_{\eps,\delta}(\be_+^\text{(c)}))\big)\,.
\]
From $t = \zeta_{\sigma,\varepsilon,r}^\text{(c)}(\bz_-^\text{(c)},W,\pN_{\eps,\delta}(\be_+^\text{(c)}))$, the system has settled to the quasi-equilibrium state and cannot undertake any combined rate- and noise-induced transitions because there is no future explicit rate-dependent parametric forcing. For $r < r_\varepsilon^\text{(cr)}$ ($r>r_\varepsilon^\text{(cr)}$), a combined rate- and noise-induced supercritical excursion (subcritical path) from $\bz_-^\text{(c)}$ occurs if and only if $T_{\sigma,\varepsilon,r}^\text{(c)}(\bz_-^\text{(c)},W) = \zeta_{\sigma,\varepsilon,r}^\text{(c)}(\bz_-^\text{(c)},W,\tilde{\pS}_\varepsilon)$ ($T_{\sigma,\varepsilon,r}^\text{(c)}(\bz_-^\text{(c)},W)=\zeta_{\sigma,\varepsilon,r}^\text{(c)}(\bz_-^\text{(c)},W,\pN_{\varepsilon,\delta}(\be_+^\text{(c)}))$), with probability less than one. By contrast, a noise-induced supercritical excursion of the sample path $\bs\phi_{\sigma,\varepsilon,I}^\text{(FHN)}(\be_I^\text{(FHN)},W,t)$ of \eqref{eqn:FHN_autonomous_slowtime_stochastic} occurs with probability one. What's more, upon returning to $\pN_{\varepsilon,\delta}(\be_{I}^\text{(FHN)})$ after visiting $\tilde{\pS}_{\varepsilon,I}$, the system will undertake another noise-induced supercritical excursion again with probability one, and this cycle repeats endlessly. To summarise, each sample path $\bs\phi_{\sigma,\varepsilon,r}^\text{(c)}(\bz_-^\text{(c)},W,t)$ of \eqref{eqn:FHN_compactified_slowtime_stochastic} will exhibit \textit{at most one} combined rate- and noise-induced transition from $\bz_-^\text{(c)}$, while each sample path $\bs\phi_{\sigma,\varepsilon}^\text{(FHN)}(\be_I^\text{(FHN)},W,t)$ of \eqref{eqn:FHN_autonomous_slowtime_stochastic} will exhibit \textit{infinitely many} noise-induced supercritical excursions from $\pN_{\varepsilon,\delta}(\be_I^\text{(FHN)})$ to $\tilde{\pS}_{\varepsilon,I}$. This relates to the following difference in what underpins these transition events. The noise-induced supercritical excursions of Section \ref{sec:stochastic_FHN_autonomous} refer to trajectories which reach $\tilde{\pS}_{\varepsilon,I}$ from $\be_I^\text{(FHN)}$ via unconstrained paths: at every time instance before entering $\tilde{\pS}_{\varepsilon,I}$, the trajectory may explore the whole state space $\bbR^2\setminus\tilde{\pS}_{\varepsilon,I}$. By contrast, the combined rate- and noise-induced supercritical excursions introduced here refer to trajectories which reach $\tilde{\pS}_\varepsilon$ from $\bz_-^\text{(c)}$ via explicitly constrained paths: for all times $t\geq (t_+-t_-)/r$, the trajectory must remain outside of $\pN_{\varepsilon,\delta}(\be_+^\text{(c)})$. Further, the combined rate- and noise-induced subcritical paths introduced here refer to trajectories which reach $\pN_{\varepsilon,\delta}(\be_+^\text{(c)})$ from $\bz_-^\text{(c)}$ via oppositely constrained paths that remain outside of $\tilde{\pS}_{\varepsilon}$ at all times.

\begin{figure}[ht]
    \centering
    \includegraphics[width=5.75in]{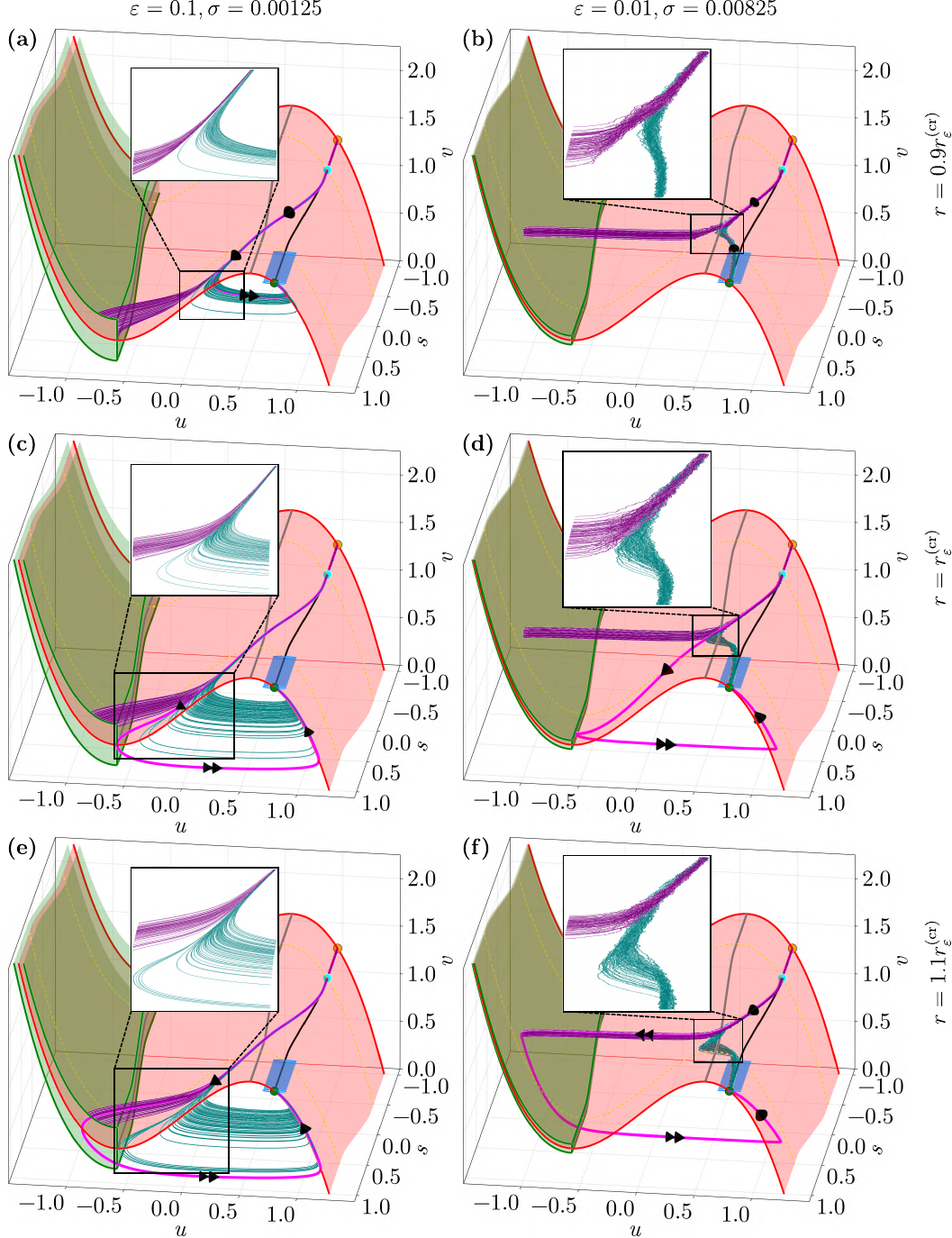}
    \caption{Sample paths $\bs\phi_{\sigma,\eps,r}^\text{(c)}(\bz_-^\text{(c)},W,t)$ of system \eqref{eqn:FHN_compactified_slowtime_stochastic} initialised at $\bz_-^\text{(c)} = (u_{I_-}^{(\text{eq})},v_{I_-}^{(\text{eq})},s_-)$ [cyan point] which terminate at $t = T_{\sigma,\eps,r}^\text{(c)}(\bz_-^\text{(c)},W)$ when entering either $\tilde{\pS}_\varepsilon$ or $\pN_{\varepsilon,\delta}(\be_+^{(\text{c})})$ [outlined in green and blue, respectively]. The parameter values considered are  $(\varepsilon,\sigma) = (0.1,0.00125),\,(\varepsilon,\sigma)=(0.01,0.00825),\,\delta=0.1$, and (deterministically) subcritical, critical, and supercritical $r$-values. Each panel shows 100 subcritical sample paths [teal flow lines], 100 supercritical sample paths [purple flow lines], the pullback attractor $A_{\varepsilon,r}(\be_-^\text{(c)})$ of the deterministic system \eqref{eqn:FHN_compactified_slowtime} [magenta curve, partly obscured by the sample paths], and the critical manifold $\pS_0$ [red surface]. The insets display a ``zoomed-in'' view of the state space where the sample paths exhibit ensemble splitting.}
\label{fig:stochastic_nonautonomous_sample_paths}
\end{figure}

Figure \ref{fig:stochastic_nonautonomous_sample_paths} presents sample paths $\bs\phi_{\sigma,\eps,r}^\text{(c)}(\bz_-^\text{(c)},W,t)$ of \eqref{eqn:FHN_compactified_slowtime_stochastic} across $t\in[0,T_{\sigma.\varepsilon,r}^\text{(c)}(\bz_-^\text{(c)},W)]$, 
showcasing realisations with
\begin{equation}\label{eqn:vist_L-metastable_set_first}
    \zeta_{\sigma,\varepsilon,r}^\text{(c)}(\bz_-^\text{(c)},W,\tilde{\pS}_\varepsilon) < \zeta_{\sigma,\varepsilon,r}^\text{(c)}(\bz_-^\text{(c)},W,\pN_{\eps,\delta}(\be_+^\text{(c)}))\,,
\end{equation}
and
\begin{equation}\label{eqn:visit_quasi-equilibrium_state_first}
    \zeta_{\sigma,\varepsilon,r}^\text{(c)}(\bz_-^\text{(c)},W,\pN_{\eps,\delta}(\be_+^\text{(c)})) < \zeta_{\sigma,\varepsilon,r}^\text{(c)}(\bz_-^\text{(c)},W,\tilde{\pS}_\varepsilon)\,,
\end{equation}
for deterministic subcritical, critical, and supercritical rates 
\begin{align*}
r = 0.9r_\eps^\text{(cr)},&&r = r_\eps^\text{(cr)},&&r = 1.1r_\eps^\text{(cr)}\,.
\end{align*}
Two timescale separations ($\varepsilon = 0.1$ and $\varepsilon=0.01$) are considered, together with practical weak noise strengths. For $r = 0.9r_\eps^\text{(cr)}$ the realisations \eqref{eqn:vist_L-metastable_set_first} and \eqref{eqn:visit_quasi-equilibrium_state_first} correspond to combined rate- and noise-induced \textit{supercritical excursions} occurring, and not occurring, respectively; for $r = 1.1r_\eps^\text{(cr)}$ the realisations \eqref{eqn:vist_L-metastable_set_first} and \eqref{eqn:visit_quasi-equilibrium_state_first} correspond to combined rate- and noise-induced \textit{subcritical paths} not occurring, and occurring, respectively. The sample paths are simulated via stochastic Monte Carlo rejection sampling using the semi-implicit Euler-Maruyama method\cite{higham2002} (see Section \begin{NoHyper}\ref{sm:numericalmethods}\end{NoHyper} of the supplementary material\cite{supplementary_material}).

\subsection{The probability of combined rate- and noise-induced transitions}\label{subsec:combined_rate_noise_probability}

Mirroring the approach of Section \ref{subsec:probability_autonomous}, one can quantify via large deviation principles\cite{freidlinwentzell,dembozeitouni} the probability that sample paths $\bs\phi_{\sigma,\eps,r}^\text{(c)}(\bz_-^\text{(c)},W,t)$ of \eqref{eqn:FHN_compactified_slowtime_stochastic}, with rates 
\begin{flalign*}
    \text{$r<r_\varepsilon^\text{(cr)}$ (subcritical) and $r > r_\varepsilon^\text{(cr)}$ (supercritical)}\,,
\end{flalign*}
undertake combined rate- and noise-induced supercritical excursions and subcritical paths, respectively, in the weak-noise limit $\sigma\rightarrow 0$. The events \eqref{eqn:supercritical_excursion_event} and \eqref{eqn:subcritical_path_event} occur when $\bs\phi_{\sigma,\eps,r}^\text{(c)}(\bz_-^\text{(c)},W,t)$ enters $\tilde{\pS}_\varepsilon$ before, respectively after, $\pN_{\eps,\delta}(\be_+^\text{(c)})$.

As a result of the globally attracting nature of the equilibrium point $\be_+^{(\tc)}$ within $\tilde{\pB}\coloneqq\bbR^2\times(-1,1]$, the energetic cost\footnote{Namely, the Freidlin-Wentzell action functional\cite{freidlinwentzell} for system \eqref{eqn:FHN_compactified_slowtime_stochastic}, \begin{gather*}
    S_T[\bs\psi]=\inf_{\bs\vartheta\in\Theta(\bs\psi)}\frac{1}{2}\int_0^T\lVert\bs\vartheta(s)\rVert^2\,\text{d}s,\, \\ 
    \Theta(\bs\psi)\coloneqq\bigg\{\bs\vartheta\in L_2([0,T],\bbR^3):\bs\psi(t)=\bz^\text{(i)}+\int_0^t \bs b_{\varepsilon,r}^\text{(c)}(\bs\psi(s))+\bs\Gamma_\varepsilon^\text{(c)}\bs\vartheta(s)\,\text{d}s,\,\forall\,t\in[0,T]\bigg\}\,.
\end{gather*}\vspace*{-0.8cm}\textcolor{white}{.}} of any path
\[
\bs\psi\in C^0\big([0,T],\pD\setminus\pN_{\varepsilon,\delta}(\be_+^\text{(c)})\big)\text{ with }\bs\psi(0) = \bz_-^\text{(c)}\,,\,T\in\bbR_+\,,
\]
which continuously lies outside of the neighbourhood $\pN_{\eps,\delta}(\be_+^{(\tc)})$ throughout $t \in [0,T]$, tends to infinity as the length $T$ of the time interval tends to infinity. The accumulation of the ``action'' required to sustain not visiting the neighbourhood must diverge on account of the fact that there are no invariant sets outside $\pN_{\eps,\delta}(\be_+^\text{(c)})$. An important consequence of this is that the least energetically costly combined rate- and noise-induced supercritical excursion (subcritical path) will be a path which visits the $L$-metastable set (neighbourhood of the stable equilibrium) in finite time, rather than only reaching it asymptotically as time approaches infinity. For this reason, the following $T$-dependent energetic costs 
\begin{gather*}
\inf_{\bs\psi\in\pP_{\eps,<T}^\text{(sup)}(\bz_-^{(\text{c})})}S_{T}[\bs\psi]\eqqcolon\Delta V_{\varepsilon,r,<T}^\text{(sup)}(\bz_-^\text{(c)})\,,\quad\inf_{\bs\psi\in\pP_{\eps,<T}^\text{(sub)}(\bz_-^{(\text{c})})}S_{T}[\bs\psi]\eqqcolon\Delta V_{\varepsilon,r,<T}^\text{(sub)}(\bz_-^\text{(c)})\,,
\end{gather*}
for
\begin{gather*}
\pP_{\eps,<T}^\text{(sup)}(\bz_-^{(\text{c})}) \coloneqq\big\{\bs\psi\in C^0([0,T],\pD):\bs\psi(0)=\bz_-^\text{(c)},\,\zeta_T(\bs\psi,\tilde{\pS}_{\eps})<\zeta_T(\bs\psi,\pN_{\eps,\delta}(\be_+^\text{(c)}))\big\}\,, \\
\pP_{\eps,<T}^\text{(sub)}(\bz_-^\text{(c)})\coloneqq\big\{\bs\psi\in C^0([0,T],\pD):\bs\psi(0)=\bz_-^\text{(c)},\,\zeta_T(\bs\psi,\tilde{\pS}_{\eps})>\zeta_T(\bs\psi,\pN_{\eps,\delta}(\be_+^\text{(c)}))\big\}\,, \\ 
\zeta_T(\bs\psi,\tilde{\pS}_{\eps}) \coloneqq \inf\{t\in[0,T]:\bs\psi(t)\in\tilde{\pS}_{\eps}\}\,, \\ 
\zeta_T(\bs\psi,\pN_{\eps,\delta}(\be_+^\text{(c)}))\coloneqq\inf\{t\in[0,T]:\bs\psi(t)\in\pN_{\eps,\delta}(\be_+^\text{(c)})\}\,,
\end{gather*}
attain their minimum values across $T\in\bbR_+$ at some (finite) $T^\text{(m)}\in\bbR_+$. This is to say that there exists some sufficiently large $T^\text{(m)}\in\bbR_+$ for which 
\begin{equation*}
\begin{array}{ccc}
\begin{aligned}
\Delta V_{\varepsilon,r,<T^\text{(m)}}^\text{(sup)}(\bz_-^\text{(c)})
&= \min_{T\in\bbR_+}
\big(\Delta V_{\varepsilon,r,<T}^\text{(sup)}(\bz_-^\text{(c)})\big)\\
&\eqqcolon \Delta V_{\eps,r}^\text{(sup)}(\bz_-^\text{(c)})\,,
\end{aligned}
&\hspace{1cm}&
\begin{aligned}
\Delta V_{\varepsilon,r,<T^\text{(m)}}^\text{(sub)}(\bz_-^\text{(c)})
&= \min_{T\in\bbR_+}
\big(\Delta V_{\varepsilon,r,<T}^\text{(sub)}(\bz_-^\text{(c)})\big)\\
&\eqqcolon \Delta V_{\eps,r}^\text{(sub)}(\bz_-^\text{(c)})\,.
\end{aligned}
\end{array}
\end{equation*}
These properties imply that the probability of combined rate- and noise-induced supercritical excursions and subcritical paths satisfy the following exponential scalings
\begin{flalign}
q_{\sigma,\eps,r}^\text{(sup)}(\bz_-^\text{(c)})&\overset{\sigma\downarrow 0}{\asymp}\exp\big(-\sigma^{-2}\Delta V_{\eps,r}^\text{(sup)}(\bz_-^\text{(c)})\big)\,,\hspace{0.5cm}\text{for } r<r_\eps^\text{(cr)}\,,\label{eqn:full_supercritical_LDP}\\ 
q_{\sigma,\eps,r}^\text{(sub)}(\bz_-^\text{(c)})&\overset{\sigma\downarrow 0}{\asymp}\exp\big(-\sigma^{-2}\Delta V_{\eps,r}^\text{(sub)}(\bz_-^\text{(c)})\big)\,,\hspace{0.5cm}\text{for } r>r_\eps^\text{(cr)}\,,\label{eqn:full_subcritical_LDP}
\end{flalign}
 with respect to $\sigma^{-2}$.

We quantify the asymptotic rate $\Delta V_{\eps,r}^\text{(sup)}(\bz_-^\text{(c)})$ with which the probability of combined rate- and noise-induced \textit{supercritical excursions} decays to zero with $\sigma^{-2}$. Indeed, this serves as a means of assessing - for instance - whether the inclusion of rate-dependent parametric forcing enhances or reduces the probability of supercritical excursions relative to the solely stochastically forced case with the same weak noise strength (c.f. Section \ref{sec:stochastic_FHN_autonomous}). It is similarly worthwhile to quantify the asymptotic rate $\Delta V_{\eps,r}^\text{(sub)}(\bz_-^\text{(c)})$ with which the probability of combined rate- and noise-induced \textit{subcritical paths} decays to zero with $\sigma^{-2}$, to investigate whether these excursion-preventing events are more or less probable than corresponding excursion-inducing events under the same weak noise strength. The asymptotic rates alone determine which of two such events is relatively rarer (whenever they are different in value) because the ratio of the probabilities of the events is logarithmically equivalent to an exponential function with an exponent proportional to the \textit{difference} of the corresponding asymptotic rates.            

\begin{figure}[ht]
    \centering
    \includegraphics[width=6.4in]{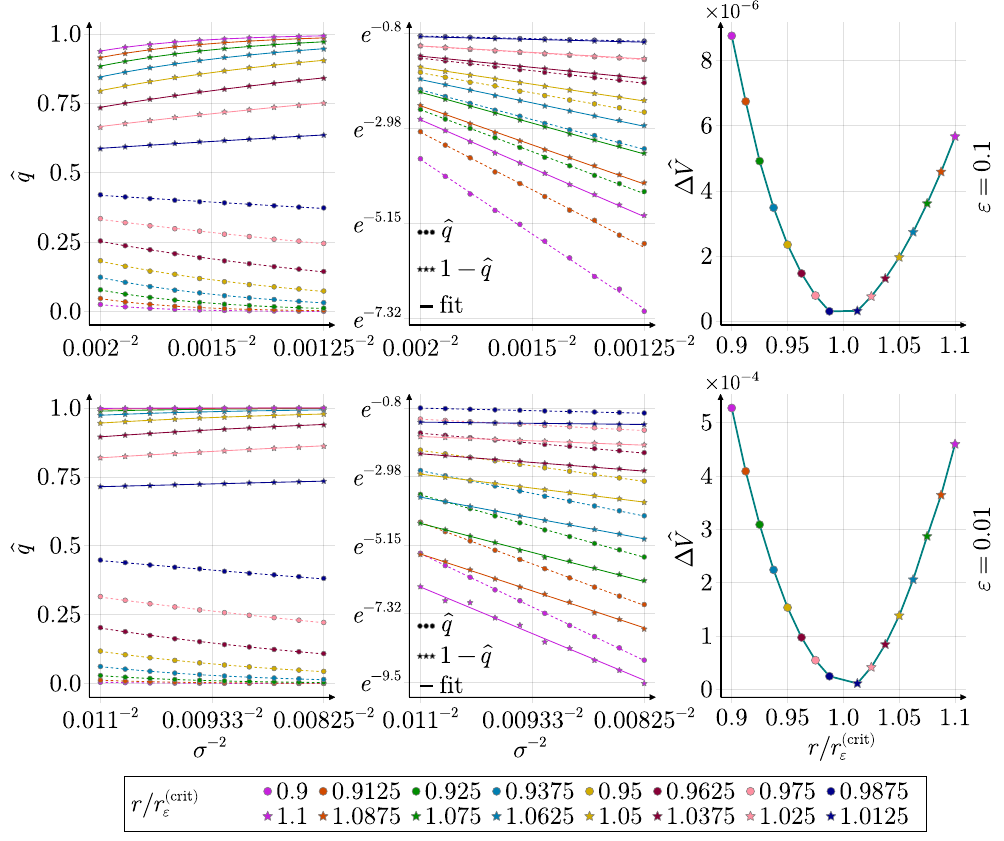}
    \caption{Vulnerability to combined rate- and noise-induced supercritical excursions and subcritical paths in the stochastic nonautonomous FitzHugh-Nagumo system \eqref{eqn:FHN_compactified_slowtime_stochastic}. Left: numerical estimates $\hat{q}^\text{(sup)}$ of the committor function $q_{\sigma,\eps,r}^\text{(sup)}(\bz_-^\text{(c)})$ against $\sigma^{-2}$ across intervals of finite noise, for a range of subcritical and supercritical $r$-values [coloured points and stars respectively]. For $r<r_\varepsilon^\text{(cr)}$ the probability $q_{\sigma,\varepsilon,r}^\text{(sup)}$ is estimated directly using the adaptive multilevel splitting (AMS) algorithm [circles]; for $r>r_\varepsilon^\text{(cr)}$ the complementary probability $q_{\sigma,\varepsilon,r}^\text{(sub)} = 1 - q_{\sigma,\varepsilon,r}^\text{(sup)}(\bz_-^\text{(c)})$ is first estimated as $\hat{q}^\text{(sub)}$ using AMS, before setting $\hat{q}^\text{(sup)} = 1-\hat{q}^\text{(sub)}$ [stars]. The coloured intersecting curves arise from weighted least-square exponential fits of $\hat{q}^\text{(sup)}$ and $\hat{q}^\text{(sub)}$: for $r<r_\varepsilon^\text{(cr)}$ fits of $\hat{q}^\text{(sup)}$ are plotted directly [dashed curves]; for $r>r_\varepsilon^\text{(cr)}$ fits of $\hat{q}^\text{(sub)}$ are made, and one minus the fitted functions are plotted [solid curves]. Centre: the estimates $\hat{q}^\text{(sup)}$ [circles] and $\hat{q}^\text{(sub)}$ [stars], for subcritical and supercritical rates respectively, plotted against $\sigma^{-2}$ using a logarithmic scale on the vertical axis. The exponential fits to $\hat{q}^\text{(sup)}$ [dashed lines] and $\hat{q}^\text{(sub)}$ [solid lines] are again shown. Right: the slopes returned by the fitting procedure as a function of $r/r_\eps^\text{(cr)}$, which provide estimates $\Delta\hat{V}$ of $\Delta V_{\eps,r}^\text{(sup)}(\bz_-^\text{(c)})$ and $\Delta V_{\eps,r}^\text{(sub)}(\bz_-^\text{(c)})$ for subcritical and supercritical rates, respectively. The matching colours correspond to eight pairs of subcritical and supercritical rates equidistant from $r_\eps^\text{(cr)}$, as listed in the legend.} 
    \label{fig:committor_function_scaling_v3}
\end{figure}

The values of $\Delta V_{\eps,r}^\text{(sup)}(\bz_-^\text{(c)})$ and $\Delta V_{\eps,r}^\text{(sub)}(\bz_-^\text{(c)})$ can be approximated numerically through a series of estimates for $q_{\sigma,\eps,r}^\text{(sup)}(\bs z_-^\text{(c)})$ and $q_{\sigma,\eps,r}^\text{(sub)}(\bs z_-^\text{(c)})$ across a range of weak-noise values, followed by a fitting of the data to extract the rates of exponential scaling. The probabilities $q_{\sigma,\varepsilon,r}^\text{(sup)}(\bz_-^\text{(c)})$ and $q_{\sigma,\varepsilon,r}^\text{(sub)}(\bz_-^\text{(c)})$ are estimated via the adaptive multilevel splitting algorithm (introduced in Section \begin{NoHyper}\ref{sm:AMS}\end{NoHyper} of the supplementary material\cite{supplementary_material}). The estimates are gathered in the left column of Figure \ref{fig:committor_function_scaling_v3} for timescale separation values $\eps=0.1$, $\eps=0.01$, and rates
$$ 
r \in R_\varepsilon^\text{(sub)}\cup R_\varepsilon^\text{(sup)}\,,
$$
\vspace{-1.875cm}
\begin{gather*}  
R_\varepsilon^\text{(sub)}\coloneqq\{(1-k)r_\eps^\text{(cr)}:k\in K\}\,, \\ 
R_\varepsilon^\text{(sup)}\coloneqq\{(1+k)r_\eps^\text{(cr)}:k\in K\}\,,
\end{gather*}
where
$$
K\coloneqq\{0.0125,0.025,0.0375,0.05,0.0625,0.075,0.0875,0.1\}
$$
i.e., for eight equally spaced values above and below the critical rate, reaching $\pm10\%$ of its value. A range of ten values for the noise strength $\sigma$ were prescribed between extreme values $\sigma_\varepsilon^\text{(min)}$ and $\sigma_\varepsilon^\text{(max)}$ (specified below in \eqref{eqn:noise_strength_extrema}). The intermediate values for $\sigma$ between $\sigma_\varepsilon^\text{(min)}$ and $\sigma_\varepsilon^\text{(max)}$ were chosen so that the corresponding values for $\sigma^{-2}$ are equally spaced, that is, 
$$ 
\sigma^{-2} \in \Sigma_\varepsilon\,,
$$
\vspace{-1.75cm}
\begin{align*}  
\Sigma_\varepsilon\coloneqq\{(\sigma_\varepsilon^\text{(max)})^{-2}(1-\kappa)+(\sigma_\varepsilon^\text{(min)})^{-2}\kappa:\kappa\in\{0,\,1/9,\,2/9,\,...,\,1\}\}\,.
\end{align*}
The following values for the upper and lower bounds for the noise strengths,   
\begin{equation}\label{eqn:noise_strength_extrema}
    \sigma_\varepsilon^\text{(min)}=\begin{cases}
        0.00125 & \text{for }\varepsilon = 0.1 \\ 
        0.00825 & \text{for }\varepsilon = 0.01
    \end{cases} \quad\text{and}\quad\sigma_\varepsilon^\text{(max)}=\begin{cases}
        0.002 & \text{for }\varepsilon = 0.1 \\ 
        0.011 & \text{for }\varepsilon = 0.01
    \end{cases}\,,
\end{equation}
were specifically chosen. This is on account of the fact that - across noise strengths $\sigma \in \{\sigma\in\bbR_+:\sigma^{-2}\in\Sigma_\varepsilon\}$ with extreme values as in \eqref{eqn:noise_strength_extrema}  - the occurrence of combined rate- and noise-induced supercritical excursions (resp. subcritical paths) are so rare as to make the estimation of their probability lie at the limit of what is practically computable, for $r = 0.9r_\varepsilon^\text{(cr)}$ (resp. $r=1.1r_\varepsilon^\text{(cr)})$. In other words, we are as close to the weak-noise limit as we can afford. 

The outcome of the fits (the procedure behind which is introduced in Section \begin{NoHyper}\ref{sm:wlsqfits}\end{NoHyper} of the supplementary material\cite{supplementary_material}) shows in particular that the probability of undertaking combined rate- and noise-induced supercritical excursions and subcritical paths in system \eqref{eqn:FHN_compactified_slowtime_stochastic} decays exponentially at rates
\begin{equation*}
\Delta V_{\varepsilon,r}^\text{(sup)}(\bz_-^\text{(c)}),\,\Delta V_{\varepsilon,r}^\text{(sub)}(\bz_-)\overset{\text{o.o.m}}{\sim}
\begin{cases}
    10^{-7} \text{ to } 10^{-6} & \text{for } \varepsilon = 0.1 \\ 
    10^{-5} \text{ to } 10^{-4}  & \text{for } \varepsilon = 0.01
\end{cases}\,,
\end{equation*}
for $r\in R_\varepsilon^\text{(sub)}$ and $r\in R_\varepsilon^\text{(sup)}$, respectively. One conclusion to draw from this is that
\begin{flalign}
\text{for all } r\in R^\text{(sub)}\,,\,\Delta V_{\eps,r,<T}^\text{(sup)}(\bz_-^\text{(c)})\ll\Delta V_{\eps,I=0.75,T}^\text{(sup)}\overset{\text{o.o.m}}{\sim}10^{-2}\,, \label{eqn:comparative_supercritical_excursions} 
\end{flalign}
for every $T\in\bbR_+$ sufficiently large such that 
\begin{equation}\label{eqn:sufficiently_large_travel_time}
\Delta V_{\eps,r,<T}^\text{(sup)}(\bz_-^\text{(c)}) \approx \Delta V_{\eps,r}^\text{(sup)}(\bz_-^\text{(c)}) 
\end{equation}
Taking $T\gg (t_+-t_-)/r$ will likely suffice to ensure \eqref{eqn:sufficiently_large_travel_time}. The result in \eqref{eqn:comparative_supercritical_excursions} implies that, with subcritical rates on the verge of criticality, combined rate- and noise-induced supercritical excursions are overwhelmingly more probable than solely noise-induced supercritical excursions across all sufficiently large bounded time intervals, under identically weak noise strengths $0<\sigma\ll 1$. 

For all pairings of rates $r_-\in R_\varepsilon^\text{(sub)},\,r_+\in R_\varepsilon^\text{(sup)}$ which are equidistant to $r_\varepsilon^\text{(cr)}$ from either side and satisfy 
\[
1-r_-/r_\varepsilon^\text{(cr)}=r_+/r_\varepsilon^\text{(cr)}-1\in\begin{cases}K\setminus\{0.0125\} & \text{for }\varepsilon=0.1 \\ 
K & \text{for }\varepsilon = 0.01\end{cases}\,,
\]
the estimates $\Delta \hat{V}_{\eps,r_-}^\text{(sup)}(\bz_-^\text{(c)})$ and $\Delta \hat{V}_{\eps,r_+}^\text{(sub)}(\bz_-^\text{(c)})$ of the asymptotic rates $\Delta V_{\eps,r_-}^\text{(sup)}(\bz_-^\text{(c)})$ and $\Delta V_{\eps,r_+}^\text{(sub)}(\bz_-^\text{(c)})$ in Figure \ref{fig:committor_function_scaling_v3} satisfy
\begin{equation}\label{eqn:comparative_subcritical_and_supercritical_excursions_v2}
\Delta \hat{V}_{\eps,r_-}^\text{(sup)}(\bz_-^\text{(c)}) > \Delta \hat{V}_{\eps,r_+}^\text{(sub)}(\bz_-^\text{(c)})\,.
\end{equation}
If both estimates are accurate, this would suggest that the probability of \textit{preventing} supercritical excursions becomes overwhelmingly larger than the probability of \textit{inducing} supercritical excursions,  under identically weak noise strengths $0<\sigma\ll 1$. This being said, the inequality \eqref{eqn:comparative_subcritical_and_supercritical_excursions_v2} becomes less clearly pronounced in a closer neighbourhood of criticality with smaller joint values of $1-r_-/r_\varepsilon^\text{(cr)}=r_+/r_\varepsilon^\text{(cr)}-1$; in fact, for the weaker timescale separation $\varepsilon=0.1$, the reverse inequality holds for $1-r_-/r_\varepsilon^\text{(cr)}=r_+/r_\varepsilon^\text{(cr)}\in\{0.0125,0.025\}$. Further statistical testing would therefore be required to more confidently make conclusions about the asymmetry of these processes (see Section \ref{sec:discuss} for suggestions in this direction). Across the entire range of noise strengths considered for estimating $\Delta V_{\eps,r_-}^\text{(sup)}(\bz_-^\text{(c)})$ and $\Delta V_{\eps,r_+}^\text{(sub)}(\bz_-^\text{(c)})$, we observed that 
\begin{flalign}
&\hat{q}_{\sigma,\eps,r_-}^\text{(sup)}(\bz_-^\text{(c)}) < \hat{q}_{\sigma,\eps,r_+}^\text{(sub)}(\bz_-^\text{(c)})\,,\hspace{0.5cm}\text{for } \eps=0.1\,, \notag 
\\
&\hat{q}_{\sigma,\eps,r_-}^\text{(sup)}(\bz_-^\text{(c)}) > \hat{q}_{\sigma,\eps,r_+}^\text{(sub)}(\bz_-^\text{(c)})\,,\hspace{0.5cm}\text{for } \eps=0.01\,, \label{eqn:comparative_probabilities_eps-0.01}
\end{flalign}
for pairings $(r_-,r_+)$ of rates which are equidistant to $r_\varepsilon^\text{(cr)}$ and satisfy
\[
1-r_-/r_\varepsilon^\text{(cr)}=r_+/r_\varepsilon^\text{(cr)}-1 \in \begin{cases}
    K\setminus\{0.0125\} & \text{for } \varepsilon = 0.1 \\ K & \text{for } \varepsilon = 0.01
\end{cases}\,.
\]
This can be seen by inspecting the central column of Figure \ref{fig:committor_function_scaling_v3}. The relation in \eqref{eqn:comparative_probabilities_eps-0.01}, which goes against the conclusion drawn after Eq. \eqref{eqn:comparative_subcritical_and_supercritical_excursions_v2}, might be an effect of the range of finite noise strengths considered (where subexponential prefactors\cite{grafke2024sharp,bouchet2016generalisation} can explain such supposedly anomalous results) and the inequality may permanently reverse in sign below sufficiently small values of $\sigma$.

Finally, note from Figure \ref{fig:committor_function_scaling_v3} that
\begin{equation}\label{eqn:comparative_asymptotic_rates_nonautonomous}
\begin{split}
    \forall\, r \in R_\varepsilon^\text{(sub)},\,&\Delta \hat{V}_{\eps=0.1,r}^\text{(sup)}(\bz_-^\text{(c)}) \overset{\text{o.o.m.}}{\sim}  10^{-7}\text{ to }10^{-6} \ll \Delta \hat{V}_{\eps=0.01,r}^\text{(sup)}(\bz_-^\text{(c)})\overset{\text{o.o.m.}}{\sim} 10^{-5}\text{ to }10^{-4}\,, \\
    \forall\, r \in R_\varepsilon^\text{(sup)},\,&\Delta \hat{V}_{\eps=0.1,r}^\text{(sub)}(\bz_-^\text{(c)}) \overset{\text{o.o.m.}}{\sim} 10^{-7}\text{ to }10^{-6} \ll \Delta \hat{V}_{\eps=0.01,r}^\text{(sub)}(\bz_-^\text{(c)})\overset{\text{o.o.m.}}{\sim} 10^{-5}\text{ to }10^{-4}\,.
\end{split}
\end{equation}
This is the opposite of what we found for purely noise-induced supercritical excursions in the stochastic autonomous system \eqref{eqn:FHN_autonomous_slowtime_stochastic}, where the quasipotential barrier heights decreased as $\varepsilon$ decreased (see \eqref{eqn:comparative_asymptotic_rates_autonomous}). For the combined rate- and noise-driven system \eqref{eqn:FHN_compactified_slowtime_stochastic}, \eqref{eqn:comparative_asymptotic_rates_nonautonomous} shows the reverse, with the following implications. 
If one fixes a weak noise strength and considers rates $r_\varepsilon\in R_\varepsilon^\text{(sub)}$ (resp., $r_\varepsilon\in R_\varepsilon^\text{(sup)})$ with $1 - r_\varepsilon/r_\varepsilon^\text{(cr)}$ held fixed, then as $\varepsilon$ decreases by an order of magnitude from $0.1$ to $0.01$, the probability of a combined rate- and noise-induced supercritical excursion (resp., a subcritical path) becomes exponentially smaller.

\section{Discussion}
\label{sec:discuss}


This paper uncovered mechanisms behind combined rate- and noise-induced transitions in a multiple timescale system, using a nonautonomous FitzHugh-Nagumo system as a paradigmatic example. The dynamics are always deterministically monostable, yet admit trajectories that undertake long-lived excursions (excitations) away from the asymptotic equilibrium. Exploiting the timescale separation, we used the framework of metastability to define \textit{supercritical excursions} as trajectories which pass through an $L$-metastable set enclosing the \textit{slow manifolds} associated with an attracting branch of the critical manifold away from equilibrium (cf. ``rivers'' in \citet{letson2018new}). In the nonautonomous case, the quasithreshold separating supercritical and subcritical behaviour is built around a family of \textit{maximal canard trajectories}, arising from a folded saddle-node type-I singularity, and their interactions with the fast dynamics (Sections \ref{subsec:geometric_singular_perturbation_theory} - \ref{subsec:nonautonomous_bifurcation_analysis}). We focused on this reversible tipping phenomenon because it cannot be observed in one-dimensional systems.

Stochastic forcing allows trajectories to cross the quasithreshold and thereby either induce or prevent supercritical excursions. Using the (trajectory) adaptive multilevel splitting (AMS) algorithm, we quantified the probability of inducing supercritical excursions under (i) stochastic forcing alone and (ii) combined stochastic and subcritical rate-dependent parametric forcing. The probability of inducing supercritical excursions increased drastically when both forcing types acted together. We likewise used AMS to quantify the probability that noise prevents a deterministically expected rate-dependent excursion.

The committor function scalings in Figure \ref{fig:committor_function_scaling_v3} suggest that the action (``energetic cost'') of combined rate- and noise-induced supercritical excursions (subcritical paths) in system \eqref{eqn:FHN_compactified_slowtime_stochastic} vanishes as the rate $r$ approaches its critical value $r_\varepsilon^\text{(cr)}$ from below (above). This can be expected since the minimal distance between the pullback attractor and the quasithreshold tends to zero in this limit. However, the action vanishes at different rates for different timescale separations: for rates equidistant from criticality, energetic costs are larger for $\eps=0.01$ than for $\eps=0.1$ (of order $10^{-4}$ and $10^{-6}$, respectively). This is because, for fixed $r/r_\eps^\text{(cr)}$, the pullback attractor $A_{\varepsilon,r}(\be_-^\text{(c)})$ tracks the quasithreshold $\pQ_{\eps,r}$ less closely as the timescale separation strengthens: stronger repulsion from the slow manifolds associated with $\pS_0^{(\tr)}$ means trajectories must come closer to this surface to track it to the same extent (Figure \ref{fig:nonautonomous_metastability_landscapes}). An alternative comparison would fix the degree to which the pullback attractor resembles canard trajectories rather than fixing $r/r_\eps^\text{(cr)}$.

In most cases\footnote{All cases excluding a close neighbourhood around criticality for $\varepsilon = 0.1$.}, Figure \ref{fig:committor_function_scaling_v3} further suggests that it is less costly for noise to \textit{prevent} a supercritical excursion (supercritical $r$) than to \textit{induce} one (equidistant subcritical $r$). We speculate that this reflects the asymmetry of the critical manifold: the optimal inducing path spends considerable time resisting convergence along increasingly attracting parts of $\pS_0^{(\ta,2)}$ away from $\pF^{(2)}$, as the evolution of $s$ slows near 1, whereas the optimal preventing path spends most of its time near $\pF^{(1)}$, where $\pS_0^{(\ta,1)}$ is least attracting. This trend depends subtly on the interplay of pullback attractor, quasithreshold, drift and diffusion, and may differ under other stochastic forcing.

Comparing the energetic costs in Figure \ref{fig:excitable_sample_paths_autonomous} ($10^{-2}$) and Figure \ref{fig:committor_function_scaling_v3} ($\leq 10^{-4}$) shows that, for the same weak noise, purely noise-induced supercritical excursions are negligibly likely compared with combined rate- and noise-induced ones, as the ratio of probabilities decays exponentially to zero with $\sigma^{-2}$. This appears to be because the pullback attractor $A_{\varepsilon,r}(\be_-^\text{(c)})$ of system \eqref{eqn:FHN_compactified_slowtime} comes closer to $\pQ_{\eps,r}$ than the equilibrium $\be_{I_+}^\text{(FHN)}$ of system \eqref{eqn:FHN_autonomous_slowtime} is to $\pQ_{\eps,I_+}$ (Figure \ref{fig:nonautonomous_metastability_landscapes}). Future work could examine how these costs compare as $r$ tends to 0 or $\infty$; co-varying the timescale separation with $r$ would yield a double limit problem, adding a stochastic component to studies of rate-induced tipping in 2-fast, 1-slow systems\cite{osullivan2023,alkhayuon2025icms}. 

Our scaling estimates are limited by (i) extrapolation to the weak-noise limit from finite noise strengths, and (ii) estimation error from a finite number of clones and realisations of the splitting algorithm. An alternative approach would consider the Hamiltonian extension of \eqref{eqn:FHN_compactified_slowtime_stochastic} and analyse its most probable combined transition paths\cite{slyman2023}.   

Our work connects to several strands of literature. Maximal canards are commonly used to separate sub-threshold and ``spiking'' behaviour in autonomous excitable systems of type B\footnote{Following \citet{wieczorek2011}: excursion amplitude depends \textit{continuously} on perturbation size, as no saddles or other unstable invariant sets are present.}$^{,}$
\cite{kuehn2022,desroches2013inflection,osullivan2023,khovanov2013noise,franovic2015activation,chen2017crossing,zhu2016subthreshold}, sometimes described as the ``ghost'' of the singular separatrix\cite{khovanov2013noise,franovic2015activation,chen2017crossing}, i.e. our $\pQ_{0}$. The autonomous system \eqref{eqn:FHN_autonomous_slowtime} exhibits the \textit{quasithreshold instability} of \citet{osullivan2023}, which generalises the threshold instability of \citet{wieczorek2023} to systems lacking a regular edge state. Whereas the \textit{R-tipping quasi-threshold} of \citet{osullivan2023} is a codimension-zero family of canards ``without head'', giving an interval of critical rates, our quasithreshold is typically codimension-one, yielding an isolated critical rate for system \eqref{eqn:FHN_compactified_slowtime}. Moreover, while singular canards associated with folded singularities (\textit{singular R-tipping edge states}\cite{szmolyan2001,wieczorek2011}) approximate the tracking/tipping threshold locally near the critical manifold\cite{perryman2014adapting,osullivan2023,vanselow24}, our quasithreshold separates subcritical and supercritical trajectories across the entire state space of the nonsingular system. While its computation was non-trivial, the proposed \textit{quasi-edge-state} may offer further insight into its structure.

By contrast, \textit{irreversible} rate-dependent transitions in bistable climate models have been studied\cite{chapman2024tipping,slyman2024tipping,slyman2025tipping}, though parametric and stochastic forcing were rarely combined; \citet{slyman2024tipping} showed combined transitions between ``storm-on'' and ``storm-off'' states but left their probability to future work, requiring extension of the methods of \citet{slyman2023}. An alternative to our large deviation perspective is direct analysis of the Fokker-Planck equation\cite{risken1989fpe}, as in \citet{ritchie2017}, which applies at any finite noise strength. Beyond the multiscale setting, our study extends \citet{slyman2023} and \citet{ritchie2017} by also quantifying \textit{subcritical paths}, where noise prevents a deterministic rate-dependent excursion.

With accelerating global warming and extreme weather, combined trend-induced and random forcing is relevant across the multiscale climate system, yet its effects remain poorly understood. This work shows that such combinations can introduce significant uncertainty into the onset of tipping, so careful modelling is needed to reliably quantify the risks of a given emissions scenario. Natural next steps include studying combined transitions in conceptual ocean circulation models\cite{alkhayuon2019basin,chapman2024tipping,soons2024optimal} and soil- or oceanic-carbon systems\cite{osullivan2023,slyman2025tipping}, and developing \textit{early warning signs}\cite{ritchie2016,morr2024,huang2024deep,lohmann2025} for these events. For instance, one could pair the methods of \citet{chapman2025} and \citet{lohmann2025choosing} with the most probable combined transition path to construct geometric, proximity-based indicators.

\section*{Supplementary Material}

See the supplementary material\cite{supplementary_material} for a description of the numerical methods applied to generate the data and plots of this study, and a more detailed analysis of the metastability landscape of the compactified forced FitzHugh-Nagumo system. In particular, details on the theory and implementation of numerical integration schemes, the adaptive multilevel splitting algorithm, weighted least-squares fits, and quasithreshold computation are provided in Section \begin{NoHyper}\ref{sm:numericalmethods}\end{NoHyper}. The bifurcation analysis of the desingularised reduced slow subsystem of the compactified forced FitzHugh-Nagumo system is presented in Section \begin{NoHyper}\ref{sm:bifofreduced}\end{NoHyper}. Explicit set-theoretic formulations of the quasithreshold of the compactified forced FitzHugh-Nagumo system in the singular ($\varepsilon = 0$) and nonsingular ($0<\varepsilon\ll 1$) cases are provided in Sections \begin{NoHyper}\ref{sm:singular_supercritical_set}\end{NoHyper} and \begin{NoHyper}\ref{sm:decomposed_quasithreshold}\end{NoHyper} respectively. Further, tables listing the data generated by the Monte Carlo sampling and adaptive multilevel algorithms implemented in this study are provided, as well as a review of the notation used throughout the paper.  

\section*{Acknowledgments}

Ryan Deeley thanks: Sebastian Wieczorek, for helpful guidance during the initial framing of the metastability landscape and for illuminating instructive connections between multiple timescale and nonautonomous dynamics; Christopher K.R.T. Jones, for continued interest and useful suggestions for this project; Val\'{e}rian Jacques-Dumas, for assistance with the adaptive multilevel splitting algorithm; Paul Ritchie, for indicating some useful lines of research; and C\'{o}nall Kelly, for highlighting convergent numerical solvers for the stochastic FitzHugh-Nagumo system. This project has received funding from the European Union’s Horizon 2020 research and innovation programme under the Marie Sk{\l}odowska-Curie Innovative Training Network CriticalEarth, grant agreement no. 956170. The simulations were performed at the HPC Cluster ROSA, located at the University of Oldenburg (Germany) and funded by the DFG through its Major Research Instrumentation Programme (INST 184/225-1 FUGG) and the Ministry of Science and Culture (MWK) of the Lower Saxony State. For the purpose of open access, the author has applied a Creative Commons Attribution (CC BY) licence to any Author Accepted Manuscript version arising.

\section*{Author Declarations}

\subsection*{Conflict of Interest}

The authors have no conflicts to disclose.

\subsection*{Author Contributions}

\textbf{Ryan Deeley}: Conceptualisation (equal); Formal analysis (equal); Investigation (lead); Methodology (equal); Software (lead); Validation (equal); Visualisation (lead); Writing - original draft (lead); Writing - review \& editing (equal). \textbf{Peter Ashwin}: Conceptualisation (equal); Formal analysis (equal); Funding Acquisition (equal); Investigation (supporting); Methodology (equal); Supervision (equal); Validation (equal); Writing - review \& editing (equal). \textbf{Ulrike Feudel}: Conceptualisation (equal); Formal analysis (equal); Funding Acquisition (equal); Investigation (supporting); Methodology (equal); Supervision (equal); Validation (equal); Writing - review \& editing (equal). 

\section*{Data Availability}

The data that support the findings of this study are available in the supplementary material\cite{supplementary_material} of this article and on Zenodo, Ref. \cite{combined_rate_noise_multiscale_zenodo}; the code that generates both this data and the plots in this paper would become in the case of acceptance (is) openly accessible on GitHub at \url{https://github.com/ryandeeley/combined_rate_noise_multiscale}.

\newpage

\bibliography{combined_rate_noise_multiscale}

\end{document}